\documentclass{gtart_h}  

\def\ifplaintex{\expandafter\ifx\csname documentclass\endcsname\relax}

\ifplaintex 
\else
\fi

\expandafter\ifx\csname beginpicture\endcsname\relax
\expandafter\ifx\csname documentclass\endcsname\relax
\input pictex \else
\input prepictex \input pictex \input postpictex \fi\fi

\def\gt{{\mathsurround=0pt\it $\cal G\mskip-2mu$eometry \&\ 
$\cal T\!\!$opology}}        

\def\gtp{{\mathsurround=0pt\it $\cal G\mskip-2mu$eometry \&\ 
$\cal T\!\!$opology $\cal P\!$ublications}}  

\def\lognumber#1{\def\thelognumber{#1}}
\def\volumenumber#1{\def\thevolumenumber{#1}}
\def\papernumber#1{\def\thepapernumber{#1}}
\def\volumeyear#1{\def\thevolumeyear{#1}}

\def\pagenumbers#1#2{\def\startpage{#1}\def\finishpage{#2}}
\def\published#1{\def\publishdate{#1}}
\def\proposed#1{\def\theproposer{#1}}
\def\seconded#1{\def\theseconders{#1}}
\def\received#1{\def\receiveddate{#1}}
\def\revised#1{\def\reviseddate{#1}}
\def\accepted#1{\def\accepteddate{#1}}
\def\asciititle#1{\def\theasciititle{#1}}

\long\def\asciiabstract#1{\long\def\theasciiabstract{#1}}
\def\asciikeywords#1{\def\theasciikeywords{#1}}

\let\\\par\let\thelognumber\relax
\let\thevolumenumber\relax\let\thepapernumber\relax
\let\thevolumeyear\relax\let\thesamplenumber\relax\let\startpage\relax
\let\finishpage\relax\let\publishdate\relax\let\receiveddate\relax
\let\reviseddate\relax\let\accepteddate\relax\let\theasciititle\relax
\let\theasciiauthors\relax
\let\theasciiabstract\relax\let\theasciikeywords\relax
\let\theasciiemail\relax\let\theshortauthors\relax\let\theshorttitle\relax

\long\def\maketitlep{   

\count0=\startpage

\gt\hfill      
\beginpicture
\setcoordinatesystem units <0.33truein, 0.33truein> point at 2.2 0.9
\setplotsymbol ({$\cal G$})
\plotsymbolspacing=9truept
\circulararc 315 degrees from 0 1 center at 0 0
\setplotsymbol ({$\cal T$})
\circulararc 315 degrees from 1 -1 center at 1 0
\endpicture
\break
{\small\ifx\thesamplenumber\relax 
Volume \else Sample
\fi\thevolumenumber\ (\thevolumeyear)
\startpage--\finishpage\nl
Published: \publishdate}
\vglue 0.5truein plus 0.4fil minus 0.1truein

{\parskip=0pt\leftskip 0pt plus 1fil\def\\{\par\smallskip}{\ifplaintex\large
\else\Large\fi\bf\thetitle}\par\medskip}   

\vglue 0pt plus 0.1fil 

{\parskip=0pt\leftskip 0pt plus 1fil\def\\{\par}{\sc\theauthors}
\par\medskip}

\vglue 0pt plus 0.1fil 

{\small\parskip=0pt\let\newline\\
{\leftskip 0pt plus 1fil\def\\{\par}{\sl\theaddress}\par}
\expandafter\ifx\theemail\relax    
\relax\else\vglue 5pt plus 0.02fil minus 2pt\def\\{\stdspace{\rm 
and}\stdspace} 
\cl{Email:\stdspace\tt\theemail}\fi
\ifx\theurl\relax                  
\relax\else\vglue 5pt plus 0.02fil minus 2pt\def\\{\stdspace{\rm 
and}\stdspace}
\cl{URL:\stdspace\tt\theurl}\fi\par}

\vglue 7pt plus 0.3fil minus 3pt

{\bf Abstract}
\vglue 5pt plus 0.1fil minus 2pt

\theabstract

\vglue 7pt plus 0.3fil minus 3pt

{\bf AMS Classification numbers}\quad Primary:\quad \theprimaryclass

Secondary:\quad \thesecondaryclass

\vglue 5pt plus 0.3fil minus 2pt

{\bf Keywords:}\quad \thekeywords

\vglue 10pt plus 0.5fil minus 5pt

{\small  Proposed: \theproposer\hfill Received: \receiveddate\nl
Seconded: \theseconders\hfill 
\ifx\reviseddate\relax                         
Accepted: \accepteddate                        
\else
Revised: \reviseddate                          
\fi}
\eject
}       

\let\maketitlepage\maketitlep
\let\maketitle\maketitlepage

\font\phead=cmsl9 scaled 950
\font=cmsl9 scaled 1050
\font\pnum=cmbx10 scaled 913
\font\lnum=cmbx10 
\font\pfoot=cmsl9 scaled 950
\font=cmsl9 scaled 1050
\ifplaintex
\headline{\vbox to 0pt{\vskip -4.5mm\line{\small\phead\ifnum
\count0=\startpage ISSN 1364-0380 (on line)
1465-3060 (printed) \hfill {\pnum\folio}\else\ifodd\count0\def\\{ }%
\ifx\theshorttitle\relax\thetitle\else\theshorttitle\fi\hfill{\pnum\folio}
\else\def\\{ and }{\pnum\folio}\hfill\ifx\theshortauthors\relax\theauthors
\else\theshortauthors\fi\fi\fi}\vss}}
\footline{\vbox to 0pt{\vglue 0mm\line{\small\pfoot\ifnum\count0=\startpage
\copyright\ \gtp\hfill\else
\gt, Volume \thevolumenumber\ (\thevolumeyear)\hfill\fi}\vss
}}
\else
\makeatletter
\def\@oddhead{{\small\lhead\ifnum\count0=\startpage ISSN 1364-0380 (on line)
1465-3060 (printed) \hfill {\lnum\number\count0}\else\ifodd\count0
\def\\{ }\ifx\theshorttitle\relax \thetitle \else\theshorttitle\fi\hfill
{\lnum\number\count0}\else\def\\{ and }{\lnum\number\count0}
\hfill\ifx\theshortauthors\relax 
\theauthors\else\theshortauthors\fi\fi\fi}}\def\@evenhead{\@oddhead}
\def\@oddfoot{\small\lfoot\ifnum\count0=\startpage\copyright\ \gtp\hfill\else
\gt, Volume \thevolumenumber\ (\thevolumeyear)\hfill\fi}
\def\@evenfoot{\@oddfoot}
\makeatother
\fi

\newwrite\gtoutfile
\long\gdef\makeheadfile{  
{\def\\{, }\def\s{ }
\immediate\openout\gtoutfile head.xxx
\immediate\write\gtoutfile{Proxy-for: \ifx\theasciiauthors\relax
\theauthors\else\theasciiauthors\fi\s<\ifx\theasciiemail\relax\theemail\else\theasciiemail\fi>}
\immediate\write\gtoutfile{\noexpand\\}
\immediate\write\gtoutfile{Authors: \ifx\theasciiauthors\relax
\theauthors\else\theasciiauthors\fi}
{\def\\{ }\immediate\write\gtoutfile{Title: \ifx\theasciititle\relax
\thetitle\else\theasciititle\fi}}
\immediate\write\gtoutfile{Subj-class: GT or SG or MG etc}
\immediate\write\gtoutfile{MSC-class: \theprimaryclass\ifx\thesecondaryclass\relax\else, \thesecondaryclass\fi}
\immediate\write\gtoutfile{Journal-ref: Geom. Topol. \thevolumenumber
(\thevolumeyear) \startpage-\finishpage}
\immediate\write\gtoutfile{Comments: Published by Geometry and Topology at}
\immediate\write\gtoutfile{\s\s http://www.maths.warwick.ac.uk/gt/GTVol\thevolumenumber/paper\thepapernumber.abs.html}
\immediate\write\gtoutfile{\noexpand\\}
\immediate\write\gtoutfile{}
\ifx\theasciiabstract\relax
\immediate\write\gtoutfile{\theabstract}\else
\immediate\write\gtoutfile{\theasciiabstract}\fi
\immediate\write\gtoutfile{}
\immediate\write\gtoutfile{\noexpand\\}
\immediate\write\gtoutfile{}
\immediate\closeout\gtoutfile}}  

\def\maketitlepage{\maketitlep\makeheadfile}
\let\maketitle\maketitlepage

\lognumber{563}
\received{8 February 2005}
\volumenumber{9}\papernumber{35}\volumeyear{2005}
\pagenumbers{1539}{1601}   
\revised{2 July 2005}
\published{10 August 2005}
\accepted{4 August 2005}
\proposed{Benson Farb}
\seconded{David Gabai, Martin Bridson}

\usepackage{amssymb, graphicx, psfrag}

\def\fref#1{\hyperlink{#1anchor}{\ref*{#1}}}
\def\figref#1{\hyperlink{#1anchor}{Figure~\ref*{#1}}}
\def\anchor#1{\noindent\hypertarget{#1anchor}{\smash{$\phantom{99}$}}\newline}

\def\psfraga <#1,#2> #3#4{%
\psfrag {#3}{\smash{\rlap{\kern #1 \raise #2\hbox{#4}}}}}

\nocolon

\newtheorem{theorem}{Theorem}[section]

\newtheorem{lemma}[theorem]{Lemma}
\newtheorem{sublemma}[theorem]{Sub-Lemma}
\newtheorem{corollary}[theorem]{Corollary}

\theoremstyle{definition}

\def\Re{{\rm Re\/} }
\def\Im{{\rm Im\/} }

\def\B{\mbox{\mathsurround0pt\boldmath{$B$}}}%
\def\C{\mbox{\mathsurround0pt\boldmath{$C$}}}%
\def\F{\mbox{\mathsurround0pt\boldmath{$F$}}}%
\def\H{\mbox{\mathsurround0pt\boldmath{$H$}}}%
\def\P{\mbox{\mathsurround0pt\boldmath{$P$}}}%
\def\R{\mbox{\mathsurround0pt\boldmath{$R$}}}%
\def\Z{\mbox{\mathsurround0pt\boldmath{$Z$}}}%

\begin{document}

\title{A better proof of the Goldman--Parker conjecture}
\asciititle{A better proof of the Goldman-Parker conjecture}
\author{Richard Evan Schwartz}

\address{Department of Mathematics, University of Maryland\\Collage Park, 
MD 20742, USA}

\email{res at math dot brown dot edu}
\urladdr{http://www.math.brown.edu/~res/}

\begin{abstract}
The Goldman--Parker Conjecture
classifies the complex hyperbolic 
$\C$--reflection
ideal triangle groups up to discreteness.
We proved the Goldman--Parker Conjecture
in \cite{S0} using a rigorous
computer-assisted proof.  In this
paper we give a new and 
improved proof
of the Goldman--Parker Conjecture.
While the proof relies on the computer
for extensive guidance, the proof
itself is traditional.
\end{abstract}

\asciiabstract{%
The Goldman-Parker Conjecture classifies the complex hyperbolic
C-reflection ideal triangle groups up to discreteness.  We proved the
Goldman-Parker Conjecture in [Ann. of Math. 153 (2001) 533--598] using
a rigorous computer-assisted proof.  In this paper we give a new and
improved proof of the Goldman-Parker Conjecture.  While the proof
relies on the computer for extensive guidance, the proof itself is
traditional.}


\asciikeywords{Hyperbolic, complex reflection group, ideal triangle group,
Goldman-Parker conjecture}

\maketitlepage

\section{Introduction}\label{sec1}

Let $\H^2$ be the hyperbolic plane.
Let $G$ denote the usual reflection ideal
triangle group acting on $\H^2$.  The
standard generators of $G$ are
$\iota_0, \iota_1,\iota_2$.

$PU(2,1)$ is the holomorphic isometry
group of $\C\H^2$,
the {\it complex hyperbolic plane\/}.
See Section~\ref{sec2} for more details.
A {\it $\C$--reflection\/} is an order $2$ element of
$PU(2,1)$ which is conjugate to the element which
has the action
$(z,w) \to (z,-w)$.  
A {\it complex hyperbolic ideal triangle group representation\/}
is a representation of $G$ which maps
the generators to $\C$--reflections,
and the products of pairs of generators 
to parabolic elements.
Let ${\rm Rep\/}(G)$ denote the set of such
representations, modulo conjugacy.  It turns out that
${\rm Rep\/}(G)$ is a half-open interval,
naturally parametrized by $s \in [0,\infty)$.
See Section~\ref{sec2}.

Define
\begin{equation}
\underline s=\sqrt{105/3};
\hskip 30 pt
\overline s=\sqrt{125/3};
\end{equation}
In \cite{GP}, Goldman and Parker introduced
${\rm Rep\/}(G)$ (using different notation) and
proved that $\rho_s$ is a discrete embedding
if $s \in [0,\underline s]$.  They
conjectured that $\rho_s$ is a discrete embedding
iff $\rho_s(\iota_0\iota_1\iota_2)$ is
not an elliptic element of $PU(2,1)$.
This corresponds to parameters
$s \in [0,\overline s]$.  
We took care of the interval
$(\underline s,\overline s]$ in
\cite{S0}, using a rigorous
computer-assisted proof, together with
some new constructions in
complex hyperbolic geometry.  However,
the proof in \cite{S0} is extremely complicated
and requires massive computations.

The purpose of this
paper is to give a new and 
improved proof of the Goldman--Parker
Conjecture. Our new proof is based on an
idea we worked out, to a limited extent,
in \cite[Sections 8--10]{S1}.  To each
of the three generators
$I_{j,s}=\rho_s(\iota_j)$ we will associate
an piecewise analytic sphere
$\Sigma_{j,s}$.   We call $\Sigma_{j,s}$ a
{\it loxodromic $\R$--sphere\/}.  Our
construction is such that
$I_{j,s}(\Sigma_{j,s})=\Sigma_{j,s}$ and
that $I_{j,s}$ interchanges the two
components of $S^3-\Sigma_{j,s}$.
The key step in our argument is
showing that $\Sigma_{i,s} \cap \Sigma_{j,s}$
is a contractible set---the union of $2$ arcs
arranged in a `T' pattern---for $i \not = j$,
and that $\Sigma_{j,s}$ is embedded.  
This sets up a version of the familiar
ping-pong lemma, and it follows
readily from this picture that
$\rho_s$ is a discrete embedding.

In \cite[Sections 8--10]{S1} we established the
intersection and embedding properties of
our spheres for
all $s \in [\overline s-\epsilon,\overline s)$,
using a perturbative argument.  However,
we couldn't get
an effective estimate on
$\epsilon$ back then.
Here, in Section~\ref{sec3}, we develop a theory
for loxodromic $\R$--spheres and use it
to establish the two desired properties
for all
$s \in [\underline s,\overline s)$.
Pictures like \figref{fig4.2} indicate that our
construction works for all $s \in [0,\overline s)$.
However, there are certain technical details
we could not overcome when trying to deal
with parameters outside the range
$[\underline s,\overline s)$.

We wrote a Java applet which illustrates this paper
in great detail and, in particular, lets the
reader plot pictures like \figref{fig4.2} for all
parameter values.  The paper is independent
of the applet, but the applet greatly
enhances the paper because it lets the
reader see visually the objects we refer to
here mainly with symbols.
We encourage the reader to use the applet while
reading the paper.  One can access the applet
from my website.  The applet provides
massive hands-on evidence that our
construction works for all 
$s \in [0,\overline s]$.  In fact, most
of our proof works for all $s \in [0,\overline s]$
but there are certain technical estimates
we rely on that do not hold over the
whole range of parameters. 

Since I wrote \cite{S0} 
$7$ years ago, there has been considerable
development of complex hyperbolic
discrete groups.  Some of us feel that
all the new technology---eg, 
\cite{AGG}, \cite{S1},
\cite{FP}, \cite{S2}---should
reprove the Goldman--Parker Conjecture
without too much pain. 
Nonetheless, a new proof has never
appeared and I thought that this
paper would be of interest. Also,
I never liked my proof in \cite{S0}
and have wanted a better proof for a
long time.

This paper divides into $2$ halves.  The first
half is organized like this:
\begin{itemize}
\item Section~\ref{sec2}: background;
\item Section~\ref{sec3}: theory of loxodromic $\R$--spheres;
\item Section~\ref{sec4}: the proof.
\end{itemize}
The proof requires a handful of
technical estimates, which we make
in Sections~\ref{sec5}--\ref{sec7}.

The technical estimates all concern
the location in $S^3$ of a certain
collection of arcs of circles.
There is a $1$--parameter family of
these arcs and one can readily compute
their positions numerically.  You
can see from my applet (or from your
own experiments) that these estimates
hold by a wide margin and are
blatantly true for parameters
in $[\underline s,\overline s]$.  The
original version of this paper had 
computer-aided estimates on the 
locations of these arcs.  At the
request of the referee of this paper,
these computer-aided proofs have 
been replaced with analytic
calculations.  

The analytic calculations done in
the paper are in part based on
a brilliant algebraic idea due to
the anonymous\footnote{Eventually I guessed that the
referee was John Parker.  You can
tell the lion by his claw.}
referee.  The idea can be
summarized by saying that one
should introduce the parameter
$$x=\frac{e^2+|e|^2+\overline e^2}{1-|e|^2}$$
and then write all relevant
quantities in terms of $x$.
(See subsections \ref{subsec5.1}--\ref{subsec5.2} for details.)
Here $e$ (which is not to be confused
with the base of the natural log) 
is one of the coordinates
of an eigenvector of the word
$\rho(\iota_1\iota_0\iota_2)$.

I would like to thank Elisha Falbel,
Bill Goldman and John Parker
for many conversations, over the years, about
complex hyperbolic geometry.
Also, I would like to thank the
University of Maryland, the
Institute for Advanced Study,
the National Science Foundation (Grant DMS-0305047) 
and the John Simon Guggenheim Memorial Foundation,
for their generous support.

\section{Background}\label{sec2}

\subsection{Complex hyperbolic geometry}\label{subsec2.1}

\cite{E} and \cite{G} are good
references for complex hyperbolic geometry. 
\cite{S2} also has a good introduction.

\subsubsection{The ball model}\label{subsubsec2.1.1}

$\C^{2,1}$ is a copy of the vector space $\C^{3}$
equipped with the Hermitian
form
\begin{equation}
\label{h1}
\langle u,v\rangle=u_1 \overline v_1
+ u_2 \overline v_2-
         u_{3} \overline v_{3} \end{equation}
$\C\H^2$ and its ideal boundary are respectively the
projective images, in the complex projective plane $\C\P^2$, 
of
\begin{equation}
N_-=\{v \in \C^{2,1}|\ \langle v,v\rangle <0\}; \hskip 15 pt
N_0=\{v \in \C^{2,1}|\ \langle v,v\rangle=0\}
\end{equation}
(The set $N_+$ has a similar definition.)
The projectivization map
\begin{equation}
(v_1,v_2,v_3) \to (v_1/v_3,v_2/v_3)
\end{equation}
takes $N_-$ and $N_0$ respectively to the open 
unit ball and unit sphere in $\C^2$.
Henceforth we identify 
$\C\H^2$ with the open unit ball.
$\C\H^2$ is called the
{\it complex hyperbolic plane\/}.
It is a symmetric space of negative curvature.

\subsubsection{Slices}\label{subsubsec2.1.2}

There are two kinds of totally geodesic
$2$--planes in $\C\H^2$:
\begin{itemize}

\item The {\it $\R$--slices\/} are
$2$--planes, $PU(2,1)$--equivalent
to $\R\H^2=\R^2 \cap \C\H^2$. 

\item The {\it $\C$--slices\/} are
$2$--planes, $PU(2,1)$--equivalent 
to $\C\H^1=\C\H^2 \cap \C^1$.
\end{itemize}
Let $\F$ stand either for $\R$ or $\C$.
The accumulation set on $S^3$, of an
$\F$--slice, is called an $\F$--{\it circle\/}.
An $\F$--reflection is an involution
in ${\rm Isom\/}(\C\H^2)$ whose fixed
point set is an $\F$--circle.
The map $(z,w) \to (z,-w)$ is a prototypical
$\C$--reflection and the map
$(z,w) \to (\overline z,\overline w)$ is
a prototypical $\R$--reflection.
The $\F$--slice determines the $\F$--reflection
and conversely.

\subsubsection{Isometries}\label{subsubsec2.2.3}

$SU(2,1)$ is the  $\langle,\rangle$ preserving
subgroup of $SL_3(\C)$, the special complex linear group.
$PU(2,1)$ is the projectivization of $SU(2,1)$.
Elements of $PU(2,1)$ act isometrically on $\C\H^2$ and
are classified according to the usual scheme
for groups acting on negatively curved spaces.
Loxodromic elements move every point
of $\C\H^2$ greater than some $\epsilon>0$;
elliptic elements fix a point in $\C\H^2$;
and the remaining elements are parabolic.

We now discuss $\C$--reflections in more detail.
Given a vector $C \subset N_+$ we define
\begin{equation}
I_{C}( U)=- U+
\frac{2 \langle  U, C\rangle}{\langle  C,
 C\rangle}   C.
\label{reflection}
\end{equation}
$I_{C}$ is an involution fixing $C$ and
$I_{C} \in SU(2,1)$.  See \cite[page 70]{G}. 
The element of $PU(2,1)$ corresponding
to $I_C$ is a $\C$--reflection.
Every $\C$--reflection is conjugate to the
map $(z,w) \to (z,-w)$ discussed above.
$\C$--reflections are also called
{\it complex reflections\/}.

\subsection{Heisenberg space}\label{subsec2.2}

\subsubsection{Basic definitions}\label{subsubsec2.2.1}

In the ball model, 
$\C\H^2$ is a ball sitting
inside complex projective space $\C\P^2$.
For this discussion we fix some $p \in S^3$,
the ideal boundary of $\C\H^2$.
There exists a complex projective
automorphism $\beta$ of $\C\P^2$ which
maps $p$ to a point in 
$\C\P^2-\C^2$ and which
identifies
$\C\H^2$ with the {\it Siegel domain\/}:
\begin{equation}
\label{siegal}
Z=\{(z_1,z_2)|\ 2\Re(z_1)< -|z_2|^2\}  \subset \C^2 \subset \C\P^2
\end{equation}
We write $\infty=\beta(p)$ in this case.
The isometries of $\C\H^2$ which fix
$\infty$ act as complex linear automorphisms
of $Z$.   The set $\partial Z$ is characterized
as the set of null vectors relative to the
Hermitian form
\begin{equation}
\label{Siegelform}
\langle u,v \rangle' = u_1\overline v_3+u_2\overline v_2+
u_3 \overline v_1.
\end{equation}
We call $\cal H$$=\C \times \R$ {\it Heisenberg space\/}.
$\cal H$ is equipped with a group law:
\begin{equation}
(\zeta_1,t_1) \cdot 
(\zeta_2,t_2)=(\zeta_1+\zeta_2,t_1+t_2+2 \Im(\overline \zeta_1 \zeta_2))
\end{equation}
There is a natural map from $\partial Z$ to $\cal H$, given
by
\begin{equation}
\mu(z_1,z_2)=(\frac{z_2}{\sqrt 2},\Im(z_1)).
\end{equation}
The inverse map is given by
\begin{equation}
(z,t) \to (-|z|^2+it,z\sqrt 2).
\end{equation}
A {\it Heisenberg stereographic projection from $p$\/}
is a map $\B\co  S^3-\{p\} \to \cal H$
of the form $\mu \circ \beta$ where
$\beta$ is as above.
We write $\infty=\B(p)$ in this case.
We will somewhat abuse terminology and speak
of elements of $PU(2,1)$ acting on
$\cal H$. We mean that the conjugate
of an element, by Heisenberg stereographic
projection, acts on $\cal H$.
If such a map stabilizes $\infty$, it
acts as an affine map of $\cal H$.

\begin{itemize}
\item The $\C$--circles in $\cal H$
which contain $\infty$ all
have the form $(\{z\} \times \R) \cup \infty$.
The remaining $\C$--circles are ellipses which
project to circles in $\C$.  The plane
containing the ellipse is the {\it contact plane\/}
based at the center of mass of the ellipse.
See below for more detail.
\item The $\R$--circles which contain $\infty$ are
straight lines.  One of these $\R$--circles
is $(\R \times \{0\}) \cup \infty$.   
The bounded $\R$--circles in
$\cal H$ are such that their
projections to $\C$ are lemniscates.
\end{itemize}

\subsubsection{The contact distribution}\label{subsubsec2.2.2}

The set of complex lines tangent to $S^3$ forms
a $PU(2,1)$--invariant 
contact distribution on $S^3$. 
The $\R$--circles are tangent to this distribution
and the $\C$--circles are transverse to it.
The image
of the contact distribution, under Heisenberg
stereographic projection, is a contact
distrubition on $\cal H$.  It is defined
as the kernel of the $1$ form
$dt+2(xdy-ydx)$, when points in
$\cal H$ are written as $(x+iy,t)$.
Compare \cite[page~124]{G},
Any element of $PU(2,1)$ acting
on $\cal H$ respects this contact
distribution.  Each plane in the distribution
is called a {\it contact plane\/}.

\medskip
{\bf Area Principle}\qua
Suppose $\alpha$ is a piecewise smooth curve
in ${\cal H\/}$, tangent to the contact
distribution, such that
$\pi(\alpha)$ is a closed loop.
Then the height difference---meaning the
difference in the $t$--coordinates---between the
two endpoints of $\alpha$ is $4$ times
the signed area of the compact
region enclosed by $\pi(\alpha)$.
This is basically
Green's theorem.
Compare \cite[Section 4]{G}.
Call this principle the {\it area principle\/}.

\subsection{Spinal spheres}\label{subsec2.3}

\subsubsection{Basic definitions}\label{subsubsec2.3.1}

Basic information about bisectors and
spinal spheres can be found in \cite{G}.
Here we recall some of the basics.

A {\it bisector\/} is a subset
of $\C\H^2$ of the form
$\{x \in \C\H^2|\ d(x,p)=d(x,q)\}$.
Here $p \not = q$ are two distinct
points in $\C\H^2$ and $d$ is the
complex hyperbolic metric. A {\it spinal
sphere\/} is the ideal boundary of
a bisector.   Every two spinal spheres
are equivalent under $PU(2,1)$,
even though this is not immediately
obvious.
Equivalently, a spinal
sphere is any set of the form
$$\B^{-1}((\C \times \{0\}) \cup \infty).$$
Here $\B$ is a Heisenberg stereographic
projection.   Thus,
$$S=(\C \times \{0\}) \cup \infty$$ is a model
in $\cal H$ for a spinal sphere.  From
the second definition we see some of the
structure of spinal spheres.
Here are some objects associated to $S$:
\begin{itemize}
\item
$S$ has a singular
foliation by $\C$--circles.  The leaves
are given by $C_r \times \{0\}$ where
$C_r$ is a circle of radius $r$ centered
at the origin.  The singular points
are $0$ and $\infty$.  We call this
the $\C$--{\it foliation\/}.
\item $S$ has a singular foliation by
$\R$--circles.  The leaves are horizontal lines
through the origin.  The singular
points are again $0$ and $\infty$.
We call this the
$\R$--{\it foliation\/}.
\item The singular points $0$ and $\infty$
are called the
{\it poles\/} of $S$.
\item The {\it spine\/} of
$S$ is defined as the $\C$--circle
containing the poles.  In our case,
the spine is $(\{0\} \times \R) \cup \infty$.
Note that the spine of $S$
only intersects $S$ at the
singular points. 
\end{itemize}
Any other spinal
sphere inherits this structure, by symmetry.
The two foliations on a spinal sphere look like
lines of lattitude and longitude on
a globe.
A spinal sphere is uniquely determined by its
poles.  
Two spinal spheres are {\it cospinal\/}
if they have the same spine.  

\subsubsection{Generic tangencies with spinal spheres}\label{subsubsec2.3.2}

\begin{figure}[ht!]\anchor{fig2.1}
\cl{\small
\psfrag {A}{$A$}
\psfrag {B}{$B$}
\psfrag {L}{$L$}
\psfrag {p}{$p$}
\psfrag {q}{$q$}
\psfrag {r}{$r$}
\psfrag {g}{$\gamma$}
\psfrag {0}{$0$}
\includegraphics[width=3.6in]{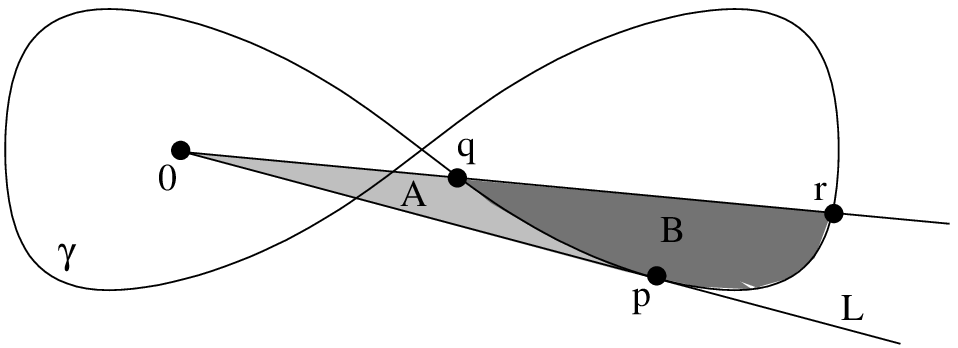}
}
\caption{}\label{fig2.1}
\end{figure}

In this section we prove a useful technical
result about how an $\R$--circle intersects a
spinal sphere.  Let $\pi\co  {\cal H\/} \to \C$
be projection.  The next result is illustrated
in \figref{fig2.1}.

\begin{lemma}
\label{Rintersect1}
Let $S=(\C \times \{0\}) \cup \infty$ as above.
Let $\gamma$ be 
a finite $\R$--circle in $\cal H$.  Suppose
\begin{itemize}
\item $\gamma$ is tangent to $S$ at
$p \not = 0$.
\item The line $L \subset \C$, 
containing $0$ and $\pi(p)$, is not
tangent to $\pi(\gamma)$ at
the double point of $\pi(\gamma)$.
\end{itemize}
Then
$\gamma$ has first but not
second order contact with $S$.
Moreover, a neighborhood of $p$ in
$\gamma$ lies on one side of $S$.
\end{lemma}

\proof
Since $\gamma$ is tangent to $S$ at $p$,
and $\gamma$ is also tangent to the contact
distribution, and the contact distribution
is not tangent to $S$ at $p$, we
see that $\gamma$ is tangent to the
$\R$--circle of $S$ which contains $p$.
This $\R$--circle is exactly $L \times \{0\}$.
But then $L$ is tangent to $\pi(\gamma)$
at $\pi(p)$.
\figref{fig2.1} shows the situation when 
a lobe of $\pi(\gamma)$ surrounds
$0$.  The other topological possibility
has the same proof.  The basic idea of
the proof is that $L$ does not have
second order contact with $\pi(\gamma)$
at $\pi(p)$, by convexity.

We first apply the area principle to the
integral curve $\alpha$ made from two
horizontal line segments and a portion
of $\gamma$, so that $\pi(\alpha)$
bounds the lightly shaded region $A$ shown
in \figref{fig2.1}.  From the area principle
we see that the height of $q$ is positive,
and also a quadratic function of the
Euclidean distance from $\pi(q)$ to
$\pi(p)$.  The quadratic dependence
comes from the strict convexity of
$\pi(\gamma)$ in a neighborhood of $\pi(p)$.

A similar argument works when we take
the relevant integral curve which
projects to the
region $A \cup B$. We see that the
height of the point $r$ is positive,
and also a quadratic function of the
Euclidean distance from $\pi(r)$ to
$\pi(p)$.   Now we know that points
on $\gamma$, on either side of
$p$, rise quadratically up and away from
$S$.
\endproof

\begin{corollary}
\label{contact}
Suppose $\gamma$ links the spine of a spinal
sphere $\Sigma$ and $\gamma$ is tangent
to $\Sigma$ at some point $p$. Then
$\gamma$ has first but not second order
contact with $\Sigma$ at $p$ and a
neighborhood of $p$ in $\gamma$ lies
to one side of $\Sigma$.
\end{corollary}

\proof
When we normalize so that $\Sigma=S$.
then the spine of $\Sigma$ projects to $0$.
One of the lobes of $\pi(\gamma)$ surrounds $0$,
and $p$ projects to some nonzero point.
In short, we have the picture in \figref{fig2.1},
and the hypotheses of the previous result
are forced.
\endproof

\subsection[Equations for C-circles]{Equations for $\C$--circles}\label{subsec2.4}

Suppose that $C$ is a $\C$--circle in
${\cal H\/}$ which links $\{0\} \times \R$.
Let $\pi\co  {\cal H\/} \to \C$ be projection as above.
Then $\pi(C)$ is a circle in $\C$ which
surrounds $0$.   As in \cite[section 2]{S1} we
study $\Psi_*(C)$, where $\Psi_*$ is the map
\begin{equation}
\Psi_*(z,t)=(\arg z,t).
\end{equation}
Define
\begin{equation}
r={\rm radius\/}(\pi(C)); \hskip 20 pt
d=|{\rm center\/}(\pi(C))|; \hskip 20 pt
A=(r/d)^2.
\end{equation}
We only define $A$ when $d>0$.
We call $A$ the {\it aspect\/} of $C$.
Note that $A>1$.

\begin{figure}[ht!]\anchor{fig2.2}
\cl{\small
\psfraga <-2pt,0pt> {p}{$\pi(C)$}
\psfraga <-2pt,0pt> {0}{$0$}
\psfraga <-2pt,0pt> {d}{$d$}
\psfraga <-2pt,0pt> {r}{$r$}
\psfraga <-2pt,0pt> {0}{$0$}
\includegraphics[width=1.5in]{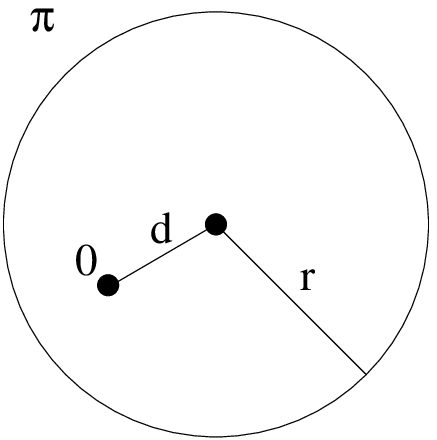}
}
\caption{}\label{fig2.2}
\end{figure} 

\begin{lemma}
Let $A$ be the aspect of $C$.
Up to scaling and rotation
$\Psi_*(C)$
is the graph of 
\begin{equation}
\label{calc1}
f_A(t)=\sin(t)(\cos(t)+\sqrt{A-\sin^2(t)})
\end{equation}
\end{lemma}

\proof
We normalize so that $(1,0)$ is the 
center of mass of $C$.  Then $d=1$ and
$C$ is contained in
the contact plane through $(1,0)$.  This plane
is spanned by $(1,0)$ and $(i,2)$.
Let $C_{\theta}$ be the point on $C$ such
that the line through $0$ and $\pi(C_{\theta})$
makes an angle of $\theta$ with the $x$ axis.
Then $\Psi_*(C)$ is the graph of the function
$\theta \to {\rm height\/}(C_{\theta})=2y$,
where $C_{\theta}=(x,y)$.
Our formula comes from solving the equations
$(x-1)^2+y^2=r^2$ and
$x=y\cot(\theta)$
in terms of $y$.
\endproof

\begin{lemma}
\label{calculus}
If $A \geq 9$ then $f_A''$ 
is negative on $(0,\pi)$ and
positive on $(\pi,2 \pi)$.
\end{lemma}

\proof
We compute that
\begin{equation}
\label{calc2}
\frac{d}{dA}f''_A(t)=
\frac{A(2-A+\cos(2t)) \sin(t)}{2 (A-\sin^2(t))^{5/2}}=
g_A(t)\sin(t),
\end{equation}
where $g_A(t)<0$.
Hence $\frac{d}{dA}f''_A$ is negative on $(0,\pi)$ and
positive on $(\pi,2 \pi)$.
We just need
to prove that
$f''_9$ is negative on $(0,\pi)$ and positive on
$(\pi,2\pi)$.
We compute that 
$f_9''(\pi/2)=-7/\sqrt 8<0$.
Thus, we just need to see that $f_9''(t)=0$ only
at $t=0$ and $t=\pi$.
Setting $u=\cos(t)$ we compute%
\footnote{We differentiate the function
$\sin(t)(\cos(t)+\sqrt{(A-1)+\cos^2(t)}),$
which is a re-writing of $f_A$, using Mathematica \cite{W}.}
\begin{equation}
f_9''=\frac{-4 h(u)\sin(t)}{(8+u^2)^{3/2}}; \hskip 15 pt
h(u)=
14+12u^2+u^4+8u\sqrt{8+u^2}+u^3\sqrt{8+u^2}.
\end{equation}
For $u \in (-1,1]$ we have
$$h(u) \geq 14+12u^2+u^4+24 u+3u^3=
(1+u)(14+10u+2u^2+u^3)>0.$$
This shows that $f_9''(t) \not = 0$ if $t \not \in\{0,\pi\}$.
\endproof

\noindent
{\bf Remark}\qua The preceding lemma is essentially
the same as Lemma 4.11 of \cite{FP}, with the
variable change $\phi=t+\pi/2$.

\subsection{Ideal triangle groups}\label{subsec2.5}
\label{idealdef}

\subsubsection{The basic definition}\label{subsubsec2.5.1}

We will use the same set-up as in \cite{S0}.
Given $s \in [0,\infty)$ we define
\begin{equation}
\label{beta}
\beta_s=\frac{s+i}{\sqrt{2+2s^2}}.
\end{equation}
Sometimes we write $\beta$ instead of
$\beta_s$, when the dependence is clear.
As we showed in \cite{S0},
every ideal triangle in $S^3$ is conjugate
to a triangle with vertices
\begin{equation}
\label{inittri}
p_0=(\beta,\overline \beta); \hskip 20 pt
p_1=(\beta,\beta); \hskip 20 pt
p_2=(\overline \beta,\overline\beta).
\end{equation}
In brief, the idea is that we can apply an element
of $PU(2,1)$ so that all three vertices of our ideal triangle
lie on the
Clifford torus $$\{(z,w)|\ |z|=|w|\} \subset S^3$$
and then we can rotate the Clifford torus until
the points are as above.

Let $I_j$ be the $\C$--reflection which
fixes $p_{j-1}$ and $p_{j+1}$.  We compute
that the elements $I_0, I_1,I_2$ are given by
\begin{equation}
\label{generators2}
 \left[\matrix{0 & -1 & 0\cr
              -1 &  0 & 0 \cr
               0 &  0 & -1}\right]; \hskip 10 pt
\left[\matrix{-1 & 0 & 0 \cr
               0 & 3 & -4\overline \beta \cr
               0 & 4\beta & -3} \right]; \hskip 10 pt
\left[\matrix{ 3 & 0 & -4\beta \cr
               0 &-1 & 0\cr
               4\overline \beta & 0 & -3}\right],
\end{equation}
Letting $g_0=I_1I_0I_2$ we compute
\begin{equation}
\label{trip}
\label{matrix}
g_0=\left[ \matrix{0&-1&0 \cr
-A_1&0&A_2 \cr
-A_2 &0&-\overline A_1}\right]; \hskip 15 pt
 A_1=\frac{s+17 i}{s+i}; \hskip 10 pt
A_2=\frac{12 \sqrt 2 i}{\sqrt{1+s^2}}.
\end{equation}
A direct computation, for example using the
result on \cite[page~201]{G}, shows that
$g$ is loxodromic for $s \in [0,\overline s)$
and parabolic for $\overline s$.  See \cite{GP}.

As we mentioned in the introduction, we are mainly
interested in the case when $s \in (\underline s,\overline s]$,
though many of our constructions work for
$s \in [0,\underline s]$ as well.   Actually,
we will carry out most of our constructions
for $s \in (\underline s,\overline s)$,
because there are now several good discreteness
proofs for the case $\overline s$.  See
\cite{S1} and \cite{FP}.

\subsubsection{Some associated objects}\label{subsubsec2.5.2}

When $s<\overline s$, the element $g_0$ is
loxodromic.  In this case $g_0$ stabilizes a pair
$(E_0,Q_0)$, where $E_0$ is a
$\C$--circle containing the fixed points
of $g_0$ and $Q_0$ is an arc of $E_0$
bounded by the two fixed points.
Of the two possible arcs, we choose
$Q_0$ so that it varies continuously
with the parameter and shrinks to
a point as $s \to \overline s$.
The curve 
\begin{equation}
\rho_0=\{(u\beta, \overline u \overline \beta)|\ u \in S^1\}
\end{equation}
is an $\R$--circle fixed by the map
$(z,w) = (\overline w,\overline z)$.  This
map interchanges $p_1$ and $p_2$ and fixes
$p_0$.  In short $\rho_0$ is an $\R$--circle
of symmetry for our configuration.

\begin{figure}[ht!]\anchor{fig2.3}
\cl{\small
\psfrag {p0}{$p_0$}
\psfraga <-1pt,0pt> {r0}{$\rho_0$}
\psfraga <-4pt,0pt> {C0}{$C_0$}
\psfraga <-2pt,0pt> {C1}{$C_1$}
\psfraga <-2pt,0pt> {C2}{$C_2$}
\includegraphics[width=1.8in]{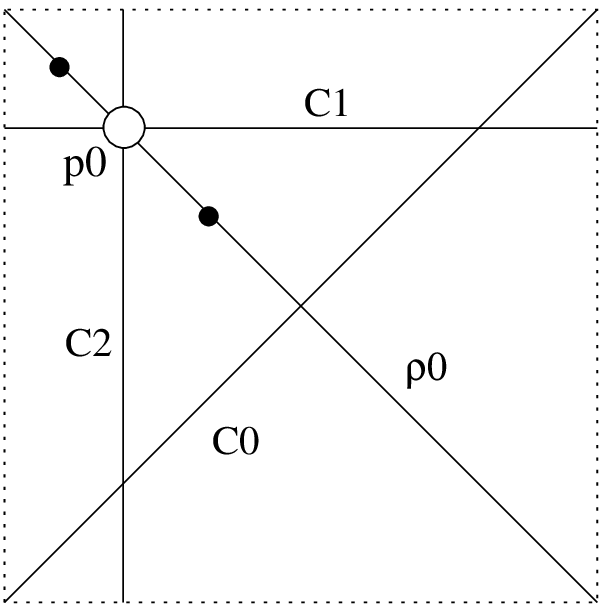}
}
\caption{}\label{fig2.3}
\end{figure}

One can define $Q_1$, $Q_2$, etc. by 
cycling the indices mod $3$.  
The objects 
\begin{equation}
(C_j,E_j,p_j,Q_j); \hskip 20 pt
j=0,1,2.
\end{equation}
are the elementary objects of interest to us.
\figref{fig2.3} shows those of the objects
which lie on the Clifford
torus, when the Clifford torus is drawn as
a square torus (in ``arg-arg coordinates'').
The black dots are the points of
$E_0 \cap \rho_0$.

\section[Loxodromic R-spheres]{Loxodromic $\R$--spheres}\label{sec3}
\setcounter{figure}{0}

\subsection{The cospinal foliation}\label{subsec3.1}

Our constructions are all
based on the pair $(E_0,Q_0)$ from subsection \ref{subsec2.5},
though we could take any pair $(E_0,Q_0)$ and
make the same definitions. Again,
$E_0$ is a $\C$--circle and $Q_0$ is a proper
arc of $E_0$.
Note that $Q_0$ determines
$E_0$ uniquely.  We include $E_0$ in
our notation for emphasis.

Let $p,q \in E_0$ be two distinct points.
The pair
$(p,q)$ is {\it harmonic\/} with respect to
$(E_0,Q_0)$ if the geodesic connecting
$p$ to $q$ in $\C\H^2$ is perpendicular
to the geodesic connecting the endpoints of
$Q_0$.  (Both these geodesics lie
in the $\C$--slice bounded by $E_0$.)
A spinal sphere $S$ is {\it harmonic\/}
w.r.t. $(E_0,Q_0)$ if the
poles of $S$ lie in $E_0$ and are harmonic
w.r.t $(E_0,Q_0)$.
An $\R$--arc $\alpha$ is
{\it harmonic\/} w.r.t. $(E_0,Q_0)$
if the endpoints of $\alpha$ are harmonic
w.r.t. $(E_0,Q_0)$.
Every harmonic $\R$--arc is contained in
a harmonic spinal sphere.

\begin{itemize}
\item Let $S(E_0,Q_0)$ denote the union of all
spinal spheres which are harmonic with
respect to $(E_0,Q_0)$.
We call $S(E_0,Q_0)$ the {\it cospinal foliation\/}.
\item Let $R(E_0,Q_0)$ denote the union of
all $\R$--arcs which are harmonic
with respect to $(E_0,Q_0)$.
\item Let $G(E_0,Q_0) \subset PU(2,1)$ denote
the group which fixes the endpoints of $Q_0$.
Then $G(E_0,Q_0)$ acts transitively 
the elements in $S(E_0,Q_0)$ and simply
transitively on the elements of $R(E_0,Q_0)$.
\end{itemize}

To see a picture
we work in ${\cal H\/}$ and normalize
so that $E_0=(\{0\} \times \R) \cup \infty$
and $Q_0$ is the unbounded arc whose
endpoints are $(0,\pm 1)$.  We call
this {\it standard position\/}. In this case
$(p,q)$ is harmonic with respect to
$(E_0,Q_0)$ iff $p=(0,r)$ and $q=(0,r^{-1})$.
We include the possibility that $r=0$, so
that $r^{-1}=\infty$.  
All the spinal spheres of interest to us
are bounded in $\cal H$, except for
$(\C \times \{0\}) \cup \infty$,
which corresponds to the case $r=0$.

\begin{lemma}
Every bounded spinal sphere in
$S(E_0,Q_0)$ is a convex 
surface of revolution.
\end{lemma}

\proof
All such spinal spheres are surfaces
of revolution, by symmetry.  Moreover,
all such spinal spheres are affine images
of the so-called {\it unit spinal sphere\/},
which has poles $(0,\pm 1)$.  The unit
spinal sphere satisfies the equation
$|z|^4+t^2=1$ and hence is convex.
See \cite[page~159]{G}.  Being affine
images of a convex set, the other spinal
spheres of interest to us are also convex.
\endproof

\begin{lemma}
\label{found1}
Every two distinct 
spinal spheres in $S(E_0,Q_0)$ are
disjoint.
\end{lemma}

\proof
Let $S_1$ and $S_2$ be two distinct spinal
spheres in $S(E_0,Q_0)$.  Using the
action of $G(E_0,Q_0)$ we can arrange
that $S_1=(\C \times \{0\}) \cup \infty$
and that the endpoints of $S_2$
are $(r,0)$ and $(r^{-1},0)$,
with $r \not = 0$.
But then $S_2$ lies either entirely
in the upper half space, or entirely the lower
half space.  In either case $S_2$ is 
disjoint from $S_1$.
\endproof

\subsection[Loxodromic R-cones and R-spheres]{Loxodromic $\R$--cones and $\R$--spheres}\label{subsec3.2}

Given $x \in S^3-E_0$ there is a unique element
$\alpha \in R(E_0,Q_0)$ such that
$x \in \alpha$.  We let
$\Sigma(E_0,Q_0;x)$ denote the portion
of $\alpha$ which connects $x$ to a point in
$Q_0$.  Given a subset $S \subset S^3-E_0$
we define
\begin{equation}
\Sigma(E_0,Q_0;S)=\bigcup_{x \in S} \Sigma(E_0,Q_0;x).
\end{equation}
We call $\Sigma(E_0,Q_0;S)$ the {\it loxodromic cone\/}
on $S$.

Let $C_1$ be a $\C$--circle which links
$E_0$.   Technical Lemma I (Lemma \ref{computation14}) establishes
this linking property when $(C_1,E_0)$ are as in subsection \ref{subsec2.5}.)
Let $(E_2,Q_2)=I_1(E_0,Q_0)$.
Here $I_1$ is the $\C$--reflection fixing
$C_1$. 
We say that a {\it loxodromic
$\R$--sphere\/} is an object of the form
\begin{equation}
\label{sig1}
\Sigma_1=\Sigma(E_0,Q_0;C_1) \cup
\Sigma(E_2,Q_2;C_1)
\end{equation}
If the $\C$--slice bounded by $C_1$ is
perpendicular to the $\C$--slice bounded
by $E_0$ then $C_1$ lies in one of
the elements of $S(E_0,Q_0)$ and
$\Sigma_1$ is a spinal sphere.
In general $\Sigma_1$ is not a spinal
sphere.  In Section~\ref{sec4} we will show that
the loxodromic $\R$--spheres of
interest to us are embedded spheres
but not spinal spheres.  

We are interested in
the case when $\Sigma_1$ is not a spinal
sphere.  Henceforth we assume that
$\Sigma_1$ is not a spinal sphere.
In this case we call $\Sigma_1$
{\it generic\/}.

\begin{lemma}
Let $\Sigma_1$ be a generic $\R$--sphere.
There exists a unique $\R$--circle $R_1$
such that $R_1$ intersects 
$E_0,C_1,E_2$ each in two points.
Also, $R_1 \subset \Sigma_1$ and
the $\R$--reflection in $R_1$ is
a symmetry of $\Sigma_1$.
We call $R_1$ the $\R$--axis
of $\Sigma_1$.
\end{lemma}

\proof
We put $(E_0,Q_0)$ in standard position.
Recall that $G(E_0,Q_0)$ is the stabilizer
subgroup of $Q_0$ which preserves the endpoints.
Using the action of $G(E_0,Q_0)$ we
can normalize so that the center of
mass of $C_1$ is $(r,0)$ for some $r>0$.
Consider the $\R$--circle
$R_1=(\R \times \{0\}) \cup \infty$.
By symmetry $R_1$ intersects
$C_1$ twice.  Also $R_1$ intersects
$E_0$ twice.  Finally,
we have $I_1(R_1)=R_1$ by symmetry.
Hence $R_1$ intersects $E_2$ twice.
The $\R$--reflection $J_1$ in
$R_1$ is an anti-holomorphic
element preserving $(E_0,Q_0)$ and $C_1$, and
hence a symmetry of
$\Sigma_1$.  
If there was some other axis $R'_1$ then
the composition of the $\R$--reflection
symmetries $J_1$ and $J_1'$ would
a non-trivial element of
$PU(2,1)$ preserving
both $(E_0.Q_0)$ and $C_1$.
But no such element exists.
Hence $R_1$ is unique.
\endproof
 
\subsection{The elevation map}\label{subsec3.3}

The set $R(E_0,Q_0)$ is topologically a
cylinder.  Since $G(E_0,Q_0)$ acts
transitively on  $R(E_0,Q_0)$, this
cylinder admits a natural family of
flat metrics.  Put another way, we can
write
$R(E_0,Q_0)=\R/2\pi\Z \times \R$.
The identification is unique up to
post-composition with a map
of the form $(x,y) \to (x+a,by+c)$.
That is, the identification is
unique up to {\it affine maps\/}.

There is a tautological map
$\Psi_0\co  S^3-E_0 \to R(E_0,Q_0)$ defined
as follows:  $\Psi_0(x)$ is the element
of $R(E_0,Q_0)$ which contains $x$.
If we identify $R(E_0,Q_0)$ with
$\R/2 \pi \times \R$ then we have
nice coordinates for this map.
Given $x \in S^3-E$ we let
$\widehat x$ be a lift of $x$. We
define
\begin{equation}
\label{elevation}
\Psi_0(x)=
\left(\arg \frac{\langle \widehat x,\widehat E_0 \rangle}
{\sqrt{\langle \widehat x,\widehat Q_1 \rangle
\langle \widehat x,\widehat Q_2 \rangle}},
\log \bigg|\frac{\langle \widehat x,\widehat Q_1 \rangle }
{\langle \widehat x,\widehat Q_2 \rangle }\bigg|\right)
\end{equation}
Here $\widehat Q_1$ and $\widehat Q_2$ are lifts
of the endpoints of $Q_0$ and $\widehat E_0$ is
a polar vector of $E_0$.   This is to say that
$\langle \widehat E_0,\widehat Y \rangle=0$ whenever
$\widehat Y$ is a lift of a point $Y \in E_0$.

\medskip
{\bf Remarks}\qua 

(i)\qua It is possible to choose a well-defined
branch of the square root in Equation \ref{elevation}.
This is basically a topological fact, deriving
from the fact that the map
$$f(x)\to \langle \widehat x,\widehat Q_1 \rangle
\langle \widehat x,\widehat Q_2 \rangle$$ induces
the map on homology
$$f_*\co  \Z=H_1(S^3-E_0) \to H_1(\C-\{0\})=\Z$$
which is multiplication by $2$. 

(ii)\qua Different choices of
lifts lead to maps which differ by
post-composition with affine maps.

(iii)\qua 
To see that Equation \ref{elevation} works as
claimed we compute that $\Psi_0$
conjugates
$G(E_0,Q_0)$ to isometries of
$\R/2 \pi \times \R$. The point here is
that $(\widehat Q_1,\widehat Q_2,\widehat E_0)$
is an eigenbasis for the elements of
$G(E_0,Q_0)$.

Henceforth we set $\Psi=\Psi_0$.
Note that
$\Psi(x)=\Psi(y)$ iff $x$ and $y$ belong to
the same element of $R(E_0,Q_0)$.  
$\Psi$ maps the
orbit $G(E_0,Q_0)(x)$ diffeomorphically
onto $\R/2 \pi\Z \times \R$.
hence $d\Psi$ has rank $2$ everywhere.

For any $x \in S^3-E_0$ let $\Pi_x$ denote the
contact plane at $x$.  Let 
\begin{equation}
L_x=d\Psi(\Pi_x)
\end{equation}
As we just mentioned, 
$d\Psi$ has full rank at at $x$ but
$d\Psi$ maps the $3$--dimensional vector space
$T_x(S^3)$ 
onto a $2$--dimensional tangent space.
The kernel of $d\Psi$ is the vector
tangent to the element of $R(E_0,Q_0)$
through $x$.  This kernel is therefore
contained in $\Pi_x$.  Hence
$L_x$ is a line.  The following
result captures some of the basic
features of this situation.

\begin{lemma}
\label{local}
Let $\gamma$ be an $\F$--circle and
let $x \in \gamma-E_0$ be a point.
\begin{itemize}
\item If $\F=\R$ and
$\gamma$ is
not tangent to $\Sigma(E_0,Q_0;x)$
then $\Psi(\gamma)$ is a nonsingular curve
at $\Psi(x)$ and the tangent line is
$L_x$. 
\item If $\F=\C$ then
$\Psi(\gamma)$ is
nonsingular at $\Psi(x)$ and
transverse to $L_x$.
\end{itemize}
\end{lemma}

\proof
When $\gamma$ is an $\R$--circle,
the tangent vector $v$ to $\gamma$ at
$x$ lies in $\Pi_x$ but is not contained
in the kernel of $d\Psi_x$.  Part
1 of our lemma follows from this fact.
When $\gamma$ is a $\C$--circle,
$v  \not \in \Pi_x$ and hence
$d\Psi_x(v) \not \in L_x$.
\endproof

\noindent
{\bf Remark}\qua
At this point, the reader anxious to see our
main construction should skip to Section~\ref{sec4}.

\subsection{More details on slopes}\label{subsec3.4}
\label{remote1}

From Lemma \ref{local} we see that
$L_x$ tells us a great
deal about what $\Psi$ does
to $\R$--circles and $\C$--circles.
We now investigate this further.
Let $\sigma_x$ denote the slope
of $L_x$.  Of course $\sigma_x$
depends on our choice of normalization,
but the general statements we make are
independent of normalization.
Let $S_0$
denote the spinal sphere whose
poles are $\partial Q_0$.

\begin{lemma}
\label{slope1}
If $x \in S_0$ then $L_x$ is a vertical line
and hence $\sigma_x$ is infinite.
Otherwise $\sigma_x$ is finite and nonzero.
\end{lemma}

\proof
Let $H \subset G(E_0,Q_0)$ denote the
$1$--parameter subgroup consisting of
the pure loxodromic elements.  These
elements do not twist at all in
the direction normal to the slice
bounded by $E_0$.  By symmetry
$\Psi$ maps the orbit $H(x)$ to
a vertical line in $\R/2 \pi\Z \times \R$.
On the other hand $H(x)$ is tangent
to $\Pi_x$ iff $x \in S_0$.  From
this we see that $\sigma_x$ is
infinite iff $x \in S_0$.
Now $\Psi$ maps the $\C$--circles
foliating the spinal spheres in
$S(E_0,Q_0)$ to horizontal lines.
From this fact, and from Lemma \ref{local},
we see that $\sigma_x \not = 0$.

\begin{lemma}
\label{slope3}
Let $x,y \in S^3-E_0$.
Then $\sigma_x=\sigma_y$ iff
$x$ and $y$ lie in the same $G(E_0,Q_0)$ orbit.
\end{lemma}

\begin{figure}[ht!]\anchor{fig3.1}
\cl{\small
\psfraga <-2pt,0pt> {hx}{$h(x)$}
\psfraga <-2pt,0pt> {hy}{$h(y)$}
\psfraga <-1pt,0pt> {x}{$x$}
\psfraga <-1pt,0pt> {y}{$y$}
\psfrag {0}{$0$}
\includegraphics[width=3.2in]{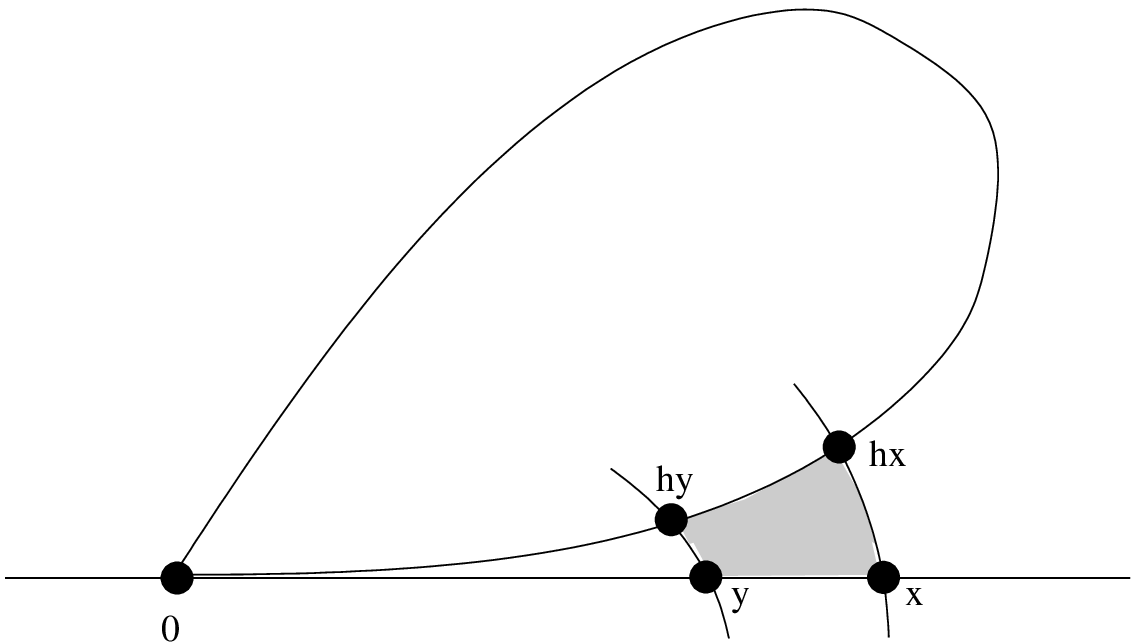}
}
\caption{}\label{fig3.1}
\end{figure}

\proof
By symmetry we have $\sigma_x=\sigma_y$ if
$x$ and $y$ are $G(E_0,Q_0)$ equivalent.
We just have to establish the converse.
Each $x \in S^3-E_0$ determines a 
$1$--parameter subgroup
$H_x \subset G(E_0,Q_0)$ which has
the property that the orbit 
$H_x(x)$ is integral to the contact
structure.  Then $\Psi$ maps
$H_x(x)$ to a geodesic on
$\R/2 \pi \Z \times \R$ which is
tangent to $L_x$.  It suffices
to show that $H_x \not = H_y$ if
$x$ and $y$ lie in different
$G(E_0,Q_0)$ orbits.   Suppose,
for the sake of contradiction,
that there are $G(E_0,Q_0)$--inequivalent
$x,y$ for which $H_x=H_y$.
Using the action of $G(E_0,Q_0)$
we can arrange that $\Psi(x)=\Psi(y)$.
Let $h \in H_x=H_y$.  By symmetry we have
$\Psi(h(x))=h(y)$.  But then
we can make a closed quadrilateral,
$\cal Q$ as follows: \begin{itemize}
\item One side of $\cal Q$ is
the portion of 
$H_x(x)$ which connects $x$ to $h(x)$.
\item One side of $\cal Q$ is the
portion of $H_y(y)$ which connects
$y$ to $h(y)$.  
\item One side of
$\cal Q$ is $\Sigma(E_0,Q_0;y)-\Sigma(E_0,Q_0;x)$.
\item One side of $\cal Q$ is
$\Sigma(E_0,Q_0;h(y))-\Sigma(E_0,Q_0;h(x))$.
\end{itemize}
Here we are choosing $x,y$ so that
$\Sigma(E_0,Q_0;x) \subset \Sigma(E_0,Q_0;y)$.
The shaded region in \figref{fig3.1} is the
projection of $\cal Q$ to $\C$.

We normalize so that $(E_0,Q_0)$ is the standard
pair and $x=(r,0)$ and $y=(s,0)$.  If we choose
$h$ close to the identity, then $\pi(\cal Q)$ projects
to an embedded quadrilateral in $\C$, as 
suggested by \figref{fig3.1}.  The point here
is that the fibers of $\Psi$ are lobes
of lemniscates which have their double points
at the origin.  Since $\cal Q$ is
integral to the contact structure, and yet
a closed loop, we
contradict the Area Principle of subsection \ref{subsubsec2.2.2}.
\endproof

\begin{corollary}[Slope Principle]
Let $\gamma_1,\gamma_2 \in S^3-E_0$ 
be two $\R$--arcs such
that $\Psi(\gamma_1)$ and $\Psi(\gamma_2)$
are nonsingular at a point
$x \in \R/2 \pi \Z \times \R$,
and tangent to each other at $x$.
Then $\gamma_1$ and $\gamma_2$
intersect at some point $y \in \Psi^{-1}(x)$.
\end{corollary}

\proof
Each fiber of $\Psi$,
including $\Psi^{-1}(x)$, intersects
each orbit of $G(E_0,Q_0)$ in one point.
Our result now follows
from Lemma \ref{local} and
Lemma \ref{slope3}.
\endproof

We normalize $\Psi$ so that $\sigma_x$ is
positive when $S_0$ separates $Q_0$
from $x$.  In this case we call $x$
{\it remote\/} from $Q_0$.
Thus we can say that $\sigma_x \in (0,\infty)$ iff
$x$ is remote from $S_0$.  

\subsection[Images of C-circles]{Images of $\C$--circles}\label{subsec3.5}

\begin{lemma}
\label{Cimage}
Suppose $C$ is a $\C$--circle which links
$E_0$.  Then $\Psi(C)$ is the graph of
a function $\psi\co  \R/2 \pi \Z \to \R$.
\end{lemma}

\proof
$\Psi(C)$ is a smooth 
loop by Lemma \ref{local}.
We need to prove that $\Psi(C)$ is
never vertical.
Let $C=C_1 \cup C_2$, where
$C_1$ is the closure of
the remote points of $C$ and $C_2$ is the
complement.
Let $I$ be a $\C$--reflection in a $\C$--circle
contained in $S_0$.  Then
$I$ interchanges the two components
of $S^3-S_0$, and $\Psi$ conjugates
$I$ to a reflection in a horizontal
line of $\R/2 \pi \Z \times \R$.
(This is seen by choosing
$\widehat Q_1$ and $\widehat Q_2$
in Equation \ref{elevation} so that these
vectors are swapped by $I$.)
Hence $\Psi(I(C_2))$ is the image
of $\Psi(C_2)$ reflected in a horizontal line.
By symmetry, then, it suffices to
show that $\Psi(C_1)$ is never vertical.

Let $x \in C_1$.
We normalize so that $(E_0,Q_0)$ is the standard
pair and
$x=(r,0)$ for some $r>0$.  Since $C$ links $E$,
we can parameterize $C$ as $C(\theta)=(z(\theta),t(\theta))$,
where $\theta$ is the angle between
the ray connecting $0$ to $\pi(C(\theta))$
and $\R$.  We are interested in $q=C(0)$.

If $t'(0)=0$ then $C$ is tangent 
to  $S=(\C \times \{0\}) \cup \infty$
at $q$.  But $S$ is a member of
$S(E_0,Q_0)$.  Hence $\Psi(C)$ is horozontal
at $\Psi(C(0))$.
Conversely, if $\Psi(C)$ is horizontal
at $\Psi(C(0))$ then
$C$ is tangent to $S$ and hence
$t'(0)=0$.

If $t'(0)>0$ then
$w=d\Psi_q(C'(0))$ 
lies in the interior of the
cone bounded by the $x$--axis and 
$L_{C(0)}$, a line whose slope is
either positive or infinite.
(The idea here is that the statement
holds when the center of mass of $C$ is
near $(0,0)$, and then Lemma \ref{local}
and the linking condition guarantee
that the cone condition holds no matter
how $C$ varies.) Hence $w$ is not vertical.

\begin{figure}[ht!]\anchor{fig3.2}
\cl{\small
\psfraga <0pt,7pt> {fiber}{fiber of $\Psi$}
\psfrag {C}{$C$}
\psfrag {C0}{$C(0)$}
\psfrag {Ce}{$C(e)$}
\psfrag {0}{$0$}
\psfraga <-2pt,0pt> {d}{$\delta$}
\psfrag {e}{$\epsilon$}
\includegraphics[width=4.5in]{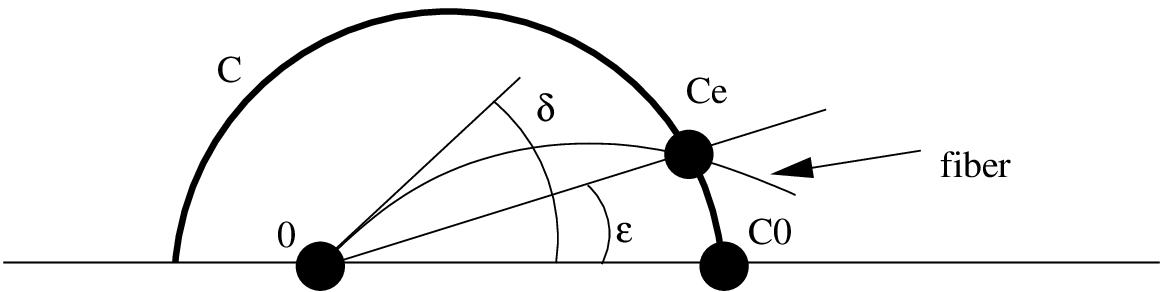}
}
\caption{}\label{fig3.2}
\end{figure}

If $t'(0)<0$ then the projection to $\C$ of the
fiber of $\Psi$ which contains $C(\epsilon)$
curves down and clockwise, as shown in \figref{fig3.2}. 
But then the horizontal component of
$w=d\Psi_q(\gamma'(0))$ exceeds $1$ because 
the angle $\delta$ in \figref{fig3.2} exceeds 
the angle $\epsilon$.
Again $w$ is not vertical. 
\endproof

We say that $C$ {\it generically links\/} $E_0$
if the $\C$--slice bounded by $C$ is not
perpendicular to the $\C$--slice bounded by
$E_0$.

\begin{lemma}
\label{onemax}
Suppose that $C$ is a $\C$--circle which 
generically links
$E_0$.  Then $\Psi(C)$ is the graph of
a function $\psi\co  S^1 \to \R$ which
has one maximum and one mimimum.
\end{lemma}

\proof
Since $C$ generically links $E_0$ the
image $\Psi(C)$ is not contained in
a horizontal line. The
$\Psi$--preimages of horizontal lines
are spinal spheres.  Since $C$ is
not contained in any of these
spinal spheres, $C$ can intersect
each of them at most twice. Hence
$\Psi(C)$ intersects each horizontal
line at most twice. Our result
follows immediately.
\endproof

\begin{corollary}
Suppose that $C_1$ generically links
$E_0$. Then
$\Sigma(E_0,Q_0;C_1)$ is an embedded topological
disk, analytic away from $C_1 \cup Q_0$.
\end{corollary}

\proof  
The set $\Sigma(E_0,Q_0;C_1)-Q_0$ is foliated
by $\R$--arcs of the form $\gamma(\theta)$,
where $\theta \in C_1$ is a point.
Two arcs $\gamma(\theta_1)$ and
$\gamma(\theta_2)$ are disjoint because
$\Psi(\gamma(\theta_j))=\Psi(\theta_j)$
and $\Psi(\theta_1) \not = \Psi(\theta_2)$.
Moreover, the arcs vary
analytically.  
Hence $\Sigma(E_0,Q_0;C_1)-Q_0$ is 
homeomorphic to an annulus (with
one boundary component deleted)
and analytic away from $C_1$.

\begin{figure}[ht!]\anchor{fig3.3}
\cl{
\includegraphics{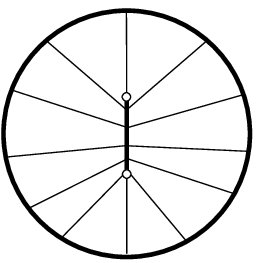}
}
\caption{}\label{fig3.3}
\end{figure}

There is a map $f\co C_1 \to Q_0$ given as follows: $f(x)$ is defined
to be the endpoint of $\Sigma(E_0,Q_0;x)$.  If $f(C_1)$ is more than
one point---as it is when $C$ generically links $E_0$---then $f$ is
generically $2$ to $1$ and $1$ to $1$ at exactly two points.  This
follows from Lemma \ref{onemax}.  Thus, $f$ has the effect of folding
$C_1$ in half over an arc of $Q_0$.  Hence the $\R$--arcs foliating
$\Sigma(E_0,Q_0;C_1)$ intersect $Q$ in pairs, with two exceptions.
Topologically, $\Sigma(E_0,Q_0;C_1)$ is obtained from an annulus by
gluing the inner circle together by a folding map, as in
\figref{fig3.3}.  From this description we see that
$\Sigma(E_0,Q_0;C_1)$ is an embedded disk, analytic off of $Q_0$.
\endproof

\subsection[Images of linked R-circles]{Images of linked $\R$--circles}\label{subsec3.6}

\begin{lemma}
\label{twotan}
Let $\gamma$ be an $\R$--circle which links
$E_0$.  Then $\gamma$ is tangent to exactly
two spinal spheres in $S(E_0,Q_0)$.
\end{lemma}

\proof
We normalize so that $(E_0,Q_0)$ is in standard
position.  Let $Q_t$ denote the unbounded
$\C$--arc whose endpoints are $(0,\pm (1+t))$.
Then $Q_t$ exits every compact subset of
$\cal H$ as $t \to \infty$.  Let $N(t)$ denote
the number of elements of $S(E_0,Q_t)$
which are tangent to $\gamma$.

We define $S(E_0,Q_{\infty})$ to be the
collection of spinal spheres of the
form $(\C \times \{s\}) \cup \infty$.
Given this definition we can define
$N(\infty)$ to be the number of
spinal spheres in $S(E_0,Q_{\infty})$
tangent to $\gamma$.   Let's analyze
$N(\infty)$ first.  Let
$\pi\co  {\cal H\/} \to \C$ be projection.
Suppose $\gamma$ is
tangent to a horizontal spinal sphere
$S_s=(\C \times \{s\}) \cup \infty$
at $x$.  Since $\gamma$ is
also tangent to the contact plane
$\Pi_x$, we see that $\gamma$ is
tangent to the line
$\Pi_x \cap S_s$.  But $\Pi_x \cap S_s$ is
a horizontal line which intersects
$E_0$.   Hence $\pi(\gamma)$ is tangent
at $\pi(x)$ to a line through the origin.
But $\pi(\gamma)$ is a lemniscate, one
of whose lobes surrounds the origin.
See \figref{fig2.1}.  Hence there are
only $2$ lines through the origin
which are tangent to $\pi(\gamma)$.
Hence $N(\infty)=2$.

Now fix some value of $t$. Since
$\gamma$ links $E_0$,
Corollary \ref{contact} applies:
If $\gamma$ is tangent to a
spinal sphere $S$ of
$\Sigma(E_0,Q_t)$ then
$\gamma$ locally lies to one
side of $S$ and has first but
not second order contact with
$S$.  Moreover,
as $t' \to t$, the spinal
spheres of $S(E_0,Q_{t'})$ converge
smoothly to the spinal spheres
of $S(E_0,Q_t)$.  These two properties
imply that the tangency points
vary continuously with $t$ and
cannot be created or destroyed
as $t$ changes.  The two
properties also hold at $t=\infty$.

From the discussion in the preceding
paragraph we see that $N(t)$ is
independent of $t$.  Since
$N(\infty)=2$ we also have
$N(1)=2$.
\endproof

\begin{corollary}[Elevation Image]
Suppose that $\gamma$ is
an $\R$--circle which links
$E_0$.  Then
$\Psi(\gamma)$ is the union of
two non-singular arcs, each having
nonzero slope at every point.
The two non-singular arcs
meet at two cusp points.
\end{corollary}

\proof
Everything but the statement about
the cusps follows from 
Lemma \ref{local}, Lemma \ref{slope1},
and Lemma \ref{twotan}.  The two
cusps appear because $\gamma$ locally lies
on one side of each of the two
spinal spheres to which it is
tangent. 
\endproof

\subsection{Linking of the poles}\label{subsec3.7}

As above
$I_1(E_0,Q_0)=(E_2,Q_2)$, where $I_1$ is the
$\C$--reflection fixing the $\C$--circle $C_1$
which generically links $E_0$.  Also
$R_1$ is the $\R$--axis of the $\R$--sphere $\Sigma_1$
given in Equation \ref{sig1}.

\begin{lemma}
\label{link5}
$E_0$ and $E_2$ are linked.
\end{lemma}

\proof
We normalize so that $(E_0,Q_0)$ is the standard pair
and the center of mass of $C_1$ is $(r,0)$ for some $r>0$.
Then $R_1=(\R \times \{0\}) \cup \infty$. 
Let $A_1 \subset R_1$ be the bounded interval
whose endpoints are $C_1 \cap R_1$.
The center of $A_1$ is $(r,0)$.
Since $C_1$ links $E_0$, we have
$0 \in A_1$.   Since the center of
$A_1$ is positive, $0$ lies in the left
half of $A_1$.  Now $I_1|_{R_1}$ 
acts as a linear fractional transformation
interchanging $A_1$ with $R_1-A_1$.
But then $I_1(0,0)=(s,0)$ with $s<0$ and
$I_1(\infty)=(r,0)$ with $r>0$.
The two points of $E_2 \cap R_1$ are
$I_1(0,0)$ and $I_1(\infty)$.  But
these points separate $(0,0)$ from
$\infty$ on $R_1$.  Hence
$E_2$ and $E_0$ are linked.
\endproof

\begin{figure}[ht!]\anchor{fig3.4}
\cl{\small
\psfrag {YE}{$\Psi(E_2)$}
\psfrag {YC}{$\Psi(C_1)$}
\psfrag {poss}{possible intersection}
\psfrag {points}{points}
\psfrag {perp}{perpendicular case}
\psfrag {int}{interlaced case}
\psfrag {Q}{\large ?}
\includegraphics[width=.94\hsize]{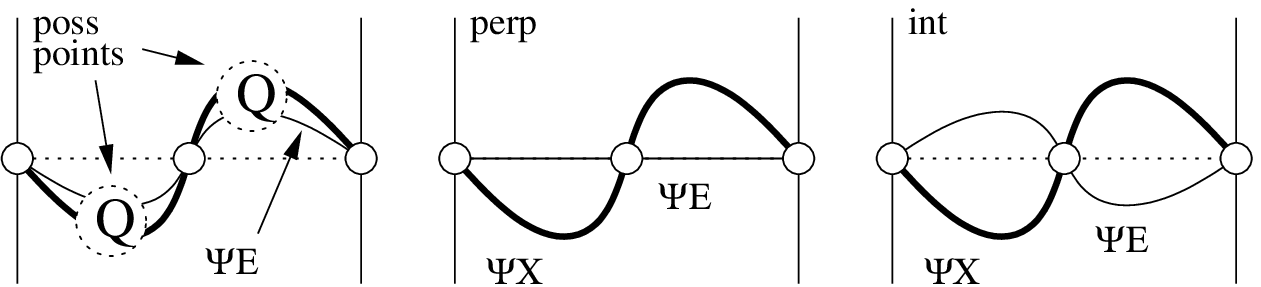}
}
\caption{}\label{fig3.4}
\end{figure}

$\Psi(C_1)$ and $\Psi(E_2)$
have some symmetry:  Let
$J_1$ be the $\R$--reflection fixing
$R_1$.  Then $\Psi$ conjugates $J_1$
to an isometric $180$ degree rotation of
$\R/2 \pi \Z \times \R$. The fixed
point set of this rotation is 
exactly $\Psi(R_1-E_0)$, which is a
pair of points on the same horizontal level
and $\pi$ units around from each other.
This rotation is a symmetry
of $\Psi(C_1)$ and also of $\Psi(E_2)$.
\figref{fig3.4} shows a
picture of the three possibilities.

We say that $\Sigma_1$ is {\it interlaced\/}
if the picture looks like the right hand side.
That is, a vertical line separates the minimum
of $\Psi(E_2)$ from the minimum of $\Psi(C_1)$.
When we normalize the interlaced case as above,
$(0,0)$ separates the center of mass of $C_1$
from the center of mass of $E_2$.  From this
we see that
the set of interlaced $\R$--spheres is connected.
In the interlaced case, the interlacing pattern
of the extrema of $\Psi(C_1)$ and $\Psi(E_2)$ forces
\begin{equation}
\label{interlace}
\Psi(C_1) \cap \Psi(E_2)=\Psi(R_1-E_0).
\end{equation}

\subsection{The Two Cusp Lemma}\label{subsec3.8}

We say that an $\R$--circle $\widehat \gamma$ is
{\it affiliated\/} with 
$\Sigma(E_0,Q_0;C_1)$ if
$\widehat \gamma$ contains an $\R$--arc
of the form $\Sigma(E_0,Q_0;x)$,
for $x \in C_1$, but $\widehat \gamma$ is
not the $\R$--axis $R_1$ of $\Sigma_1$.
The purpose of this section is to prove
the following result.

\begin{lemma}[Two Cusp]
\label{twocusp}
Suppose that $\Sigma_1$ is interlaced.
Suppose that $\widehat \gamma$ is an $\R$--circle
affiliated to $\Sigma(E_2,Q_2;C_1)$.  Then
$\Psi(\widehat \gamma)$. is the union of two
nonsingular arcs, each of which has everywhere
nonzero slope.  The two arcs are joined at
two cusps.
\end{lemma}

The Two Cusp Lemma is an immediate consequence of
Lemma \ref{linked} below and the
Elevation Image Lemma.

\begin{lemma}
\label{linked}
Suppose that $\Sigma_1$ is interlaced.
Then every $\R$--circle affiliated to
$\Sigma(E_0,q_0;C_1)$ links $E_2$.
\end{lemma}

\proof
By symmetry, every
$\R$--circle affiliated to
$\Sigma(E_0,q_0;C_1)$ links $E_2$ if and only if
every $\R$--circle affiliated to 
$\Sigma(E_2,Q_2;C_1)$ links $E_0$.

Let $\widehat \gamma$ be an $\R$--circle affiliated
to $\Sigma(E_2,Q_2;C_1)$.  We claim that
that $\widehat \gamma \cap E_0=\emptyset$.
Once we know this, we see that either
all affiliates of $\Sigma(E_2,Q_2;C_1)$
link $E_0$ or all affiliates fail to link
$E_0$. By continuity, the link/unlink
option is independent of the choice of
interlaced $\R$--sphere.  We check explicitly,
for one interlaced $\R$--sphere---eg,
the one in \figref{fig4.5}---that the
link option holds for some of the
affiliates.  Hence the link option
always holds.

It remains to establish our claim.
By symmetry $\widehat \gamma'=I_{C_1}(\gamma)$ is
affiliated with $\Sigma(E_0,Q_0;C_1)$. 
By construction, $\Psi(\gamma'-E_0)$
is a single point of $\Psi(C_1)$ and
$\Psi(\gamma'-E_0) \in \Psi(R_1-E_0)$
iff $R'=R_1$.  The point is that
$\Psi$ is injective on $C_1$.
Therefore 
\begin{equation}
\label{linked1}
 \Psi(\gamma'-E_0) \cap \Psi(R_1-E_0)=\emptyset.
\end{equation}
We have by hypotheses and Equation \ref{interlace} that
\begin{equation}
\label{linked2}
\Psi(\gamma'-E_0) \cap \Psi(E_2) \subset
\Psi(C_1) \cap \Psi(E_2)=
\Psi(R_1-E_0).
\end{equation}
Combining Equations \ref{linked1} and
\ref{linked2} we have
\begin{equation}
\Psi(\gamma'-E_0) \cap \Psi(E_2)=\emptyset.
\end{equation}
Hence $(\gamma'-E_0) \cap E_2=\emptyset$.
Since $E_0 \cap E_2=\emptyset$ we
conclude that $\widehat \gamma' \cap E_2=\emptyset$.
Hence $\widehat \gamma \cap E_0=\emptyset$.
This establishes our claim.
\endproof

\subsection{Asymmetry}\label{subsec3.9}

Let $\Delta_0$ denote the $\C$--slice which bounds
$E_0$.  Let $S_0$ denote the spinal sphere
whose poles are $Q_0$.  Let $R_1$ denote
the $\R$--axis of the $\R$--sphere $\Sigma_1$
given in Equation \ref{sig1}.  Again
recall that $I_1(E_0,Q_0)=(E_2,Q_0)$.  We
we will assume explicitly that 
$E_1$ and $E_2$ are generically linked.
Hence $E_i$ and $E_j$ are also generically linked.

\begin{figure}[ht!]\anchor{fig3.5}
\cl{\small
\psfrag {h7}{$\eta(Q_0)$}
\psfrag {h6}{$\eta(x)$}
\psfrag {h5}{$\eta(y)$}
\psfrag {h7}{$\eta(Q_0)$}
\psfrag {h4}{$\eta(E_2)$}
\psfrag {h3}{$\eta(Q_2)$}
\psfrag {Q0}{$\Theta_0$}
\psfrag {Q1}{$\Theta_1$}
\psfrag {Q2}{$\Theta_2$}
\psfraga <-2pt,0pt> {g}{$\gamma$}
\psfraga <-2pt,0pt> {g1}{$\gamma_1$}
\psfraga <0pt,2pt> {g0}{$\gamma_0$}
\includegraphics[width=2.5in]{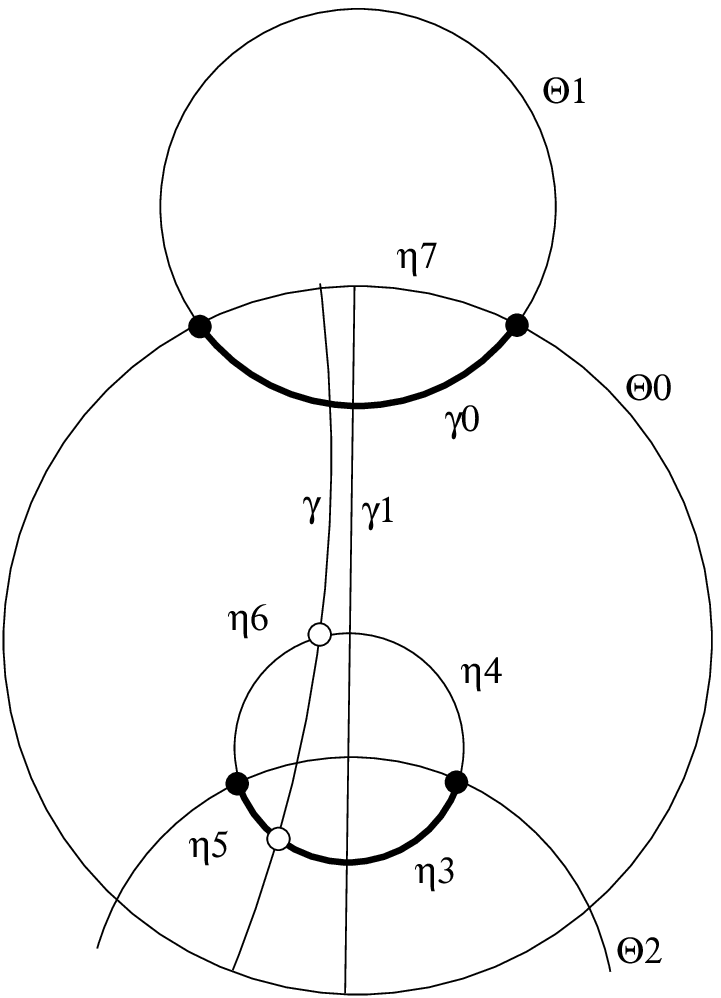}
}
\caption{}\label{fig3.5}
\end{figure}

Let $\eta\co  S^3 \to \Delta_0$ denote
orthogonal projection. 
The generic linking condition implies
that $\eta(E_j)$ is a circle (rather than
a point) for $j=1,2$.
Let $\Theta_1=\eta(E_0)$.  
For $j=0,2$ let $\Theta_j$ be the
circle which is perpendicular to
$\eta(E_j)$ and contains the
endpoints of $\eta(Q_j)$.  
Note that $\Theta_0$ and $\Theta_1$
intersect at right angles, since
$\Theta_0=\eta(S_0)$ and
$\Theta_1=\eta(E_0)$.Some of
these objects are drawn in \figref{fig3.5}.
Say that
$\Sigma_1$ is {\it asymmetric\/}
if $\Theta_0 \cap \Theta_1 \cap \Theta_2=\emptyset$.
This is the generic case.  The goal of this
section is to prove:

\begin{lemma}[Asymmetry Lemma]
\label{rise}
Suppose $\Sigma_1$ is asymmetric
and interlaced.
Let $x,y \in E_2$ be two points which are
harmonic with respect to $Q_2$.  Then
$\Psi(x)$ and $\Psi(y)$ lie on the same
horizontal line in $\R/2 \pi \Z \times \R$
iff $x,y \in R_1$.
\end{lemma}

\proof
We first list some basic properties of the
map $\eta$.

\begin{itemize}
\item $\eta(E_0)$ is the boundary of $\Delta_0$ and
$\eta(Q_0)$ is an arc of $\eta(E_0)$.
This follows from the fact
that $\eta$ is the identity on $E_0$.

\item If $C$ is a $\C$--circle which is
disjoint from $E_0$ then $\eta(C)$ is a
circle contained in the interior of $E_0$.
The restriction of $\eta$ to $C$ is
a linear fractional transformation.
This property comes from the fact
that $\eta$ is holomorphic on complex lines.

\item If $S$ is a spinal sphere whose spine
is $E_0$, then $\eta(S)$ is a geodesic in $\Delta_0$.
In particular, $\gamma_0=\eta(S_0)$ is the geodesic whose
endpoints are $\eta(Q_0)$.  Indeed, an
alternate definition of a spinal
sphere is the preimage of such a geodesic
under $\eta$.  See \cite{G}.

\item $\eta$ maps each spinal sphere in the
cospinal foliation to
geodesics perpendicular to $\eta(S_0)$.
This follows from symmetry:
namely that $\eta$ conjugates $G(E_0.Q_0)$ to
isometries of $\Delta_0$ which fix both
endpoints of $\eta(Q_0)$.

\item $\eta$ maps the $\R$--axis of
$\Sigma_1$ to a geodesic $\gamma_1$ which is
simultaneously perpendicular to
$\eta(S_0)$ and $\eta(C_1)$ and $\eta(Q_2)$.
Again this follows from symmetry:
The $\R$--reflection in the $\R$--axis
of $\Sigma_1$ preserves both $C_1$ and $Q_2$.
Indeed, the isometric reflection in
$\gamma_1$ stabilizes $\eta(Q_0)$ and
$\eta(Q_2)$ and $\eta(C_1)$ and
$\Delta_0$.
\end{itemize}

\noindent
{\bf Remark}\qua The reader can see all these
properties in action using my Applet.

Now we turn to the main argument in the proof
of the Asymmetry Lemma.
If $x,y \in R_1$ then $\Psi(x)$ and $\Psi(y)$
are precisely the two symmetry points of
$\Psi(E_2)$ and $\Psi(C_1)$ discussed in 
Lemma \ref{linked}, and then $\Psi(x)$ and
$\Psi(y)$ lie on the same horizontal line.

Suppose, conversely, that $\Psi(x)$ and
$\Psi(y)$ lie on the same horizontal line.
This means that $\eta(x)$ and $\eta(y)$
lie on the same geodesic $\gamma$ of $\Delta_0$,
where $\gamma$ is perpendicular to both
$\Theta_0$ and $\Theta_1$.
By assumption $x,y \in E_2$ are harmonic
with respect to $Q_2$.  
Since the restriction of $\eta$ to
$E_2$ is a linear fractional transformation,
the points $\eta(x)$ and $\eta(y)$
are in harmonic position with respect
to $\eta(Q_2)$.
So, we conclude that the geodesic $\gamma$
has the following properties:
\begin{itemize}
\item $\gamma$ is perpendicular to
$\Theta_0$ since $\gamma$ is a geodesic.
\item $\gamma$ is perpendicular to
$\gamma_0$, the geodesic connecting
the endpoints of $\gamma(Q_0)$.  Hence
$\gamma$ is perpendicular to $\Theta_1$.
\item $\gamma$ intersects
$\eta(E_2)$ in two points which are
in harmonic position with respect to
$\eta(Q_2)$.  But then $\gamma$ is
perpendicular to $\Theta_2$, the
circle which is perpendicular to 
$\eta(E_2)$ and contains $\eta(Q_2)$
in its endpoints.
\end{itemize}
\figref{fig3.5} shows a picture.
The lemma below says that there is
only one geodesic which has this property,
and this geodesic is
$\gamma_1=\eta(R_1)$.   Hence
$\Psi(x)$ and $\Psi(y)$ lie on
the same horizontal line as
$\Psi(R_1-E_0)$. This forces
$x,y \in R_1$. 
\endproof

\begin{lemma}
\label{circle}
Let $\Theta_0,\Theta_1,\Theta_2$ be $3$ circles
in $\C \cup \infty$.  Suppose that
$\Theta_0 \cap \Theta_1$ is a pair of points
and $\Theta_0 \cap \Theta_1 \cap \Theta_2=\emptyset$.
Then there is at most one circle which is
simultaneously perpendicular to
$\Theta_j$ for $j=0,1,2$.
\end{lemma}

\proof
We normalize by a Moebius transformation
so that $\Theta_0$ and $\Theta_1$ are lines
through the origin.  Then a circle in $\C$
is perpendicular to $\Theta_0$ and $\Theta_1$ iff
this circle is centered at the origin.
By assumption,
$\Theta_2$ is a finite circle in $\C$ which does
not contain the origin.  From here it is easy
to see that at most one circle, centered at the
origin, can be perpendicular to $\Theta_2$.
\endproof

We end this chapter with a result which
relates symmetry and asymmetry to the
image of the arc $\Psi(Q_2)$.

\begin{lemma}
\label{comp3}
$\Sigma_1$ is symmetric if and only if
the maximum and minimum heights of
$\Psi(E_2)$ occur at the endpoints
of $\Psi(Q_2)$.
\end{lemma}

\proof
If $\Sigma_1$ is symmetric then the geodesics
in $\Delta$ perpendicular to $\Theta_0$ and $\Theta_1$
and containing the endpoints of $\eta(Q_2)$ are
tangent to $\eta(E_2)$.  Recalling the
$4$th property of $\eta$ mentioned above,
the tangency property translates exactly
into the statement that the height of 
$\Psi|_{E_2}$ takes on its maxima and minima at
the endpoints of $Q_2$.  The converse is
proved simply by running the argument in
reverse.
\endproof

\subsection{Remoteness}\label{subsec3.10}

We say that $\Sigma_1$ is {\it remote\/} if
every point of
$\Sigma(E_2,Q_2,C_1)$ is remote from
$Q_0$.  In this case $\alpha=\Psi(\gamma)$
has everywhere positive slope, when
$\gamma$ is an $\R$--arc of $\Sigma(E_2,Q_2;C_1)$.
In this section we give a technical
criterion for remoteness.  We work in $\cal H$.

\begin{lemma}
\label{remote2}
$\Sigma_1$ is remote provided:
\begin{itemize}
\item $C_1$ has aspect at least $9$.
\item The endpoints of $Q_0$ have
the form $(0,\pm u)$, with $u \geq 4$.
\item $\pi(C_1)$ has radius $1$.
\item $C_1$ is centered on the real axis.
\end{itemize}
\end{lemma}

\proof
Given our bound on the aspect, the center of
$C_1$ is at most $1/3$ from $0$.  Hence
every point of $C_1$ is at most $2/3$ from
$\C \times \{0\}$.
Also, there is a
spinal sphere $S_1$, containing $C_1$ such that
$$S_1 \subset \pi(C_1) \times [-5/3,5/3].$$
To obtain $S_1$ we simply take the spinal
sphere with poles $(0,\pm 1)$ and left
translate by less than $1/3$ along 
$\R \times \{0\}$.

The poles of $S_0$ are
$(0,\pm u)$, with $u>4$.
It is easy to that the bounded
portion of this huge (and convex) set
contains $S_1$
in its interior.   
Let $\gamma$ be an $\R$--arc of $\Sigma(E_0,Q_0;C_1)$.
Then $I_1(\gamma)$ is an 
$\R$--arc of $\Sigma(E_2,Q_2;C_1)$.
We claim that $\gamma$ intersects $S_1$ only
at its endpoint.  In this case $I(\gamma)$
is contained in the bounded portion of $S_1=I_1(S_1)$,
which is in turn contained in the bounded
portion of $S_0$.  Hence $I(\gamma)$ is remote.

To finish our proof we need to establish our
claim. Let $\widehat \gamma$ be the fiber of
$\Psi$ containing $\gamma$.  Then
$\widehat \gamma$ is harmonic with respect to $(E_0,Q_0)$
and $\pi(\widehat \gamma)$ is one lobe of a lemniscate.
Let $(0,t_1)$ and $(0,t_2)$ be the two endpoints
of $\widehat \gamma$, with $t_1<t_2$.
Without loss of generality assume that $t_1>0$.
Then $\widehat \gamma$ rises up from its lower
endpoint until it intersects $C_1$. Hence
$t_1 \in [0,2/3]$.
Also $t_1t_2=u^2>14$.  Hence
$t_2>21$. Hence $\widehat \gamma$ rises up
at least $20$.  The projection $\pi(\widehat \gamma)$
is a huge lemniscate. From all this information
we can see that $\gamma$ only intersects
$S_1$ at its endpoint:  The only points
$x \in \widehat \gamma$ with $\pi(x) \in \pi(C_1)$
have height greater than $5/3$.
\endproof

\section{The proof}\label{sec4}
\setcounter{figure}{0}

\subsection{Main construction}\label{subsec4.1}

Let $s \in [\underline s,\overline s)$.
Let $(C_j,E_j,Q_j,R_j,p_j)$ be as in subsection \ref{subsec2.5}.
(Actually, we defined $R_j$ in Section~\ref{sec3}, as the
$\R$--axis of $\Sigma_j$.)
Also define
\begin{equation}
Q_{21}=Q_2 \cap \Sigma_1; \hskip 30 pt
Q_{12}=Q_1 \cap \Sigma_2.
\end{equation}
All our objects depend on a parameter $s$,
though we typically suppress $s$ from our notation.

Our proof includes several technical lemmas
whose proofs will be given in Sections \ref{sec5}--\ref{sec7}.
It is to be understood that
these results are only proved for
parameters in $[\underline s,\overline s]$.

\begin{lemma}[Technical Lemma I]
\label{computation14}
The following is true:
\begin{enumerate}
\item
$E_0$ and $C_1$ are linked.
\item
$E_0$ is normalized to be
$(\{0\times \R\}) \cup \infty$ in Heisenberg space
then the aspect $A$ of $C_1$ is at least $9$.
\item  Suppose we normalize so that
$C_1 \cap C_2=(1,0)$ in $\cal H$ and
$E_0=\{0\} \times \R$ and the map
$(z,t) \to (\overline z,-t)$ swaps
$C_1$ and $C_2$.  Then
$0$ is closer to $1=\pi(p_0)$ than it is
to the other intersection point of
$\pi(C_1)$ and $\pi(C_2)$. 
\end{enumerate}
\end{lemma}

\begin{lemma}[Technical Lemma II]
The following is true:
\begin{enumerate}
\item $\Sigma_1$ satisfies the
criteria of Lemma \ref{remote2}.
Hence $\Sigma_1$ and $\Sigma_2$ are remote.
\item The curve $\Psi(Q_2)$ has 
negative slope, even at the endpoints.
\end{enumerate}
\end{lemma}

\begin{lemma}[Technical Lemma III]
The following is true:
\begin{enumerate}
\item
A horizontal line
in $\R/2 \pi \Z \times \R$
separates $\Psi(Q_{21})$ from $\Psi(C_2)$,
with $\Psi(Q_{21})$ lying on top.
\item If we normalize as in item $3$ of
Technical Lemma I then the center of
$C_1$ lies above all points of $C_2$.
Likewise the center of $C_2$ lies below
all points of $C_1$.
\end{enumerate}
\end{lemma}

Now we are ready for our main construction.
By symmetry $C_i$ and $E_j$ are linked for $i \not = j$.
We define
\begin{equation}
\Sigma_{ij}=\Sigma(E_j,Q_j;C_i).
\end{equation}
Then $\Sigma_i=\Sigma_{ij} \cup \Sigma_{ik}$ is an
$\R$--sphere because $I_i$ interchanges
$(E_j,Q_j)$ and $(E_k,Q_k)$.  Here
$i,j,k$ are distinct indices.
Let $\Psi=\Psi_0$ be the
map from Section~\ref{sec3}.

\begin{figure}[ht!]\anchor{fig4.1}
\cl{
\includegraphics{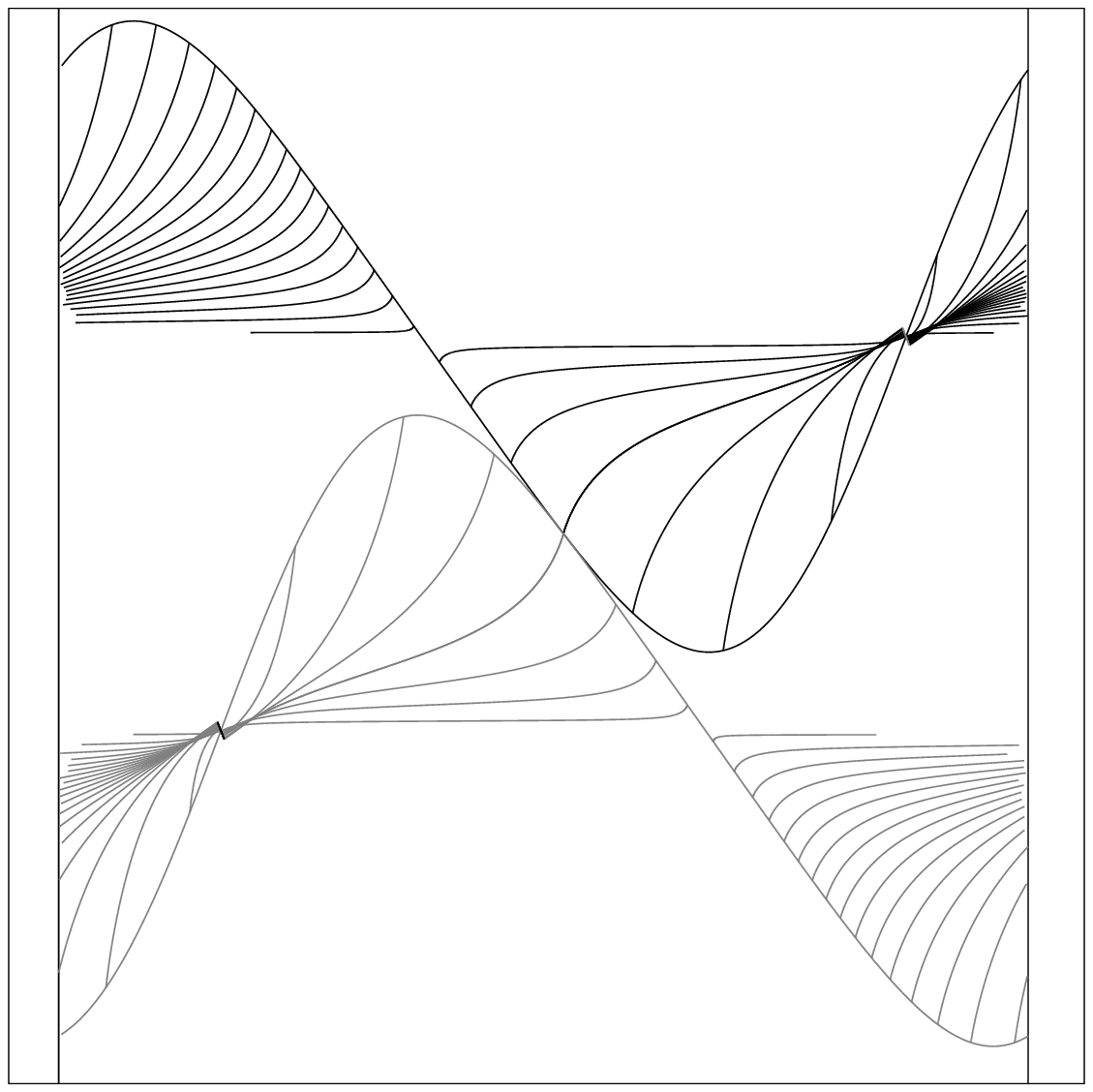}
}
\caption{}\label{fig4.1}
\end{figure}

\figref{fig4.1} shows $\Psi(\Sigma_1-E_0)$ in black and
$\Psi(\Sigma_2-E_0)$ in grey for the parameter $\underline s$.
All the black curves terminate on the tiny grey arc
$\Psi(Q_{21})$ and all the grey curves terminate
on the tiny black arc $\Psi(Q_{12})$.  The
Technical Lemma III says that the tiny grey
arc lies above the grey curve
$\Psi(C_2)$.  
This is obvious from the picture.

\begin{figure}[ht!]\anchor{fig4.2}
\cl{
\includegraphics{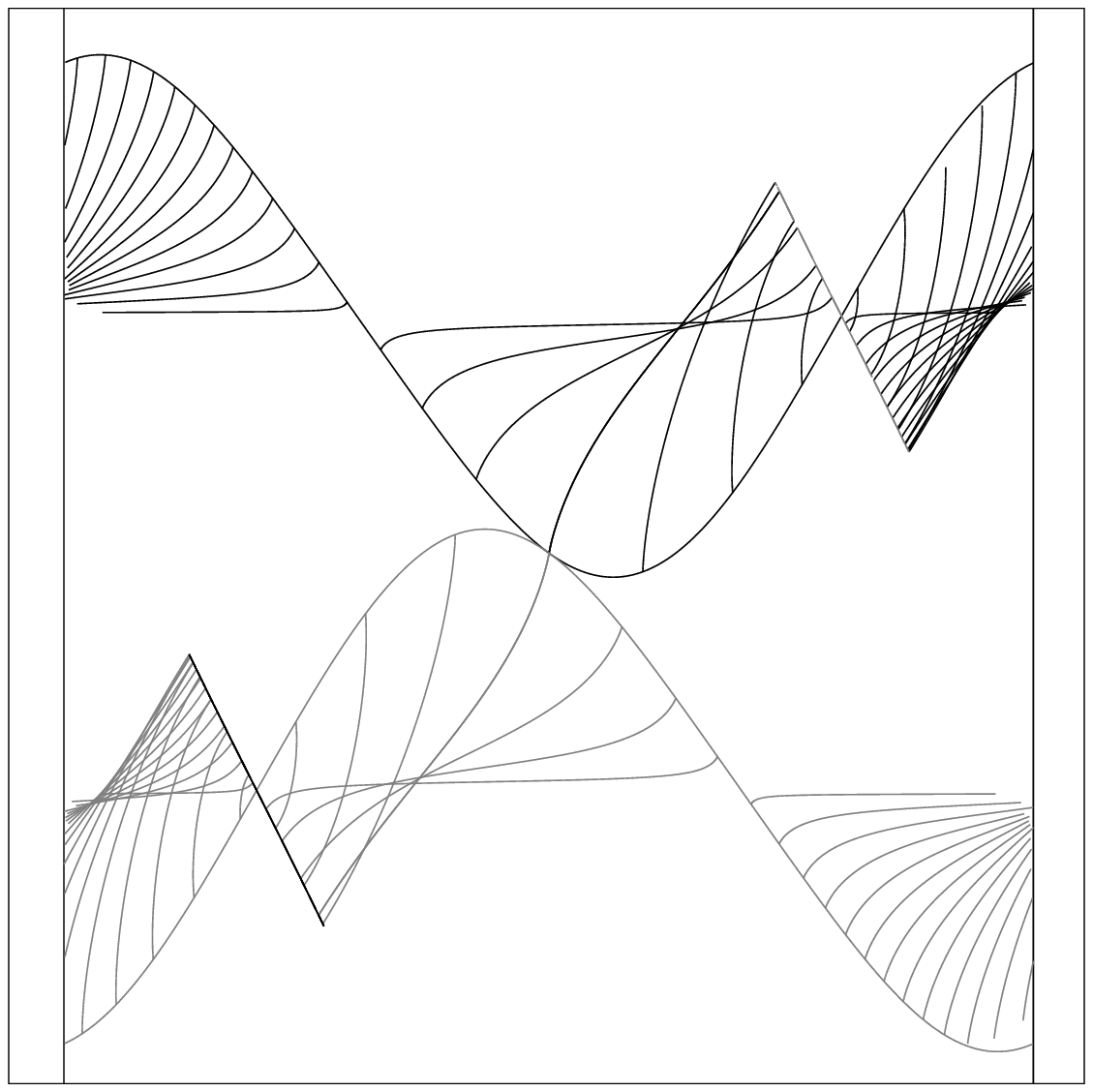}
}
\caption{}\label{fig4.2}
\end{figure}

\figref{fig4.2} shows the same picture for 
$s=1$.  Even though the parameter value $s=1$ is
outside the interval of interest to us, we
include the picture because the main features
are more dramatic.  Notice that the Technical
Lemma III remains true even at this parameter.
(However, our proof breaks down.)

Recall that $\Sigma_1=\Sigma_{10} \cup \Sigma_{12}$.
The $\R$--arcs foliating $\Sigma_{10}$
are mapped to individual points on $\Psi(C_1)$.
The $\R$--arcs foliating $\Sigma_{12}$
are mapped to black curves connecting
$\Psi(C_1)$ to the black arc
$\Psi(Q_{21})$.  Here, in general
\begin{equation}
Q_{ij}=Q_i \cap \Sigma_{ji}; \hskip 30 pt
i \not = j.
\end{equation}
The point here is that the $\R$--arcs foliating
$\Sigma(E_2,Q_2;C_1)$ start on $C_1$ and
end on a subset of $Q_2$.
The point $\Psi(p_0)$ is the point of tangency
between the black curve $\psi(C_1)$ and the
grey curve $\psi(C_2)$.  

\medskip
{\bf Remark}\qua Using our applet, the reader
can see the picture at any parameter, and can
recolor each individual object, so as to see
in a more direct way what part of the
picture corresponds to what object.

We are going to draw two conclusions from the
pictures.  The main thrust of our proof is
verifying that the pictures have the claimed
property.

\medskip
{\bf Contractible Intersections}\qua
Figures \fref{fig4.1} and \fref{fig4.2} indicate that $\Psi(\Sigma_1-E_0)$ and
$\Psi(\Sigma_2-E_0)$ only intersect in a
single point at the parameters
$s=\underline s$ and $s=1$ respectively.
Hence
$A=(\Sigma_1-E_0) \cap (\Sigma_2-E_0)$ is a single arc.
Note that $\Sigma_j \cap E_0$ is a connected
arc $Q_j^* \subset Q_j$.  The intersection
$A^*=Q_1^*  \cap Q_2^*$ is another arc, which intersects
$A$ in a single point.  Hence
$\Sigma_1 \cap \Sigma_2=A \cup A^*$, the union
of two arcs arranged in a $T$ pattern---a
contractible set.

\medskip
{\bf Embedded Spheres}\qua
Figures \fref{fig4.1} and \fref{fig4.2}  also suggest that
$\Sigma_1$ and $\Sigma_2$ are embedded.
We work this out here.
First, we have $\Psi(\Sigma_{10}-E_0) \subset \Psi(C_1)$
whereas the nontrivial arcs of $\Psi(\Sigma_{12}-E_0)$ only
intersects $\Psi(C_1)$ at one point.
Therefore  \begin{equation}
\label{embedd1}
\Sigma_{10} \cap \Sigma_{12} \subset C_1 \cup E_0.
\end{equation}
By symmetry
\begin{equation}
\label{embedd11}
\Sigma_{10} \cap \Sigma_{12} \subset C_1 \cup E_2.
\end{equation}
Lemma \ref{link1} below shows that
$E_0 \cap E_2=\emptyset$.
Hence $\Sigma_{10} \cap \Sigma_{12}=C_1$.
We already know that each hemisphere
$\Sigma_{1j}$ is an embedded disk,
and we've just seen that these hemispheres
just intersect along the equator.
Hence $\Sigma_1$ is an embedded sphere.
The same result holds for $\Sigma_2$ by
symmetry.

\medskip
{\bf Conclusions}\qua
If the pictures are right---namely
if $\Psi(\Sigma_1-E_0) \cap \Psi(\Sigma_2-E_0)$ is
a single point---then
$\Sigma_0,\Sigma_1,\Sigma_2$ are embedded spheres
which have pairwise contractible intersections.
But then these spheres bound balls
$B_0,B_1,B_2$ with pairwise disjoint interiors.
Moreover $I_j$ interchanges the two components
of $S^3-\Sigma_j$.  This picture, at the
parameter $s$, easily implies
that $\rho_s$ is a discrete embedding.
So, to prove the Goldman--Parker conjecture,
we just need to show that
 $\Psi(\Sigma_1-E_0) \cap \Psi(\Sigma_2-E_0)$ is
a single point for each $s \in [\underline s,\overline s)$.

\subsection{Containing the image}\label{subsec4.2}

\begin{figure}[ht!]\anchor{fig4.3}
\cl{\small
\psfrag {S12B}{$S_{12B}$}
\psfrag {S12A}{$S_{12A}$}
\psfrag {S21A}{$S_{21A}$}
\psfrag {S21B}{$S_{21B}$}
\includegraphics[width=3.5in]{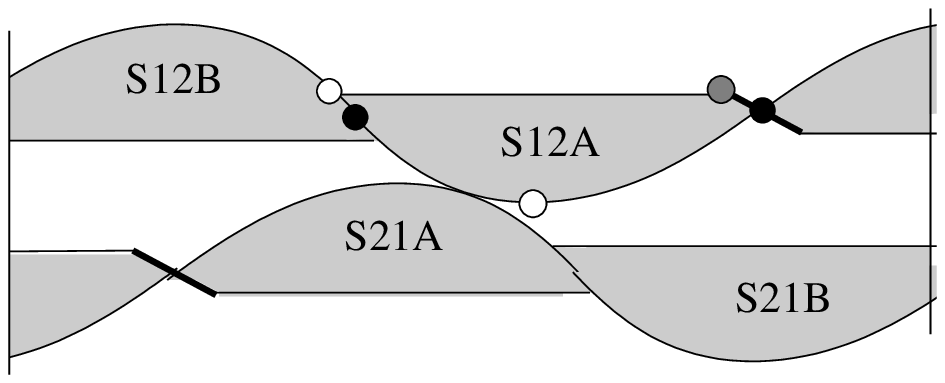}
}
\caption{}\label{fig4.3}
\end{figure}

We define $S_{12A}$ to be the
region of $\R/2 \pi \Z \times \R$ bounded by $\Psi(C_1)$,
$\Psi(Q_{21})$ and the horizontal line
through the highest point of
$\Psi(Q_{21})$, as shown in \figref{fig4.3}. 

From Technical Lemma III, 
the horizontal edge through the lowest
point of $\Psi(Q_{21})$ is disjoint from
$\Psi(C_2)$.  
We define $S_{12B}$ to be the region of
$\R/2 \pi \Z \times \R$ bounded by
$\Psi(C_1)$, $\Psi(Q_{21})$ and the
horizontal line through the lowest
point of $\Psi(Q_{21})$.  

We define
$S_{21A}$ and $S_{21B}$, similarly,
switching the roles of the indices
$1$ and $2$.

The $\R$--axis $R_1$ divides $\Sigma(E_2,Q_2;C_1)$
into two halves, which we denote by
$\Sigma_{12A}$ and $\Sigma_{12B}$.  
\figref{fig4.4} shows schematically how $\Psi(\Sigma_{12A})$ sits 
inside $S_{12A}$.  Compare \figref{fig4.2}.
Given Technical Lemma III, the Goldman--Parker Conjecture
comes down to proving the following two items,
which we establish in turn.

\begin{itemize}
\item 
Let $\alpha=\Psi(\gamma)$, where
$\gamma$ is an $\R$--arc of $\Sigma_{12A}$.
Then $\alpha$ does not intersect $\partial S_{12A}$ at
an interior point.

\item $\Psi(C_1)$ and $\Psi(C_2)$ intersect in
one point
\end{itemize}

Item 1 implies that
$\Psi(\Sigma_{12A}) \subset S_{12A}$ and also
that $\Sigma_{12A} \cap \Sigma_{10} \subset C_1$.
The same holds
for $\Sigma_{12B}$ and hence
$\Sigma_{12} \cap \Sigma_{10}=C_1$.  Hence
$\Sigma_1$ is embedded.
Item 2 combines with Item 1 and with the
Technical Lemma III to show that
$\Psi(\Sigma_1) \cap \Psi(\Sigma_2)$ is
a single point.

Our proof of the Goldman--Parker Conjecture
boils down to establishing Items 1 and 2.

\subsection{The height property}\label{subsec4.3}

\begin{figure}[ht!]\anchor{fig4.4}
\cl{\small
\psfrag {S}{$S_{12A}:$}
\psfrag {Y}{$\Psi(C_1)$}
\psfrag {YQ}{$\Psi(Q_{21})$}
\psfrag {hor}{horizontal edge}
\psfrag {neg}{negative edge}
\psfrag {pos}{positive edge}
\includegraphics[width=4in]{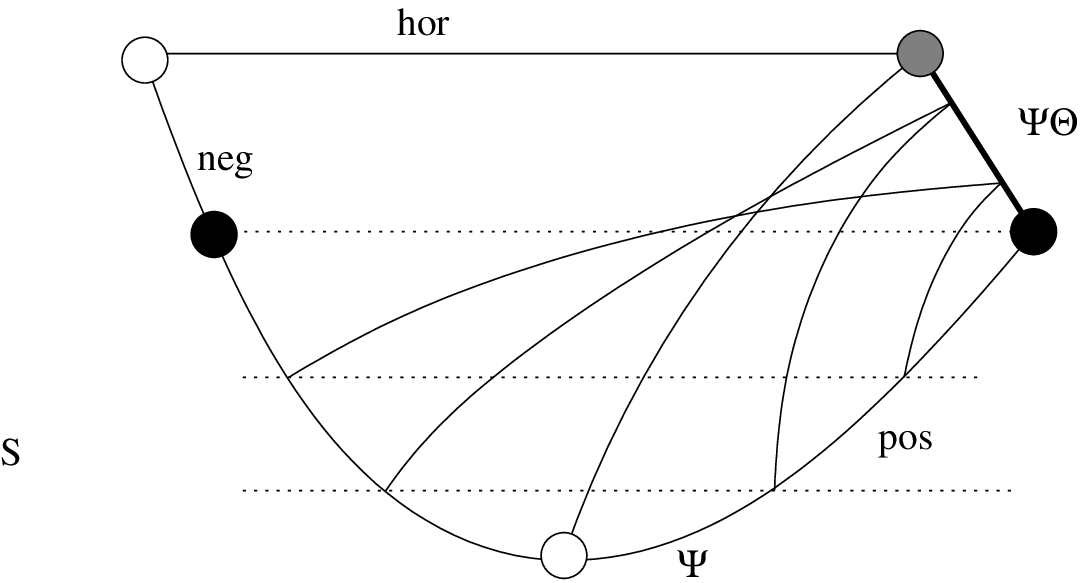}
}
\caption{}\label{fig4.4}
\end{figure}

\figref{fig4.4} shows an enlargement of part of
\figref{fig4.3}. 
As in \figref{fig4.4} we have 
distinguished $5$ points on $\partial S_{12A}$:
The $2$ black points are 
$\Psi(R_1-E_1)$, where $R_1$ is the
$\R$--axis of $\Sigma_1$.  
The $2$ white points are
the points of maximum and minimum height on
$\Psi(C_1)$.  
The grey point is the point of 
$\Psi(\Sigma_{12A})$ which lies on the same
horizontal line as the top white point.

\begin{lemma}
\label{heights}
Two arcs of $\Psi(\Sigma_{12A})$, which start at
the same horizontal level on $\Psi(C_1)$, end
at the same point of $\Psi(Q_{21})$.   In particular,
there is an arc of $\Psi(\Sigma_{12A})$ which
connects the bottom white point to the grey
point.
\end{lemma}

\proof
The second statement follows as a limiting
case of the first statement, so we will
concentrate on the first statement.
Let $\alpha_1$ and $\alpha_2$ be two
$\R$--arcs of $\Sigma_{12A}$.  Let
$a_j=\Psi(\alpha_j)$ and let
$\beta_j=I_1(\alpha_j)$.  \begin{enumerate}
\item
$a_1$ and $a_2$ have endpoints
at the same horizontal level of $\Psi(C_1)$.
\item Hence $\alpha_1 \cap C_1$ and $\alpha_2 \cap C_1$
lie in the same spinal sphere of the cospinal
foliation $S(E_0,Q_0)$.
\item Hence $\beta_1 \cap C_1$ and $\beta_2 \cap C_1$
lie in the same spinal sphere of the cospinal
foliation $S(E_0,Q_0)$.
\item
Hence $\beta_1$ and $\beta_2$ lie in the
same spinal sphere of the cospinal
foliation $S(E_0,Q_0)$. 
\item Hence $\beta_1 \cap Q_0=\beta_2 \cap Q_0$. 
\item Hence $\alpha_1 \cap Q_2=\alpha_2 \cap Q_2$. 
\item Hence $\Psi(\alpha_1 \cap Q_2)=
\Psi(\alpha_2 \cap Q_2)$. 
\item Hence $a_1$ and
$a_2$ contain the same point of $\Psi(Q_2)$.
\end{enumerate}
This is what we wanted to prove.
\endproof

Lemma \ref{heights} is borne out in
Figures \fref{fig4.1} and \fref{fig4.2} and can be seen in
more detail on our applet.

\subsection{The interlacing property}\label{subsec4.4}

Here we recall \figref{fig3.4} for convenience.
Recall that $\Sigma_1$ is {\it interlaced\/}
if the picture looks like the right hand side
of \figref{fig3.4}.
That is, a vertical line separates the minimum
of $\Psi(E_2)$ from the minimum of $\Psi(C_1)$.
In this section we prove that
$\Sigma_1$ is an interlaced
$\R$--sphere for all $s \in [\underline s,\overline s)$.
This result lets us use all the machinery from
Section~\ref{sec3}. 

\begin{lemma}
\label{link1}
The pairs $(C_1,E_0)$ and
$(E_0,E_2)$ are each generically linked.
\end{lemma}

\proof
We already know that $(E_0,C_1)$ are
linked.
Lemma \ref{link5} now says that
$(E_0,E_2)$ are linked.  We just need
the genericity condition.
If $C_1$ and $E_0$ bound perpendicular $\C$--slices
then $I_1$ stabilizes
$E_0$.  By symmetry $I_2$ stabilizes $E_0$.
But then $I_1I_2$ stabilizes a $\C$--slice.
This does not happen. 
If $E_0$ and $E_2$ bound perpendicular
$\C$--slices then $E_0$ and $E_1$ bound
perpendicular slices, by symmetry.  But then
the slices bounded by $E_1$ and $E_2$ are
disjoint.  Hence $E_1$ and $E_2$ are
unlinked.  This is a contradiction.
\endproof

Comparing \figref{fig3.4} with \figref{fig4.5} we see 
that $\Sigma_1$ is interlaced for the
parameter $\underline s$. 
By continuity, and Lemma \ref{link1},
we get that $\Sigma_1$ is interlaced
for all $s \in [\underline s,\overline s)$.  Lemma
\ref{link1} prevents the picture
from switching from 
the right hand side of \figref{fig3.4} to
the left hand side of \figref{fig3.4}
as $s$ varies.

\subsection{No interior cusps}\label{subsec4.5}

\begin{figure}[ht!]\anchor{fig4.5}
\cl{
\includegraphics{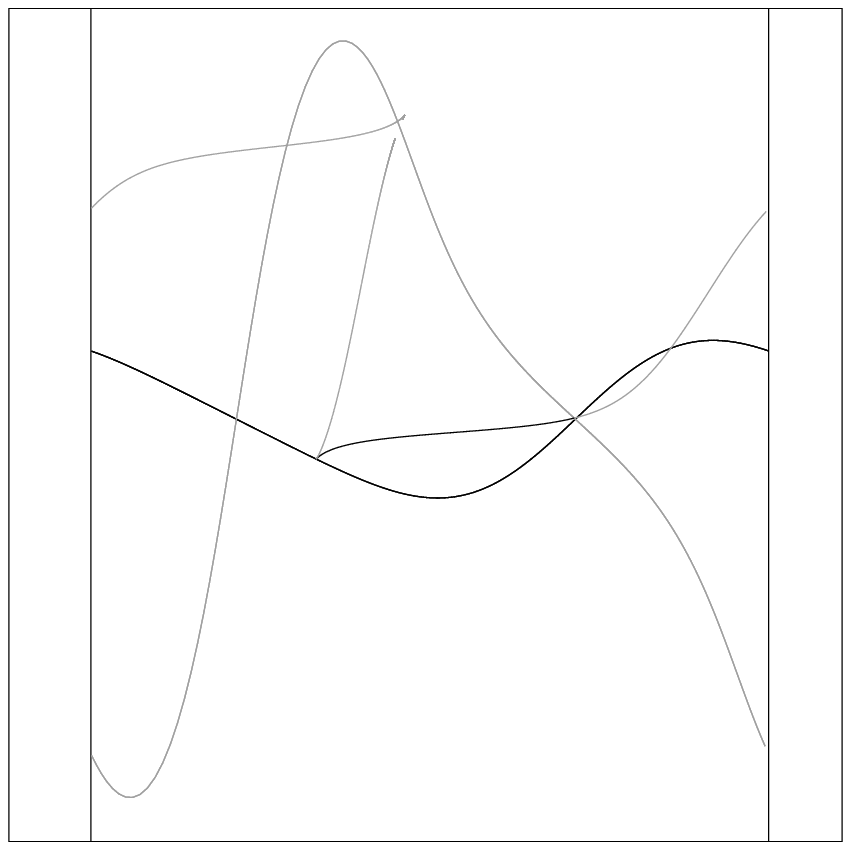}
}
\caption{}\label{fig4.5}
\end{figure}

Let $\gamma$ be an $\R$--arc of
$\Sigma(E_2,Q_2;C_1)$ and let
$\widehat \gamma$ be the $\R$--circle which
contains $\gamma$.  By definition, $\widehat \gamma$ is
an affiliate of $\Sigma(E_2,Q_2;C_1)$.
\figref{fig4.5} shows a particular example at the 
parameter $\underline s$.  Here
$\Psi(\gamma)$ in black
the rest of $\Psi(\widehat \gamma)$ in grey.
The loop $\Psi(C_1)$ is drawn in black and
the loop $\Psi(E_2)$ is drawn in grey.
The vertical direction in \figref{fig4.5}  is scaled 
differently than the vertical direction in \figref{fig4.1},
because otherwise $\Psi(E_2)$ would be quite a tall curve.
(On my
applet one can use many more colors.)

In general, let $s \in [\underline s,\overline s)$.
Suppose, for this parameter,
that $\gamma$ is an $\R$--arc of
$\Sigma(E_2,Q_2;C_1)$, not
contained in the $\R$--axis of $\Sigma_1$.  Let
$\widehat \gamma$ be the $\R$--circle
which contains $\gamma$. 
Let $\widehat \alpha=\Psi(\widehat \gamma)$ and
let $\alpha=\Psi(\gamma)$.   Let
$\alpha_0=\Psi(\gamma \cap C_1)$.

\begin{lemma}
$\alpha_0$ is a cusp of $\widehat \alpha$.
\end{lemma}

\proof
Let $\theta \in C_1$ be the endpoint
of $\gamma$.
Let $C_{\theta}$ be the contact plane
at $\theta$. Then $I_1$ rotates 
$C_{\theta}$ by $180$ degrees.
This means that $\gamma_{\theta}$
and $I_1(\gamma_{\theta})$ are
tangent at $\theta$.  But
$I_1(\gamma_{\theta})=\Sigma(E_0,q_0;\theta)$
is a fiber of $\Psi$.
Hence $\gamma_{\theta}$ is tangent to
a fiber of $\Psi$ at $\theta$.
Hence $\alpha_0=\Psi(\theta)$ is a
cusp.
\endproof

\begin{lemma}
$\alpha$ has does not contain the second cusp of $\widehat \alpha$.
\end{lemma}

\proof
Write $\alpha=\Psi(\gamma)$ and let
$\widehat \gamma$ be the $\R$--circle containing
$\gamma$.
\figref{fig4.5} shows an example where 
$\alpha$ does not contain the second cusp.
If there are examples where $\alpha$ does
contain the second cusp then there is
an example where the two endpoints of
$\alpha$ are both cusps. But then
this second cusp is the point of
maximum height of $\alpha$.  
Note that $\widehat \gamma$ intersects
$E_2$ in two points.  For the example
under consideration, we therefore have
$${\rm height\/}(\Psi(\widehat \gamma \cap Q_2))>
{\rm height\/}(\Psi(\widehat \gamma \cap (E_2- Q_2))).$$
However, for the example in \figref{fig4.5} 
$${\rm height\/}(\Psi(\widehat \gamma \cap Q_2))<
{\rm height\/}(\Psi(\widehat \gamma \cap (E_2- Q_2))).$$
But then, by continuity, we have an example where:
$${\rm height\/}(\Psi(\widehat \gamma \cap Q_2))=
{\rm height\/}(\Psi(\widehat \gamma \cap (E_2- Q_2))).$$
However $\gamma \cap Q_2$ and
$\gamma \cap (E_2-Q_2)$ are in harmonic
position with respect to $Q_2$.
Thus by Lemma \ref{rise}, $\Sigma_1$ is
symmetric.
However, Technical Lemma II
and Lemma \ref{comp3} combine to say that
$\Sigma_1$ is asymmetric.  This
is a contradiction. 
\endproof

\begin{corollary}[Rising Property]
Let $\gamma$ be an $\R$--arc of
$\Sigma(E_2,Q_2;C_1)$.  If
$\gamma$ is not contained in the
axis then $\alpha=\Psi(\gamma)$ is
nonsingular on its interior and
has nowhere vanishing slope.
\end{corollary}

\proof
This follows immediately from
Lemma \ref{local} and Lemma \ref{slope1},
combined with the fact that $\alpha$
has no cusps in its interior.
\endproof

\subsection{Proof of Item 1}\label{subsec4.6}

At the end of subsection \ref{subsec4.2} we reduced the whole
Goldman--Parker Conjecture to the verification
of two items. Here we prove the first of
these items.

Let $\gamma$ be an $\R$--arc of $\Sigma_{12A}$.  Let
$\alpha=\Psi(\gamma)$.  We want to show that
$\alpha$ does not intersect $S_{12A}$ at
an interior point.   We break $\partial S_{12A}$ into
three edges, as shown in \figref{fig4.4} 
From the Rising Property, the height of $\alpha$ attains
its maximum at the endpoint of $\alpha$ which
lies on $\Psi(Q_{21})$.  But the horizontal
edge only intersects $\Psi(Q_{21})$ at its
point of maximum height by Technical Lemma III.
Hence $\alpha$ does not hit the horizontal edge.
Since $\Sigma_1$ is remote $\alpha$ has positive
slope.  Hence $\alpha$ cannot hit the negative edge.
We just have to worry about the positive edge.

Before we deal with the problem of hitting the
positive edge, we
want to divide the arcs of
$\Psi(\Sigma_{12A})$ into two categories.
Say that a {\it positive arc\/} is one
whose endpoint lies on the positive edge
of $\partial S_{12A}$.  Likewise
define {\it negative arcs\/}.  Say that
the {\it middle arc\/} is the arc
which contains the lowest point of $\Psi(C_1)$.
This arc connects a white point to a grey point
in \figref{fig4.4}. 
Looking at \figref{fig4.2}
we can see that two positive arcs appear
never to cross each other whereas two
negative arcs always cross each other.

\begin{figure}[ht!]\anchor{fig4.6}
\cl{\small
\psfrag {a1}{$\alpha_1$}
\psfrag {a2}{$\alpha_2$}
\psfrag {a3}{$\alpha_3$}
\psfrag {detail}{detail}
\psfrag {x}{$x$}
\includegraphics{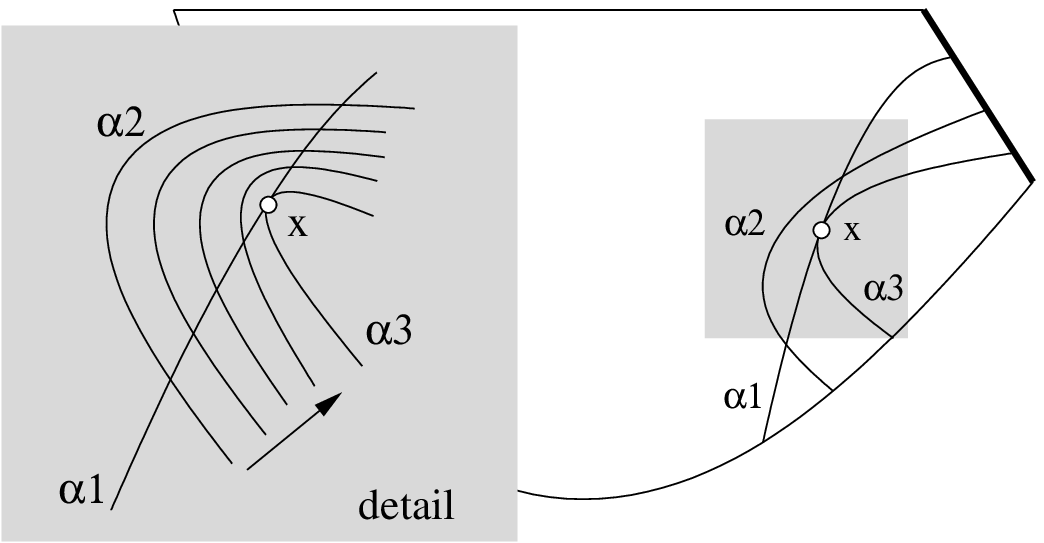}
}
\caption{}\label{fig4.6}
\end{figure} 

\begin{lemma}
\label{posedge}
Two positive arcs never cross each other
at interior points.
\end{lemma}

\proof
Looking at the ordering of the endpoints of
the positive arcs, on $\Psi(C_1)$ and $\Psi(Q_2)$
respectively, we see that two positive arcs
$\alpha_1$ and $\alpha_2$
must cross in at least one pair of oppositely oriented
crossings, if they cross at all.
By varying $\alpha_2$ away
from $\alpha_1$, as in \figref{fig4.6},
we can produce a positive arc
$\alpha_3$ which is tangent to $\alpha_1$
at some point $x$.  We write
$\alpha_j=\Psi(\gamma_j)$.  By the
Slope Principle of subsection \ref{subsec3.3}, the
$\R$--arcs $\gamma_1$ and $\gamma_3$ intersect
at some point of $\Psi^{-1}(x)$.  But
this contradicts the fact that $\Sigma_{12}$
is an embedded disk.
\endproof

\begin{corollary}
A positive edge does not intersect the
positive arc of $C_1$ in a point which
is interior to the positive edge.
\end{corollary}

\begin{figure}[ht!]\anchor{fig4.7}
\cl{\small
\psfrag {S}{$S_{12A}$}
\psfrag {e0}{$e_0$}
\psfrag {et}{$e_t$}
\psfrag {a0}{$\alpha_0$}
\psfrag {at}{$\alpha_t$}
\psfrag {xt}{$x_t$}
\psfrag {x0}{$x_0$}
\includegraphics[width=3in]{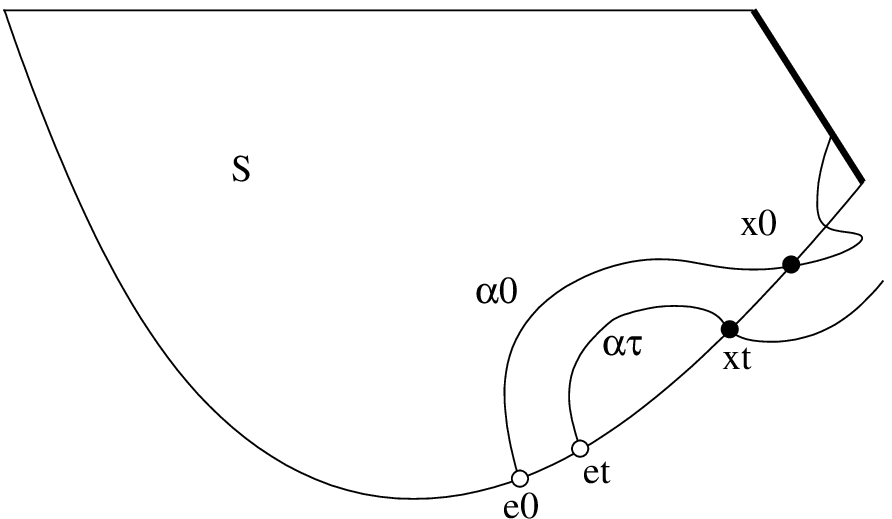}
}
\caption{}\label{fig4.7}
\end{figure}

\proof 
Suppose $\alpha_0$ is a positive arc which intersects
the positive edge at a point $x_0 \in \Psi(C_1)$.
By Lemma \ref{local} the arc
$\alpha_0$ is transverse
to $\Psi(C_1)$ at the endpoint.  Hence
an initial open segment of $\alpha_0$ either
is contained in the interior of
$S_{12A}$ or is contained in the complement of
$S_{12A}$.  \figref{fig4.7} shows the former
option, which turns out to be the true option.
Which option obtains
is independent of the choice of $\alpha_0$,
and for $\alpha_0$ very near the middle
arc, the former option obtains by the
Rising Property.  Hence, the former
option always obtains.  In summary,
some initial portion of $\alpha_0$ is
contained in the interior of $S_{12A}$.

We can choose $x_0$ to be the first point
where $\alpha_0$ crosses the positive edge.
Let
$$\{\alpha_t|\ t \in [0,1]\}$$ denote the family of
positive arcs such that the endpoint $e_t$
of $\alpha_t$ on $\Psi(C_1)$ interpolates between
the endpoint $e_0$ of $\alpha_0$ and the point
$e_1=x_0$.  
The arcs $\alpha_0$ and $\alpha_t$ cannot
cross, by the previous result.  Also,
an initial open segment of $\alpha_t$ is
contained in the interior of $S_{12A}$.
The only possibility is that
$\alpha_t$ crosses the positive edge at
some first point $x_t$ and,
in order along $\Psi(C)$ the points
come as $$e_0,e_t,x_t,x_0.$$
Since $\{e_t\}$ moves all the way
from $e_0$ to $x_0$ and
$e_t$ comes before $x_t$ (by the Rising
Property) we must have some value
$s$ for which $e_s=x_s$, but this
contradicts the fact that some
initial open segment of
$\alpha_s$ is contained in the
interior of $S_{12A}$.

\begin{lemma}
Suppose that $\alpha_1$ is a negative arc and
$\alpha_2$ is a positive arc.  Suppose
the endpoint of $\alpha_1$ on $\Psi(C_1)$ is
lower than the endpoint of $\alpha_2$ on
$\Psi(C_2)$.  Then $\alpha_1$ and $\alpha_2$
do not intersect.  
\end{lemma}

\proof
The ordering of the endpoints
for $\alpha_1$ and $\alpha_2$ is
the same as in Lemma \ref{posedge}.
The same argument as in 
Lemma \ref{posedge} works here.
\endproof

\begin{corollary}
A negative arc cannot intersect the positive
edge of $C_1$ in a point which is interior
to the negative arc.
\end{corollary}

\proof
Suppose $\alpha_1$ is a negative arc which
intersects the positive edge. Let
$x$ be the endpoint of $\alpha_1$
on $\Psi(C_1)$.  Then, from the Rising Property,
$\alpha_1$ intersects the positive edge at a point
$y$ which lies above $x$.  But then
$\alpha_1$ must cross some positive
arc $\alpha_2$ whose endpoint lies
above $x$ and below $y$.  This
contradicts the previous result.
\endproof

Our lemmas cover all the cases.
This establishes Item 1.

\subsection{Proof of Item 2}\label{subsec4.7}

Our goal is to show that $\Psi(C_1)$ and
$\Psi(C_2)$ intersect in a single point,
for all parameters $s \in [\underline s,\overline s)$.
We work in ${\cal H\/}$ and normalize so
that $(E_0,Q_0)$ is in standard position and
$\rho_0=(\R \times \{0\}) \cup \infty$.
Then $C_1$ and $C_2$ are swapped by the
$\R$--reflection $(z,t) \to (\overline z,-t)$
which fixes $\rho_0$.  To fix the scale,
we arrange that $C_1 \cap C_2=(1,0)$.
Both $C_1$ and $C_2$ have aspect greater than $9$
by Technical Lemma I, part 3.

\begin{lemma}
Define $\Psi_*(z,t)=(\arg z,t)$. Then
$\Psi_*(C_1) \cap \Psi_*(C_2)=(0,0)$,
with $\Psi_*(C_1)$ lying on top.
\end{lemma}

\proof
Let $s \in [\underline s,\overline s)$ be some parameter.
$\Psi_*(C_j)$ is the graph
of a function $\psi_j$.  Up to
rotations and scaling, $\psi$
satisfies the equation in Lemma \ref{calculus},
for some $A>9$.
Hence, by Lemma \ref{calculus},
the function $\psi_j$ is convex on an interval
of length $\pi$ and concave on the complementary
interval of length $\pi$.   We have parametrized
$\psi$ so that $\psi(0)=\Psi(p_0)$.  Here
$p_0=(1,0)=C_1 \cap C_2$.

\begin{sublemma} 
$\psi''(0)>0$ for all
$s \in [\underline s,\overline s)$.
\end{sublemma}

\proof
We compute that $\psi_1''(0)>0$ at the
parameter $\overline s$.  Suppose
that $\psi_1''(0)=0$ for some parameter $s$.
Then $\Psi_*(C_1)$ and $\Psi_*(C_2)$ are
tangent at their inflection point.
Then, by symmetry, $\Psi_*(C_1)=\Psi_*(C_2)$.
Here we are using Lemma \ref{calculus}
and the fact that $A(C_j)>9$.
Since $C_1$ and $C_2$ are generically
linked, this forces $C_1=C_2$,
a contradiction.  Hence
$\psi''(0)>0$ for all
$s \in [\underline s,\overline s)$.
\endproof

\begin{figure}[ht!]\anchor{fig4.8}
\cl{\small
\psfraga <-2pt,0pt> {t1}{$t_1$}
\psfraga <-2pt,0pt> {t2}{$t_2$}
\psfrag {L1}{$L_1$}
\psfrag {L2}{$L_2$}
\psfrag {Y1}{$\Psi(C_1)$}
\psfrag {Y2}{$\Psi(C_2)$}
\includegraphics[width=3.8in]{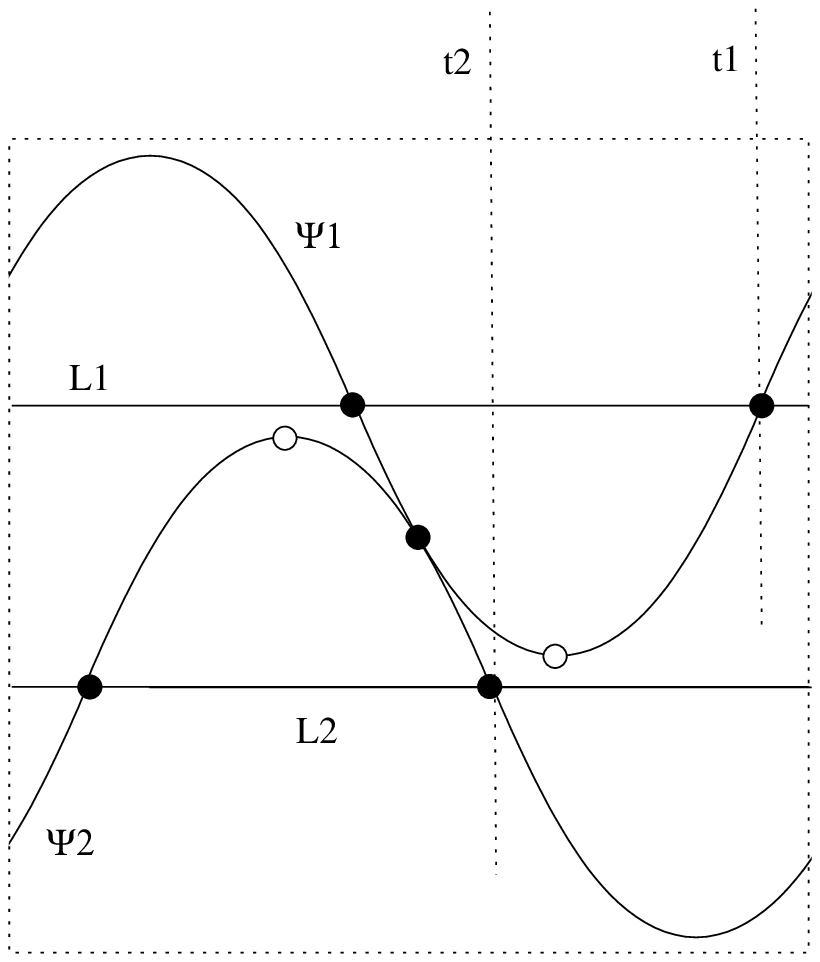}
}
\caption{}\label{fig4.8}
\end{figure}

We have $\Psi_*(p_0)=(0,0)$.
By symmetry we get
\begin{equation}
\label{initial}
\psi_1(0)=\psi_2(0);
\hskip 15 pt
\psi_1'(0)=\psi_2'(0);
\hskip 15 pt
\psi_1''(0)=-\psi_2''(0).
\end{equation}
The sublemma says that $\psi_1''(0)>0$.
Hence $\psi_2''(0)<0$ by symmetry.
To establish Item 2 it suffices to show
that $\psi_1(t)>\psi_2(t)$ for $t \in (0,\pi]$.
There are values $t_1,t_2 \in (0,\pi)$ such
that $\psi_j''(t_j)=0$.  The point
$(t_j,\psi_j''(t_j))$ is one of the points
of $\Psi_*(C_j) \cap L_j$.  Here $L_j$ is
the horizontal line through the inflection
points of $\Psi_*(C_j)$.
 We suppose
$t_2 \leq t_1$, as indicated in \figref{fig4.8}.
The other case is similiar.
Then on $(0,t_2)$ we have
$\psi_1>\psi_2$ because of our
initial conditions at $0$, and the
fact that $\psi_1''>0$ on $(0,t_1)$ and
$\psi_2''<0$ on $(0,t_1)$.
For $t \in [t_2,\pi)$ the curve
$\Psi_*(C_2)$ lies below $L_2$, 
and $\Psi_*(C_1)$ lies completely
above $L_2$ by Technical Lemma III (Section \ref{sec7}).
(The point here is that $L_j$ contains the
point $\Psi_*(c_j)$, where $c_j$ is the
center of $C_j$.)
This does it.
\endproof

\begin{figure}[ht!]\anchor{fig4.9}
\cl{\small
\psfrag {R0}{$\rho_0$} 
\psfrag {Z}{$Z$}
\psfrag {A}{$A$}
\psfraga <10pt,-2pt> {A2}{$A_2^*$}
\psfrag {g}{$\gamma$}
\psfrag {a}{$a$}
\psfrag {0}{$0$}
\psfrag {V}{$V$}
\psfrag {C1}{$C_1$}
\psfrag {C2}{$C_2$}
\psfrag {p0}{$p_0$}
\includegraphics[width=4.6in]{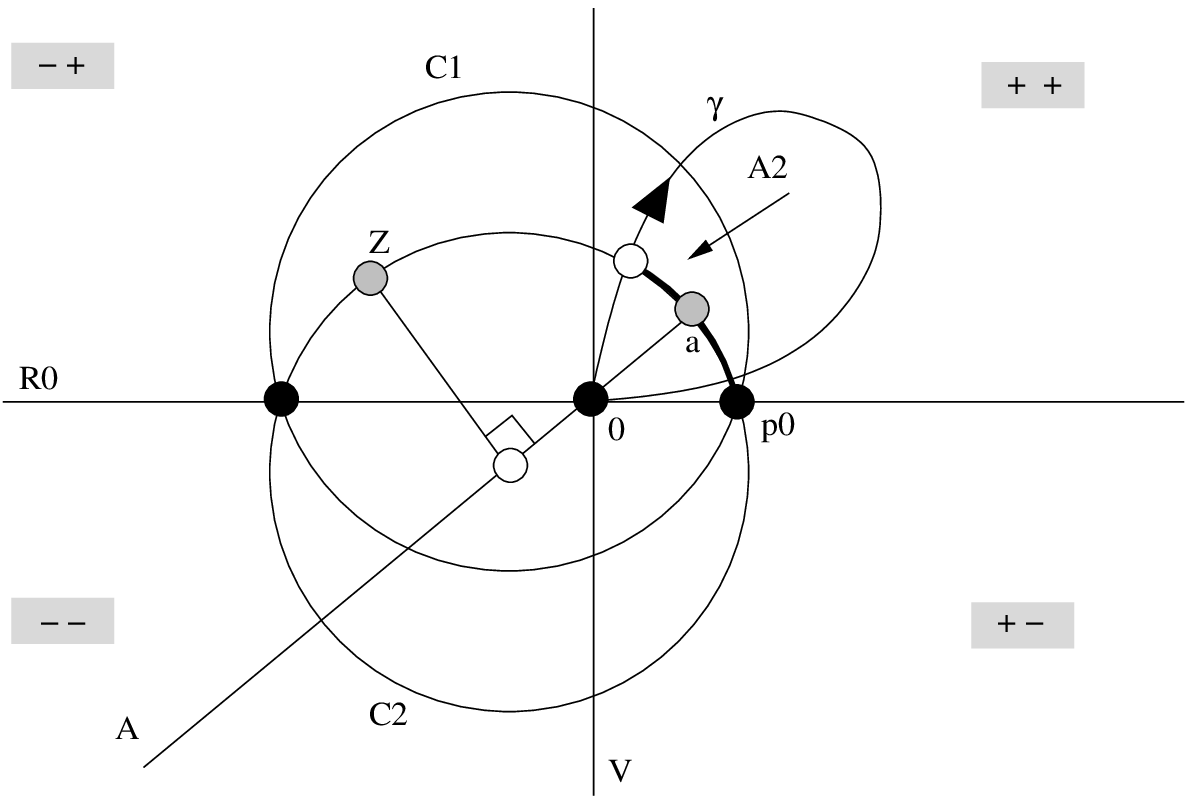}
}
\caption{}\label{fig4.9}
\end{figure}

To finish our proof
we compare the loxodromic
elevation map $\Psi$ with $\Psi_*$.
Let $A_2^* \subset C_2$
denote the arc such that $\Psi_*(A_2^*)$ 
ranges between $0$ and $\pi/2$ in the
$S^1$ direction on $S^1 \times \R$.
See \figref{fig4.9}.  Let $\pi$
be projection into $\C$.

\begin{lemma}
$\Psi_*(A_2^*)$ has negative slope.
\end{lemma}

\proof
Let $m_2^*$ denote the point on
$\Psi_*(C_2)$ having minimum
height.  We can locate 
$Z=\pi(\Psi_*^{-1}(m_2^*))$ as follows:
Let $A$ be the line through $0$ which
contains the point $a$ on $\pi(C_2)$ closest
to $0$.  
Then $Z$ is
obtained from $A$ by rotating
$90$ degrees about the center of $\pi(C_2)$,
as shown.
The key observations are that
$\pi(C_2)$ lies more in the lower half
plane than the upper half plane.
This property is
true for one parameter and cannot change
as the parameter varies, because the
$\pi(C_1)$ and $\pi(C_2)$
are symmetrically placed with respect
to the real axis, and never coincide.
Compare \figref{fig4.9}. 
From this observation, and item 3 of Technical Lemma I 
(subsec \ref{subsec5.4})
we conclude that $Z$ lies in the $(-,+)$
quadrant, as shown.

Let $A_2^{**}$ denote the arc of $C_2$
which projects to the $(+,+)$ quadrant.
Since $\Psi^{-1}(m_2^*)$ projects
to a point in the $(-,+)$ quadrant,
we conclude that $\Psi_*(A_2^{**})$ has
negative slope. Hence all points of
$A_2^{**}$ lie in the lower half-space of
$\cal H$.  But then the fibers of $\Psi$,
which contain points on
$A_2^{**}$, project to lemniscate lobes
which curve around counter-clockwise.
\figref{fig4.9} shows one such fiber $\gamma$.
In particular, if $\gamma$ contains the
endpoint of $A_2^*$ then $\pi(\gamma)$ is
tangent to $V$ at $0$.  (This is the
$\pi/2$ condition, which defines
$A_2^*$.)   But then $\pi(\gamma)$ intersects
$\pi(C_2)$ in a point which lies to the
right of $V$.  Hence $A_2^* \subset A_2^{**}$.
\endproof

\begin{lemma}
Suppose $x_2 \in C_2$ is a point such that
$\Psi(x_2) \in \Psi(C_1)$.  Then $x$
cannot be an interior point of $A_2^*$.
\end{lemma}

\begin{figure}[ht!]\anchor{fig4.10}
\cl{\small
\psfrag {C1}{$C_1$}
\psfrag {C2}{$C_2$}
\psfraga <0pt,2pt> {x1}{$x_1$}
\psfraga <-2pt,1pt> {x2}{$x_2$}
\psfrag {x2*}{$x_2^*$}
\psfraga <-2pt,0pt> {0}{$0$}
\psfrag {1}{$1$}
\psfraga <-2pt,0pt> {L}{$L$}
\psfrag {V}{$V$}
\psfrag {y}{$y$}
\psfrag {g}{$\gamma$}
\psfrag {detail}{detail}
\includegraphics[width=.9\hsize]{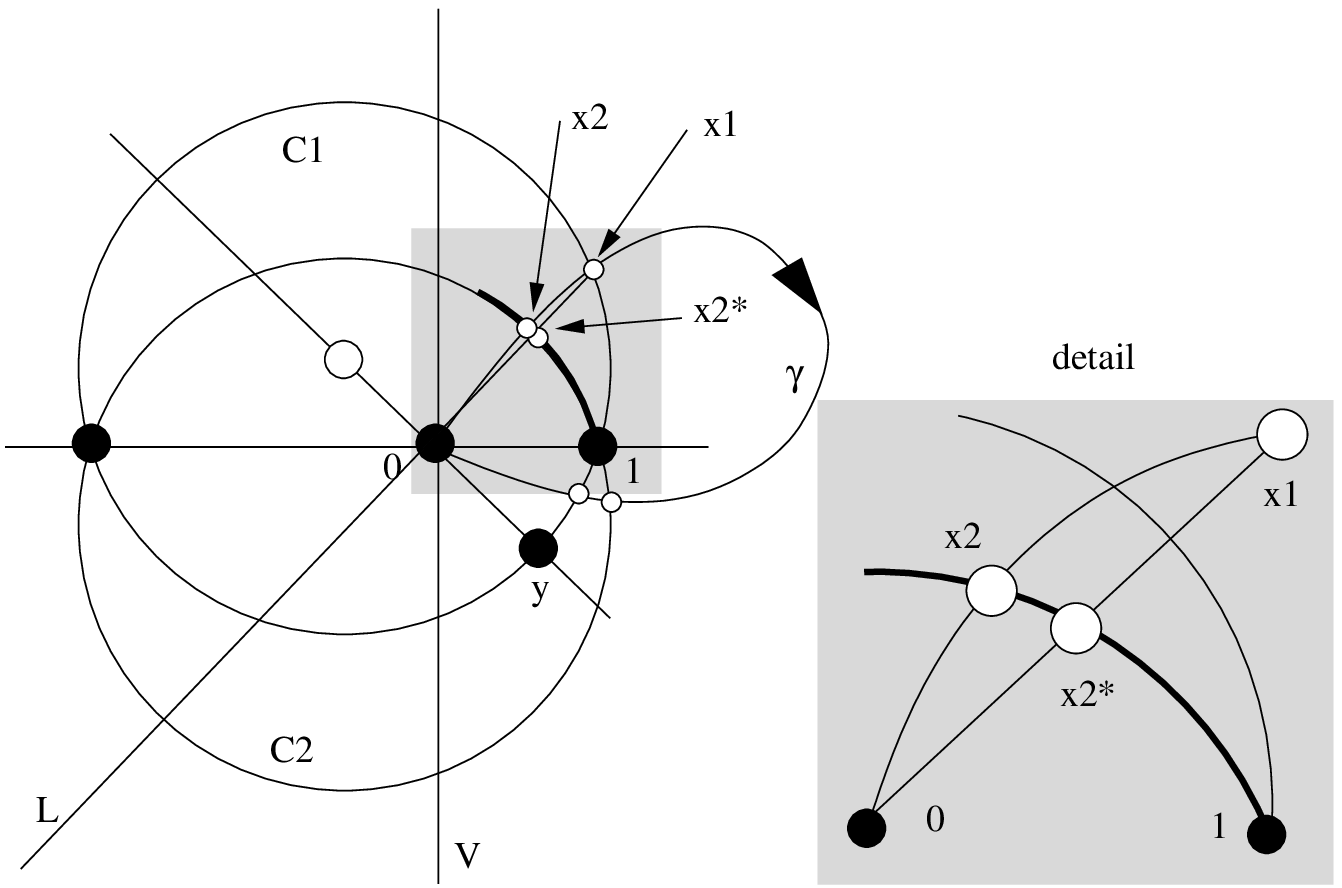}
}
\caption{}\label{fig4.10}
\end{figure}

\proof
We will suppose this is false and derive
a contradiction.  Let $x_2=(z_2,t_2)$.
There is some $x_1=(z_1,t_1) \in C_1$
and a fiber $\gamma$ of $\Psi$ such that
$x_1,x_2 \in \gamma$.
By the previous result $t_2<0$.
Hence $\pi(\gamma)$ curves around
clockwise, as shown in \figref{fig4.10}.

We trace counterclockwise around $\gamma$ as
indicated by the arrow in \figref{fig4.10}.  Suppose
for the moment that we encounter $x_2$ before
we encounter $x_1$, as shown in \figref{fig4.10}.
The height of $\gamma$ in $\cal H$ is
monotone decreasing.  Hence $t_2>t_1$.
The line $L$ through $0$ and $\pi(x_1)$
intersects $\pi(C_1)$ at some point 
between $1$ and $\pi(x_2)$.  Let
$x_2^*=(z_2^*,t_2^*)$ be the corresponding
point on $C_2$.  Then $x_2^*$ lies
between $(1,0)$ and $x_2$.  In
particular $x_2^* \in A_2^*$.
Since $\Psi_*(A_2^*)$ has negative
slope, we have $t_2^*>t_2$.
Hence $t_2^*>t_1$.
But then $\Psi_*(x_2^*)$ lies above
$\Psi_*(x_1)$, on the same vertical
line.  This is a contradiction.

To show that $\gamma$ encounters $x_2$ before
$x_1$ we first introduce some terminology.
For $j=1,2$ we say that a
{\it type $j$ arc\/} is a
portion of a fiber of $\Psi$ which
connects a point of $E_0-Q_0$ to 
$C_j$.   The type $j$ subarc of $\gamma$
projects to the portion of $\pi(\gamma)$
which connects $0$ to $x_j$.
Thus, we want to show that
the type $1$ sub-arc $\gamma_1 \subset \gamma$
contains the type $2$ sub-arc 
$\gamma_2 \subset \gamma$.
We are going to make
an argument based on a computer plot,
but we say in advance that we only
use topological features of the plot.
We simply need the plot to draw
the reader's attention to the
relevant details.

\begin{figure}[ht!]\anchor{fig4.11}
\cl{
\includegraphics{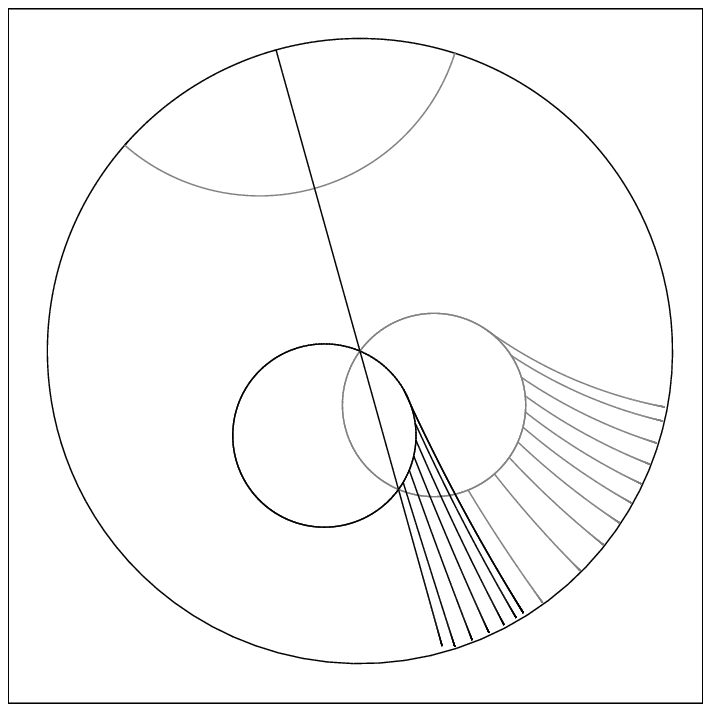}
}
\caption{}\label{fig4.11}
\end{figure}

We use the projection map $\eta$ from
subsection~\ref{subsec3.9}.  \figref{fig4.11} shows
the relevant objects at the parameter
$\underline s$.  The picture looks similar
at other parameters.
\begin{itemize}
\item The large black circle
is $\eta(E_0)$.  
The small black circle
is $\eta(C_1)$. 
The small grey circle
is $\eta(C_2)$.  
The grey geodesic has
$\eta(\partial Q_0)$ for endpoints.
\item The black geodesic, which
contains $\eta(C_1) \cap \eta(C_2)$,
is $\eta(\rho_0)$.
The black geodesic arcs are images
of the type $1$ arcs of interest to us.
$\eta(\gamma)$ must project onto the
same side of $\eta(R_1)$ as these black arcs,
and of course $\eta(\gamma)$ must intersect
the black circle. 
\item The grey arcs are images
of the type $2$ arcs of interest to us.
The last grey arc, the one tangent to
the grey circle, is the type $2$ subarc
which contains the point $z \in C_2$
such that $\Psi(z)$ has minimum height
on $\Psi(C_2)$.   Thus $\eta(\gamma_2)$
must project into the region indicated
by these grey arcs.
\end{itemize}

The type 1 and type 2 subarcs of $R_0$ coincide.
As we move away from $\eta(R_0)$ the 
projections of the type 1 subarcs grow
longer in comparison to the corresponding
projections of the type 2 subarcs.
In particular $\eta(\gamma_2)$ would
be shorter than $\eta(\gamma_1)$.
Hence $\gamma_2\subset \gamma_1$.
\endproof

In summary, if $\Psi(C_1)$ and $\Psi(C_2)$ intersect
in a second point $(x,y)$ then
$x \not \in (0,\pi/2]$.  By symmetry
$x \not \in [-\pi/2,0)$.  
Recall that $L_j$ is the horizontal line
in $S^1 \times \R$ which contains the
points of symmetry of $\Psi(C_j)$.
Then the two points of $L_1 \cap \Psi(C_1)$
are $\pi$ apart.   Let $s_1$ be this
symmetry point.  As we trace around
$\Psi(C_1)$ from $(0,0)$ to $s_1$ we
remain above $\Psi(C_2)$.  Once we reach
$s_2$ we remain above $L_1$ for another
$\pi$ radians of travel in the $S^1$
direction.  It follows from Technical Lemma III
that $L_1$ lies above every
point of $\Psi(C_2)$.
Depending on which option obtains, 
we have shown either
that $x \not \in (0,3 \pi/2]$
of $x \not \in [-3 \pi/2,0)$.
Either case implies the other by
symmetry.  But then $x=0$ and we
are done.

\section{Technical Lemma I}\label{sec5}
\setcounter{figure}{0}
\label{aspect}

\subsection{A resume of formulas}\label{subsec5.1}

Here we introduce the formulas which we use
for our technical estimates.  As a double-check,
we verified all the numbered equations
computationally, for the parameter $\overline s$.
  We will usually
suppress the parameter $s$ from our
notation. One should view this section
as a continuation of subsection \ref{subsec2.5}.
 First, some quantities from
 subsection \ref{subsec2.5}:

\begin{equation}
\label{parker0}
\beta=\frac{s+i}{\sqrt{2+2s^2}};
\hskip 30 pt A_1=\frac{s+17 i}{s+i}; \hskip 30 pt
A_2=\frac{12 \sqrt 2 i}{\sqrt{1+s^2}}.
\end{equation}
One verifies easily that
\begin{equation}
\label{alg1}
\overline \beta=\frac{3}{4} \times
\frac{A_1-1}{A_2}.
\end{equation}
The matrix $g_0=I_1I_0I_2$ has a positive
eigenvector $(e,\overline e,1)$.  The
quantity $e$ (not to be confused with
the base of the natural log) figures
heavily in our estimates.  As suggested
by the referee, we introduce the quantity:
\begin{equation}
\label{x}
x=\frac{e^2+|e|^2+\overline e^2}{1-|e|^2}.
\end{equation}
Here are $4$ equations, all due to the referee:
\begin{equation}
\label{parker5}
\label{PARKER}
|e|^2=\frac{x^2}{9(x-1)}
\end{equation}
\begin{equation}
\label{parker6}
|e-\overline e|^2=\frac{x(x-3)^2}{9(x-1)}
\end{equation}
\begin{equation}
\label{parker8}
|2\overline \beta e -1|^2=
\frac{(x-3)^2}{18(x-1)}
\end{equation}
\begin{equation}
\label{parker99}
\frac{288}{1+s^2}=|A_2|^2=\frac{9(x-1)(x-3)^2}{x}
\end{equation}
We will give the derivations below.  

Note
that $x>1$ by Equation \ref{PARKER}.
It follows readily from Equation \ref{parker99} that
\begin{equation}
\label{range}
s \in [\underline s,\overline s] \hskip 30 pt
\Longrightarrow \hskip 30 pt x \in [4,\frac{817}{200}].
\end{equation}
Recall that
$\rho_0$ is the $\R$--circle $\{(z,\overline z)\} \subset S^3$.
Recall also that $E_0$ is the $\C$--circle stabilized
by the element $g_0$.  We have
$E_0 \cap \rho_0=\{q_0,q_0'\}$ where
\begin{equation}
q_0=\left[\matrix{a \cr \overline a}\right];
\hskip 15 pt
q'_0=\left[\matrix{b \cr \overline b}\right];
\hskip 30 pt
\label{coeff1}
\label{coeff}
a=\frac{1+ i \sqrt{2|e|^2-1}}{2 \overline e} \hskip 15 pt
b=\frac{1- i \sqrt{2|e|^2-1}}{2 \overline e}.
\end{equation}
This last equation follows from the fact that
$$|a|^2=|b|^2=1/2; \hskip 30 pt
\langle (a,\overline a,1),(e,\overline e,1)\rangle=
\langle (b,\overline b,1),(e,\overline e,1)\rangle=0.$$
The endpoints of the arc $Q_0 \subset E_0$ are given by
$$
\partial_1Q_0=(c,d); \hskip 20 pt
\partial_2 Q_0=(\overline d,\overline c).
$$
We label so that $|c|^2 \geq 1/2$ and
$|d|^2 \leq 1/2$.
We have 
\begin{equation}
\label{norm}
|c|^2|d|^2=\frac{1}{x(x-3)^2};
\hskip 30 pt |c|^2+|d|^2=1.\end{equation}
We will give the derivation below.

It follows readily from Equation \ref{PARKER} that:
\begin{equation}
\label{parker199}
{\rm Re\/}(a\overline b)=\frac{-9+9x-x^2}
{2x^2}
\end{equation}
Below we will derive:
\begin{equation}
\label{massive}
{\rm Re\/}(a\overline \beta)=
\frac{3(x-1)(9+x)-(x-3) \sqrt{(x-3)(2x-3)(9-x)}}{16x^2}
\end{equation}
A similar derivation, which we omit, shows that:
\begin{equation}
\label{massive3}
{\rm Re\/}(\beta \overline b)=
\frac{3(x-1)(9+x)+(x-3) \sqrt{(x-3)(2x-3)(9-x)}}{16x^2}
\end{equation}
Finally:
\begin{equation}
\label{finale}
{\rm Im\/}(a\overline b + 
b \overline \beta + \beta \overline a)=
\frac{(x-3)\sqrt{(x-3)(2x-3)(x-1)}}{8x^2}
\end{equation}
I'm proud to say that I found these last
$3$ equations myself, though of course the referee
had a much better derivation for the last one,
which he communicated to me in his final comments
on the paper.

\subsection{The derivations}\label{subsec5.2}

The first $4$ derivations, as well as the last
one, are essentially due to the
referee.  We begin with a technical lemma:

\begin{lemma}
\begin{equation}
\label{parker4}
A_1=\frac{e^2-\overline e^4}{|e|^2-|e|^4} \hskip 30 pt
A_2=\frac{e^3-\overline e^3}{|e|^2-|e|^4} \hskip 30 pt
\frac{A_1-1}{A_2}=\frac{e}{x}+\overline e
\end{equation}
\end{lemma}

\proof
The vector $\widehat E_0$ is an eigenvalue for the
matrix $g_0$ from Equation \ref{matrix}. That is:
$$
\left[ \matrix{0&-1&0 \cr
-A_1&0&A_2 \cr
-A_2 &0&-\overline A_1}\right]
\left[\matrix{e \cr \overline e \cr 1}\right]=
\lambda \left[\matrix{e \cr \overline e \cr 1}\right]
$$
Reading across the top row we see that
$\lambda=-\overline e/e$.
Reading across the second and third rows,
using the relation $\overline A_2=-A_2$,
multiplying through by powers of $e$, and
conjugating when necessary, we get:
$$
-A_1 e^2+A_2 e = -\overline e^2
\hskip 30 pt
-A_1\overline e + A_2 \overline e^2=-e 
$$
The first two equations in 
Equation \ref{parker4} are now derived by the
usual method of elimination.  The third
equation is verified by expanding out both sides,
using the definition of $x$, and the
identity $e^3-\overline e^3=(e-\overline e)(e^2+|e|^2+\overline e^2)$.
\endproof

\subsubsection{Derivation of equation \ref{parker5}}\label{subsubsec5.2.1}

Inspecting Equation \ref{parker0}, we get the
equation
$9|A_1-1|^2=8|A_2|^2$.  Therefore:
\begin{equation}
\label{mess}
9|e^2-\overline e^4-|e|^2-|e|^4|^2=8|e^3-e^3|^2=
8|e-\overline e|^2\ |e^2+|e|^2+\overline e^2|^2
\end{equation}
The left hand side equals
$$9|e-\overline e|^2|e(1-|e|^2)+\overline e(e^2+|e|^2+\overline e^2)^2.$$
Therefore
$$9|e(1-|e|^2)+\overline e(e^2+|e|^2+\overline e^2)|^2=
8|e^2+|e|^2+\overline e^2|^2.$$
Dividing through by $(1-|e|^2)^2$ and cancelling the
term $|e-\overline e|^2$ which occurs on both sides
of Equation \ref{mess}, we get
$9|e+\overline e x|^2=8x^2$.  Expanding:
$$9(|e|^2+|e|^2x^2+e^2x+\overline e^2 x)=8x^2.$$
The left hand side can be re-written as
$$9(|e|^2+|e|^2x^2-|e|^2 x+ x(e^2+|e|^2+\overline e^2))=
9|e|^2(1+x^2-x)+x^2(1-|e|^2).$$
Therefore
$$9|e|^2(1+x^2-x))+9x^2(1-|e|^2)=8x^2.$$
Solving this last equation for $|e|^2$
yields Equation \ref{parker5}.

\subsubsection{Derivation of equation \ref{parker6}}\label{subsubsec5.2.2}

$$|e-\overline e|^2=
(e-\overline e)(\overline e -e)=
-e^2-\overline e^2 + 2|e|^2=
-(e^2+|e|^2+\overline e^2)+3|e^2|=$$
$$-x(1-|e|^2)+\frac{x^2}{3(x-1)}=
-x(1-\frac{x^2}{9(x-1)})+\frac{x^2}{3(x-1)}.$$
The last equation is equivalent to the
right hand side of Equation \ref{parker6}.

\subsubsection{Derivation of equation \ref{parker8}}\label{subsubsec5.2.3}

From Equation \ref{alg1} and the third part of
Equation \ref{parker4} we get
$$
2\overline \beta e
=\frac{3}{2x} \times (e^2+x|e|^2).$$
Substituting in Equation \ref{PARKER} and grouping
terms we get
\begin{equation}
\label{parker3009}
\frac{2x}{3} \times (2\overline \beta e -1)=
e^2+P; \hskip 30 pt
P=\frac{6x - 6 x^2+x^3}{9x-9}.
\end{equation}
Therefore
$$\frac{4x^2}{9} \times |2\overline \beta e-1|^2=
|e|^4+(e^2+\overline e^2)P + P^2=$$
\begin{equation}
\label{parker9}
|e|^4+(e^2+|e|^2+\overline e^2)P+P^2
-P|e|^2=
|e|^4+x(1-|e|^2)P + P^2-P|e|^2.
\end{equation}
When we simplify Equation \ref{parker9},
using Equation \ref{PARKER}, 
we get Equation \ref{parker8}.  

\subsubsection{Derivation of equation \ref{parker99}}\label{subsubsec5.2.4}

The second equality in Equation \ref{parker99} is
the nontrivial one.
From Equation \ref{parker4} we get
$$A_2=(e-\overline e) \frac{x}{|e|^2(1-|e|^2)}; \hskip 50 pt
|A_2|^2=\frac{x^2 |e-\overline e|^2}{|e|^2(1-|e|^2)}.$$
Now we plug in our equations for
$|e-\overline e|^2$ and $|e|^2$ and simplify.

\subsubsection{Derivation of equation \ref{norm}}\label{subsubsec5.2.5}

The values
$c,\overline d,e$ are all roots of
$$
z^3+\left[
\frac{2 \overline{A}_1}{A_2}\right] z^2+
\left[\frac{\overline{A}_1^2}{A_2^2}-\frac{A_1}{A_2^2}\right] z+
\frac{1}{A_2}
$$
This polynomial comes from solving the system of
equations entailed by the equation
$g_0(c,d,1)=\lambda(c,d,1)$.
The last term is the product of the roots.  Hence
$$
|c|^2|d|^2=\frac{1}{|A_2|^2|e|^2}=\frac{1}{x(x-3)^2}.
$$
Since $(c,d,1)$ is a null vector we also have
$|c|^2+|d|^2=1$.

\subsubsection{Derivation of equation \ref{massive}}\label{subsubsec5.2.6}

We will be a bit sketchy here, to avoid a huge mess
which is best done symbolically.
Note that
\begin{equation}
\label{massive2}
a\overline \beta=\frac{(2a\overline e)(2 \overline \beta e)}
{4|e|^2}=
\frac{(1+i \sqrt{2|e|^2-1})(\frac{3}{2x} \times (e^2+P)+1)}
{4|e|^2}.
\end{equation}
Here $P$ is as above.
The only term on the right hand side of
Equation \ref{massive2} which is
not readily expressible in terms of $x$ is the
$e^2$ term.  However, when we take twice the
real part of the right hand side of Equation
\ref{massive2}, which amounts to adding
this formula to the conjugate of itself,
the only terms not expressed in terms of
$x$ are the real and imaginary parts of $e^2$.
We have
\begin{equation}
\label{parker2009}
2 {\rm Re\/}(e^2)=e^2+\overline e^2=x(1-|e|^2)-|e|^2;
\hskip 30 pt
({\rm Re\/}(e^2))^2+({\rm Im\/}(e^2))^2=(|e|^2)^2.
\end{equation}
Using these equations and Equation \ref{PARKER}
we get expressions for the real and imaginary
parts of $e^2$ in terms of $x$. When we simplify the
massive expression we arrive at Equation \ref{massive}.

\subsubsection{Derivation of equation \ref{finale}}\label{subsubsec5.2.7}

Noting that $a\overline e = e \overline b$ we can write:
$$
a\overline b+b\overline \beta + \beta\overline a=
\frac{(2a\overline e)^2+
2\overline a e (2\beta \overline e + 2\overline \beta e)}{4|e|^2}
$$
Using Equations \ref{coeff} and \ref{parker3009},
\ref{parker2009} we find that this last
expression equals:
$$\frac{(1+i \sqrt{2|e|^2-1})^2+
\frac{3}{2}(1-i \sqrt{2|e|^2-1}) (1+|e|^2-|e|^2/x)}{4|e|^2}$$
Using Equation \ref{PARKER} and expanding, we get
Equation \ref{finale}.

\subsection{Items 1 and 2}\label{subsec5.3}

We turn now to the proof of the Technical Lemma I.
Here $A$ denotes the aspect of $C_1$ relative to $E_0$.
A {\it polar vector\/} to a $\C$--circle $C$ is a
vector $X$ such that $\langle X,\widehat y\rangle=0$
whenever $\widehat y$ is a lift of a point on $C$.
The polar vector for $C$ is unique up to scale.

\begin{lemma}
\label{parker11}
Let $\widehat E_0$ and $\widehat C_1$ be polar vectors
to $E_0$ and $C_1$ respectively.  Let $A$ be the
aspect of $C_0$ when $E_1$ is normalized as above.
Then:
\begin{equation}
\label{parker1}
A=\frac{|\langle \widehat E_0,\widehat E_0 \rangle |
|\langle \widehat C_1,\widehat C_1 \rangle |} 
{|\langle \widehat C_1, \widehat E_0 \rangle |^2}
\end{equation}
Moreover, the two $\C$--circles are linked
provided that $A>1$.
\end{lemma}

\proof
The statement about the linking comes straight from
\cite[subsection 3.3.2]{G}.  Now for the second statement:
If $E_0$ and $C_1$ are normalized as in Lemma
\ref{computation14} then the
polar vectors to $E_0$ and
$C_1$, in the Siegel model, are:
\begin{equation}
\label{parker2}
\widehat E_0=\left[\matrix{0\cr 1\cr0}\right] \hskip 30 pt
\widehat C_1=\left[\matrix{r^2-d^2 \cr \sqrt 2 du \cr 1}\right]
\end{equation}
Here $u$ is a unit complex number which is a real
multiple of the projection of the center of $C_1$ to
$\C$. 
For $E_0$ this is obvious and for $C_1$ it is proved
by showing that any two distinct points on $C_1$ are
$\langle,\rangle'$ orthogonal to $\widehat C_1$.
From here an easy computation shows
that the left hand side of Equation \ref{parker1},
computed with the Siegel Hermitian form, 
yields $r^2/d^2=A$.
\endproof

In the ball model we have
\begin{equation}
\widehat E_0=(e,\overline e,1); \hskip 30 pt
\widehat C_1=(0,2 \overline \beta,1).
\end{equation}
In fact $\widehat E_0$ is the positive eigenvector
for $g_0$ and the eigenvalue is $-\overline e/e$.
From Lemma \ref{parker11} and Equations \ref{PARKER}
and \ref{parker8} we have
\begin{equation}
\label{part2}
A=\frac{2|e|^2-1}{|2\overline \beta e-1|^2}=
\frac{2(2x-3)}{x-3}.
\end{equation}
For $x$ as in Equation \ref{range} we
compute that $A \in [9.5,10]$.  This establishes
Items 1 and 2 of the Technical Lemma I.

\subsection{Item 3}\label{subsec5.4}

Here, as in Section~\ref{sec2},
$$C_1 \cap C_2=p_0=(\beta,\overline \beta,1).$$
Let $q_0$ and $q_0'$ be as in subsection \ref{subsec5.1}.
Let $\B\co  S^3 \to \cal H$ be the
Heisenberg stereographic projection which
normalizes as in item 3 of Technical Lemma I 
(subsec \ref{subsec5.4}).
Recall that $\pi(z,t)=z$.

\begin{lemma}
Let $q_0''=(-a,-\overline a)$.  Then
$\pi \circ \B(q_0'')$ is the
second intersection point of
$\pi(C_1) \cap \pi(C_2)$.
\end{lemma}

\proof
Let $V$ denote the $\C$--circle such
that $\B(V)$ is vertical and
$\pi(\B(v))$ is the second intersection
point of $\pi(C_1) \cap \pi(C_2)$.
By symmetry $V$ lies on the
Clifford torus and is stabilized
by the map $(z,w) \to (\overline w,\overline z)$.
Also $V$ contains
$(a,\overline a)$.  Hence
$V$ also contains $q_0''$.
\endproof

We introduce the cross ratio
\begin{equation}
\chi(z_1,z_2,z_3,z_4)=\frac{(z_1-z_3)(z_2-z_4)}{(z_1-z_2)(z_3-z_4)}.
\end{equation}
Let $\zeta \in (-\infty,0)$.  
Note that $\chi(\infty,\zeta,0,1)>2$
iff $\zeta<-1$.  Note also that
the restriction of $\Psi_*$ to $\rho_0$ is
a Mobius transformation.  Hence
$\zeta=\Psi_*(-a,-\overline a)$ is closer
$0=\Psi_*(b,\overline b)$ than is
$1=\Psi_*(\beta,\overline \beta)$ iff
\begin{equation}
f(s):=\chi(a,-a,b,\beta)=\frac{(a-b)(-\beta-a)}{2a(b-\beta)}>2.
\end{equation}
Equivalently (since $|a|^2=|b|^2=|\beta|^2=1/2$) 
it suffices to show that
\begin{equation}
\frac{|a-b|^2|\beta+a|^2}{|\beta-b|^2}=
\frac{(1-2\  {\rm Re\/}(a\overline b))
(1+2\ {\rm Re\/}(a\overline \beta))}
{1-2\  {\rm Re\/}(b \overline \beta) }
>8
\end{equation}
It is an exercise in calculus to show that the
the quantities
${\rm Re\/}(a\overline b)$ and
${\rm Re\/}(a\overline \beta)$ and
${\rm Re\/}(b\overline \beta)$
are all monotone for $x$ as in
Equation \ref{range}.
To sketch the idea, let $f_1$ be the function
from Equation \ref{massive}.  
We compute explicitly that $f_1'(4)=-1/16$ and we easily
get the crude bound $|f_1''|<1$ on $[4,5]$.  Hence
$f_1'<0$ for $x$ in our range.

Computing at the endpoints of $[\underline s,\overline s]$
we have:
\begin{equation}
\label{explicit}
{\rm Re(a\overline b)\/}  \in [.331,.344]; \hskip 20 pt
{\rm Re(a\overline \beta)\/} \in [.432,.438]; \hskip 20 pt
{\rm Re(b\overline \beta)\/} \in [.474,.477].
\end{equation}
Hence 
\begin{equation}
\frac{(1-2\  {\rm Re\/}(a\overline b))
(1+2\ {\rm Re\/}(a\overline \beta))}
{1-2\  {\rm Re\/}(b \overline \beta)} \geq
\frac{(1-2(.334))(1+2(.438))}{1-2(.474)}=
11.9775.
\end{equation}
This completes our proof.

\medskip
{\bf Remark}\qua All our points in $\C$ are confined to
a single quadrant.  Hence, the imaginary parts
also vary monotonically.  We compute:
\begin{equation}
\label{impart}
{\rm Im\/}(b \overline a) \in [-.374,-.363] \hskip 15 pt
{\rm Im\/}(a \overline \beta) \in [.242,.252] \hskip 20 pt
{\rm Im\/}(\beta \overline b) \in [.151,.171]
\end{equation}

\section{Technical Lemma II}\label{sec6}
\setcounter{figure}{0}

\subsection{Estimating $Q_0$}\label{subsec6.1}

Recall that $E_0$ intersects the Clifford torus
in points $(a,\overline a)$ and $(b,\overline b)$,
where $a$ and $b$ are as in Equation \ref{coeff}.
The vectors $(a,\overline a,1)$ and $(b,\overline b,1)$
have real Hermitian inner product.
Therefore the point $P_r$, represented
by the vector
\begin{equation}
r \left[\matrix{a\cr \overline a \cr 1}\right]+
i \left[\matrix{b \cr \overline b \cr 1} \right]
\end{equation}
lies in $E_0$.   The idea here is that the chosen
lift of $P_r$ is null and lies in the span
of $(a,\overline a,1)$ and $(b,\overline b,1)$.
As $r \to \infty$ the point $P_r$ converges to $(a,\overline a)$
and $P_0=(b,\overline b)$.  The fixed points of
$g_0$ have the form $P_r$ and $P_{-r}$ where
$r=r(s)$.  We define $r(\overline s)=\infty$,
in keeping with the analysis just made.

\begin{lemma}
\label{rbound}
$r(s)>3$ for all $s \in [\underline s,\overline s]$.
\end{lemma}

\proof
Equation \ref{norm} gives us
$|d|^2(1-|d|^2)>.2079$.  It follows easily that
\begin{equation}
\label{dbound}
|d|^2>.294
\end{equation} when $x$ is as
in Equation \ref{range}.  This also
holds for the smaller range
$s \in [\underline s,\overline s]$.

We will suppose that there is some $s \in [\underline s,\overline s]$
such that $r(s)=3$ and we will derive a contradiction.
If $r(s)=3$ then we have
$$d=\frac{3\overline a + i \overline b}{3+i}.$$
Since $|3+i|^2=10$ we have
$$|3\overline a + i\overline b|^2>2.94$$
Using $|a|^2=|b|^2=1/2$ and expanding:
$$5+3i(a\overline b-\overline a b)=5-6\ {\rm Im\/}(a\overline b)>2.94.$$
Hence
\begin{equation}
\label{parker66}
{\rm Im\/}(a\overline b)<.343333.
\end{equation}
This contradicts Equation \ref{impart}, which
says that
${\rm Im\/}(a\overline b)>.36$ throughout
$[\underline s,\overline s]$. 
\endproof

\subsection{Item 1}\label{subsec6.2}
\label{height}

There exists a (unique) Heisenberg stereographic
projection $\B$ 
which maps $E_0$ to $(\{0\} \times \R) \cup \infty$,
and has the following properties:
\begin{itemize}
\item $\B(C_1)$ is centered at a point
on the positive real axis.
\item $\B(C_1)$ projects to a circle
in $\C$ of radius $1$.
\item The endpoints of $\B(Q_0)$ are
symmetrically located on $\B(E_0)$.
Hence the $\R$--axis for $\Sigma_1$ is
just the real line in $\C \times \{0\}$.
\end{itemize}
The Siegel model polar vector for $\B(C_1)$ is:
\begin{equation}
\label{polar3}
V_3=\left[\matrix{D \cr \sqrt 2 d \cr 1}\right]
\hskip 30 pt D=1-d^2
\end{equation}
The Siegel model vectors representing
the endpoints of $\B(Q_0)$ are given by:
\begin{equation}
\label{polar1}
V_1=\left[\matrix{iu \cr 0 \cr 1}\right] \hskip 50 pt
V_2\left[\matrix{-iu \cr 0 \cr 1}\right]
\end{equation}
In this section we prove the following result,
which implies item 1 of the
Technical Lemma II.

\begin{lemma}
\label{wide}
$u>4.2$ for all $s \in [\underline s,\overline s]$.
\end{lemma}

Given $3$ vectors $V_1,V_2,V_3 \in \C^{2,1}$ we
define:
\begin{equation}
\label{ai}
\delta(V_1,V_2,V_3)=\frac{|{\rm Im\/}(\tau)|}{|{\rm Re\/}(\tau)|}
\hskip 30 pt \tau =
\langle V_1,V_2 \rangle
\langle V_2,V_3 \rangle
\langle V_3,V_1 \rangle
\end{equation}
Using the form in Equation \ref{Siegelform} we
compute readily that:
\begin{equation}
\label{parker1009}
\delta(V_1,V_2,V_3)=
\frac{|D^2-u^2|}{2Du}
\end{equation}
$V_1,V_2,V_3$ correspond to the following
vectors in the ball model:
\begin{equation}
P_r \hskip 30 pt -P_{-r} \hskip 30 pt \widehat C_1=
\left[\matrix{0 \cr 2 \overline \beta \cr 1}\right]
\end{equation}
(We prefer to use $-P_{-r}$ in place of $P_{-r}$.)
It is convenient to define:
\begin{equation}
\label{Ua}
U_a=\langle \left[\matrix{a \cr \overline a \cr 1}\right],
\widehat C_1 \rangle=2 \overline a \beta-1 \hskip 30 pt
U_b=\langle \left[\matrix{b \cr \overline b \cr 1}\right],
\widehat C_1 \rangle=2\overline b \beta-1 \hskip 30 pt
\end{equation}
Then
$$ \langle P_r,\widehat C_1 \rangle=r U_a + i U_b; \hskip 10 pt
 \langle \widehat C_1,-P_{-r} \rangle
=r \overline U_a + i \overline  U_b; \hskip 10 pt
 \langle P_r,-P_{-r} \rangle \in i\R.$$
From this we get:
\begin{equation}
\label{delta2}
\delta(V_1,V_2,V_3)=
\bigg|\frac{r|U_a|^2-r^{-1}|U_b|^2}
{{2\rm Re\/}(U_a \overline U_b)}\bigg|
\end{equation}
From Equation \ref{explicit} we have:
$$|U_a|^2=|1-2\overline a \beta|^2=
2-4\ {\rm Re\/}(a\overline \beta)>.248
$$
$$|U_a|^2=|1-2\overline b \beta|^2=
2-4\ {\rm Re\/}(b\overline \beta)<.104
$$
Using Equations \ref{parker199}, \ref{massive},
\ref{massive2}  and \ref{Ua} we compute that:
$$2{\rm Re\/}(U_a \overline U_b)=
\frac{-3(3-4x+x^2)}{2x^2}<.301$$
for $x$ as in Equation \ref{range}.
Therefore
$$
\delta(V_1,V_2,V_3)=\delta(P_r,-P_{-r},\widehat C_1)
\geq \frac{3(.248)-(1/3)(.104)}{(.301)}> 2.35
$$
Combining our last result with
Equation \ref{parker1009} we get
$|-D^2+u^2| \geq 4.7Du.$
Since $u>0$ for all parameters,
the quantity $-D^2+u^2$ cannot change sign.
For otherwise we would have $0>Du$.
Hence $u^2>4.7Du$.  From Technical
Lemma I we have:
$$D=1-\frac{1}{d^2}>1-\frac{1}{9.5}>.894$$
Hence 
$u>(4.7)(.894)>4.2.$  This completes our proof.

\subsection{Item 2}\label{subsec6.3}

We normalize by the map $\B$
so that all our calculations take place in $\cal H$.
So, $C_1$ is a $\C$--circle, centered on a point
$c_1 \in \R^+$ and $C_1$ projects to a circle
of radius $1$.  Since the aspect of $C_1$ is
at least $9$ we conclude that $c \in (0,1/3].$
Now $I_1$ is a $\C$--reflection in $C_1$.  The
$x$--axis is precisely $R_1$, the axis
for $\Sigma_1$.  Note that $R_1$ intersects
$C_1$ twice, at points which are precisely
$2$ units away from each other.
The restriction of $I_1$
to $R_1$ is a Mobius transformation.

\begin{lemma}
$E_2=I_1(E_0)$ intersects $R_1$ in two points.
One of the points is $(c,0)$ and the
other one is $(-1/c+c,0)$.
\end{lemma}

\proof
This follows from the fact that the
restriction of $I_1$ to $\R^1$ is an
inversion in the segment of radius $1$
centered at $0$.
\endproof

Thus we see that $E_2$ projects to a circle
diameter $1/c \geq 3.$
Moreover, the center of $E_1$ is 
$(-1/2c+c)<-7/6.$
Now $E_1$ is contained in the contact
plane based at its center.  From this we
see that $E_2$ is contained in a
contact plane of slope at least $7/3$.
(Here, and below, {\it slope\/} means
vertical rise divided by horizontal run.)
Compared to $C_1$, the
$\C$--circle $E_2$ is a big and tall
set.

\begin{lemma}
\label{extend1}
Both endpoints of $Q_2$ are at
most $1/4$ from 
the horizontal plane
$\C \times \{0\}$.
\end{lemma}

\proof
Working in the Siegel model, the endpoints of
$Q_0$ have lifts $V_1$ and $V_2$, as Equation \ref{polar1}
and the polar vector $V_3$ for $C_1$ is as in
Equation \ref{polar3}.
Using Equation \ref{reflection} relative to the
Siegel form we have
$$I_{C_1}(V_1)=
-\left[\matrix{iu \cr 0 \cr 1}\right]+
2 \frac{iu+D}{2D-2d^2} \left[\matrix{D \cr \sqrt 2 d \cr 1}\right]=
\left[\matrix{D^2+i(uD-u) \cr * \cr D-1+iu}\right].$$
(The starred entry is irrelevant to us.)
Multiplying through by $D-1-iu$ we see that
$I_{C_1}(V_1)$ is a scalar multiple of the matrix
\begin{equation}
\label{schwartz1}
\left[\matrix{P(u,D)+iu(1-2D) \cr * \cr u^2+(D-1)^2}\right].
\end{equation}
Here $P(u,d)$ is a real polynomial in $u$ and $D$ whose
form is not important to us.
If $(z,t)$ is the point represented by $I_{C_1}(V_1)$,
then 
$$|t|=\frac{u|1-2D|}{u^2+(D-1)^2}<
\frac{|1-2D|}{u} \leq 1/u \leq 1/4.$$
But the point $(z,t)$ is one of the endpoints
of $Q_2=I_1(Q_0)$.  
The same argment works for the other endpoint
of $Q_2$.
\endproof

Here are three geometric facts:
\begin{enumerate}
\item
$Q_2$ is centered at $0$ in $\cal H$ and
extends upwards at most $1/4$ in either
direction.
Since $E_2$ is
contained in a contact plane of slope
at least $7/3$, the circle
$E_2$ has slope $7/3$ at the origin
by symmetry.  Since $Q_2$ rises up only
$1/4$ in either direction away from the
origin, we see easily that every point
of $Q_2$ has slope greater than $1$.
\item Given that every point of $Q_2$
has slope greater than $1$, we see that
the projection $\pi(Q_2)$ is contained in
the dist $\Delta_{1/4}$ of radius $1/4$ about the origin.
\item Any harmonic $\R$--arc which intersects
$Q_2$ has slope less than $1/2$ at the
intersection point, because this $\R$--arc
is integral to the contact structure
and the intersection point projects
inside $\Delta_{1/4}$.
\end{enumerate}
These properties together imply that
$\Psi(Q_2)$ does not contain any
extrema of the height function, even
at the endpoints.  This establishes Item 2 of
Technical Lemma II.

\medskip
{\bf Remark}\qua A more intuitive way to see the
same result is that $E_2$ is extremely tall
in comparison to $Q_2$, and the fibers of
$\Psi$ are fairly straight near $Q_2$ (because
they have large diameter) and hence there is
no way the height of $\Psi|_{E_2}$ takes on its
extrema on $Q_2$.

\section{Technical Lemma III}\label{sec7}
\setcounter{figure}{0}

\subsection{Estimating $Q_{21}$}\label{subsec7.1}

We continue our work from the previous chapter.
In this section, we continue to normalize
using $\B$.  However, our main estimate
is independent of the normalization.

Recall that $\Psi$ is the loxodromic
elevation map, the map of interest to us.
Given a subset $X$, which is disjoint
from $E_0$, let
$2\|X\|$ denote the ``vertical diameter''
of the set $\Psi(X)$.  By this we mean
that $2\|X\|$ denotes the maximum difference
in heights between two points
of the form $\Psi(p_1)$ and $\Psi(p_2)$,
where $p_1,p_2 \in X$.  We call
$\|X\|$ the {\it vertical $\Psi$--radius\/}
of $X$.  This quantity is not quite
canonical; it depends on us choosing a
scaling factor for the image of $\Psi$.
However, we shall always be interested
in quantities of the form
$\|X\|/\|Y\|$, and this ratio is independent
of the way we scale $\Psi$.
The main goal of this section is to prove

\begin{lemma}
\label{hard1}
\begin{equation}
\frac{\|Q_{12}\|}{\|C_1\|}<\frac{1}{15}.
\end{equation}
for all $s \in [\underline s,\overline s]$.
\end{lemma}

So far we have estimated quantities in $\cal H$ and
we need to translate the information we have
gathered into terms related to $\Psi$.  Here is
an outline of how we will do this.
For $X \subset {\cal H\/}$ let $\|X\|'$ denote half
the width of the smallest vertical slab which
contains $X$.  At the parameter $s=\overline s$ we
have $\|X\|'=\|X\|$.   We want to compare
$\|X\|'$ and $\|X\|$ in general.  For this
purpose we let $[X] \subset E_0-Q_0$ denote
the set of all (lower) endpoints of
harmonic $\R$--arcs which contain points of
$X$.  These $\R$--arcs are meant to be
harmonic with respect to $(E_0,Q_0)$. They
are the arcs used in the definition of $\Psi$.
Here are the estimates we will prove:
\begin{enumerate}
\item $\|Q_{12}\|'/\|C_1\|'<1/17$.
\item $\|[Q_{12}]\|'<\|Q_{12}\|'$ and
$\|[C_1]\|'>(15/17)\|C_1\|'$.
\item $\|Q_{12}\|/\|C_1\|<\|[Q_{12}]\|'/\|[C_1]\|'$.
\end{enumerate}
These estimates combine, in a straightforward way,
to establish Lemma \ref{hard1}.

\subsubsection{First estimate}\label{subsubsec7.1.1}

We have $Q_{21}=I_1(Q_{01})$.  
Here $Q_{01} \subset Q_0$
consists of points $x$ such that $x$ is the endpoint
of an $\R$--arc which is harmonic with respect
to $(E_0,Q_0)$ and which intersects $C_1$.
The endpoints of $Q_{01}$ are represented
by vectors of the form:
\begin{equation}
V_1'=\left[\matrix{iu' \cr 0 \cr 1}\right] \hskip 30 pt
V_2'=\left[\matrix{-iu' \cr 0 \cr 1}\right]
\end{equation}
The two endpoints 
$t_1$ and $t_2$
of any $\R$--arc, harmonic with respect to
$(E_0,Q_0)$, satisfy $$t_1t_2=u^2>17,$$ and
also $|t_1|<\|C_1\|'$.  Hence
$u'>17\|C_1\|'$.
The same argument as in Lemma \ref{extend1}
now shows that both endpoints of $Q_{21}$
are at most $$1/(17/\|C_1\|')=\|C_1\|'/17$$ away from
$\C \times \{0\}$.  Hence $\|Q_{12}\|'<\|C_1\|'/17$.

\subsubsection{The second estimate}\label{subsubsec7.1.2}

Our harmonic arcs vary monotonically in
height from their lower to their upper
endpoints.  Hence
$\|[Q_{12}]\|' \leq \|Q_{12}'\|$.
This estimate is true for any set $X$,
actually.  The other estimate is the
interesting one.  

\begin{figure}[ht!]\anchor{fig7.1}
\cl{\small
\psfrag {h}{$\eta(\alpha)$}
\psfrag {L}{$\Lambda$}
\psfrag {h1}{$h_1$}
\psfrag {h2}{$h_2$}
\psfrag {x1}{$(1-x_1)^2$}
\psfraga <15pt,3pt> {x2}{$(1+x_1)^2$} 
\includegraphics{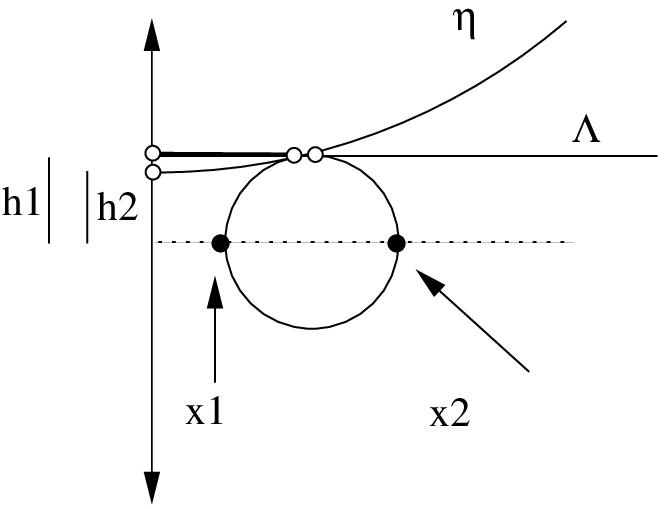}
}
\caption{}\label{fig7.1}
\end{figure}

Let $\eta$ be the map described in subsection~\ref{subsec3.9}.
To do our analysis efficiently we post-compose
$\eta$ with a Mobius transformation so that
$\eta$ is the identity on $E_0$ and
$\eta({\cal H\/})$ is the left half plane.
(The HYP1 coordinate system in our applet
is precisely this map.)
 $\eta$ conjugates the Heisenberg
automorphism $(z,t) \to (rz,r^2t)$ to the
hyperbolic isometry $z \to rz$.  Hence,
the restriction of $\eta$ to the 
$\R$--circle
$R_1=R \times \{0\}$ is the map $r \to r^2$.

let $(c,0)$ denote the center of $C_1$.
Here $c \in (0,1/3]$ as in subsection \ref{subsec6.3}.
Since $C_1$ projects to a circle of
radius $1$ we have $\|C_1\|'=2c$.
Now, $C_1$ intersects $R_1$ in the two points
$c-1$ and $c+1$.  Hence
$\eta(C_1)$ is the circle of radius
$2c$ which intersects the positive
real axis in the points
$(1-c)^2$ and $(1+c)^2$.
\figref{fig7.1} shows the picture.

Let $\alpha$ be a harmonic arc
which contains a point of $C_1$.
Then $\eta(\alpha)$ is a semicircle
which connects two points
of the $y$--axis which are
harmonic with respect to the
endpoints of $Q_0$.  The upper
endpoint of $Q_0$ is at least $4$
away from the origin.
Also, the lower endpoint 
of $\eta(\alpha)$ is at most $2c$
away from $0$.  Hence the other
endpoint of $\eta(\alpha)$ is
at least $4^2/(2c)>8/c$
away from the origin.  
Hence $\eta(\alpha)$ has
radius at least \begin{equation}
\label{rad}
R(c)=
\frac{4}{c}-c
\end{equation}
Note that $R(c)>11$ for $c \in [0,1/3]$.

We want to compare the point $h_2$ where
$\eta(\alpha)$ intersects the $y$--axis
with the point $h_1$ where the horizontal line $\Lambda$
tangent to $\eta(C_1)$ intersects the $y$--axis.
These two points are shown in \figref{fig7.1}.
Our goal is to show that $h_2/h_1>15/17$.
Let $\lambda$ be the length of the
portion of $\Lambda$
contained inside $\eta(C_1)$.
This segment is drawn thickly in
\figref{fig7.1}.
Using a familiar fact from high school
geometry we have
\begin{equation}
\label{hard3}
(h_1-h_2)(2R(c)-(h_1-h_2)) \leq \lambda^2<(c+(1-c)^2)^2.
\end{equation}
The right hand side of the equation comes from
the fact that the intersection $\Lambda \cap \eta(C_1)$
occurs closer to the $y$--axis than does the
center of $\eta(C_1)$.  

To simplify this equation note that $2R(c)>22$ and
certainly $h_1-h_2<1$. Hence
$$(2R(c)-(h_1-h_2))>\frac{21}{22}R(c).$$
Therefore (using the fact that $h_1=2c$) we have
\begin{equation}
\label{strong1}
1-\frac{h_2}{h_1}=
\frac{h_1-h_2}{h_1}=
\frac{h_1-h_2}{2c}<
\frac{22}{21} \times \frac{(100/81)}{4c R(c)} \leq
\frac{1}{12}<\frac{2}{17}.
\end{equation}
In particular $h_2/h_1>15/17$, as claimed.
(We will use the stronger estimate in
Equation \ref{strong1} below.)

\subsubsection{The third estimate}\label{subsubsec7.1.3}

At this point we need to look at the loxodromic
elevation map $\Psi$ geometrically.
We need to quantify the
difference in heights between the two points
$\Psi(p_1)$ and $\Psi(p_2)$.   Let
$\alpha_j$ be the harmonic arc
which contains $p_j$.  Let $x_j \in E_0-Q_0$
be the relevant endpoint of $\alpha_j$.
Let $u_1=u$ and $u_2=-u$ be the two endpoints
of $Q_0$.  Then, up to a constant,
the difference in heights of
$\Psi(p_1)$ and $\Psi(p_2)$ is:
\begin{equation}
\label{logx}
\log \frac{(u_1-x_2)(u_2-x_1)}{(u_1-x_1)(u_2-x_2)}
\end{equation}
This quantity is none other than the
hyperbolic radius of the segment 
$[x_1,x_2]$, when the segment $[-u,u]=
E_0-Q_0$ is identified with the
hyperbolic line.
This interpretation comes from the fact that
the quantity in Equation \ref{logx} is
invariant under Heisenberg automorphisms
of $(E_0,Q_0)$, and $\Psi$ conjugates such
maps to vertical translations.

After rotating the picture sideways (so as to
work with the $x$--axis rather than the $y$--axis)
our third estimate reduces to the following
situation.  We have points $u$ and $-u$ as
reference points.  We have positive 
$0<v<w<u$ and we consider the intervals
$I_v=[-v.v]$ and $I_w=[-w,w]$.  We let
$\|I_v\|_u$ denote the hyperbolic
radius of the segment $I_v$ when the
interval $[u,-u]$ is identified with
the hyperbolic line.  Likewise
we define $\|I_w\|_u$.

To make our third estimate, we have in mind that
$$[Q_{12}]=I_v; \hskip 50 pt [C_1]=I_w.$$
$$\|[Q_{12}]\|'=\|I_v\|_{\infty}; \hskip 30 pt
\|[C_1]\|'=\|I_w\|_{\infty};\leqno{\hbox{We have}}$$
$$\|Q_{12}\|=\|I_v\|_{u}; \hskip 30 pt
\|C_1\|=\|I_w\|_{u}.$$
Here our specific choice of $u$ is as
in Equation \ref{polar1}.
Thus, to establish our third estimate, it
suffices to show that:
$$\frac{v}{w}>\frac{\|I_v\|_u}{\|I_w\|_u}=
\frac{\log(\frac{(u+v)^2}{(u-v)^2})}
{\log(\frac{(u+w)^2}{(u-w)^2})}=
\frac{\log(u+v)-\log(u-v)}{\log(u+w)-\log(u-w)}$$
This is an easy exercise in calculus.
This completes our proof of the third estimate.

\subsection{Estimating the gap}\label{subsec7.2}

At this point we have done the hardest part of
the estimate, which involved controlling
$\Psi(Q_{21})$.  Given Lemma \ref{hard1} we
are back to a problem involving the
two curves $\Psi(C_1)$ and $\Psi(C_2)$.
It is difficult to draw these curves well,
so we will consider $\eta(C_1)$ and $\eta(C_2)$,
where $\eta$ is the map considered in
the previous section.  Given the work done
in the previous section, we will see
readily how to translate back and forth
between $\eta$ coordinates and
$\Psi$ coordinates.

\begin{figure}[ht!]\anchor{fig7.2}
\cl{\small
\psfrag {g1}{$\gamma_1$}
\psfrag {g2}{$\gamma_2$}
\psfrag {g3}{$\gamma_3$}
\psfrag {g4}{$\gamma_4$}
\psfrag {g5}{$\gamma_5$}
\psfrag {h1}{$\eta(C_1)$}
\psfrag {h2}{$\eta(C_2)$}
\psfraga <-4pt,4pt> {hQ}{$\eta(Q_{21})$}
\includegraphics[width=2.7in]{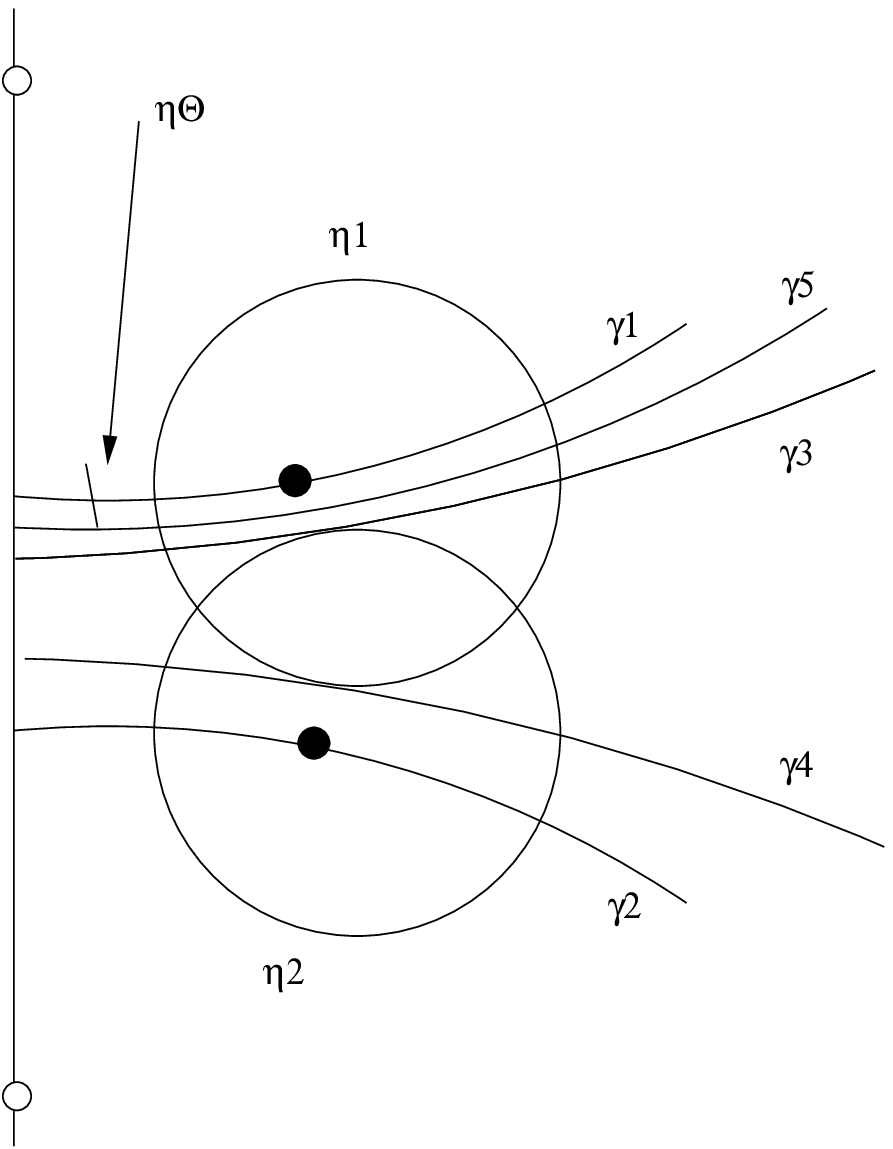}
}
\caption{}\label{fig7.2}
\end{figure}

\figref{fig7.2} shows a schematic (and fairly accurate)
picture of $\eta(C_1)$ and $\eta(C_2)$. 
Here $\eta$ is the map we used in the previous
section.

Regarding this picture:
\begin{itemize}
\item The two white dots on the $y$--axis are
the endpoints of $\eta(Q_0)$.
\item The two black dots are the hyperbolic
centers of $\eta(C_1)$ and $\eta(C_2)$.
\item All the geodesics drawn are orthogonal
to the geodesic joining the endpoints of
$\eta(Q_0)$.  Indeed, these geodesics are
all certain images of harmonic $\R$--arcs.
\item For $j=1,2$, we have
$\gamma_j=\eta(R_j)$, where $R_j$ is
the $\R$--axis of $C_j$.
\item $\|C_1\|$ is the hyperbolic distance
between $\gamma_1$ and $\gamma_2$ and
equally well the hyperbolic distance
between $\gamma_2$ and $\gamma_4$.
\item The geodesic $\gamma_5$ contains
the lower endpoint of $\eta(Q_{21})$.
\end{itemize}
To finish the Technical Lemma III we just
need to show that $\eta(Q_{21})$ lies above
the geodesic $\gamma_3$.  Let $d_{ij}$ denote
the hyperbolic distance between $\gamma_i$ and
$\gamma_j$.  Lemma \ref{hard1}, interpreted
in terms of $\eta$, says that:
\begin{equation}
\frac{d_{14}}{d_{15}}<\frac{1}{15}
\end{equation}
To prove item 1 of Technical Lemma III it suffices
to establish the inequality:
\begin{equation}
\label{crunch}
\frac{d_{13}}{d_{14}}>\frac{1}{15}
\end{equation}
Using a bit of algebra, we see that
Equation \ref{crunch} is equivalent to the
more symmetric:
\begin{equation}
\label{crunch2}
\frac{d_{34}}{d_{12}}<\frac{7}{8}
\end{equation}

{\bf Remark}\qua
Our normalization here is
slightly different than in the previous section.
We are now normalizing so that reflection in
the $x$ axis interchanges $\eta(C_1)$ and
$\eta(C_2)$.  The reason we have changed
normalizations is that previously we were
just concentrating on $C_1$ and now we
need to treat both $C_1$ and $C_2$.
Note that Lemma \ref{hard1} is independent
of normalization, as we mentioned above. 
\medskip

Let $d'_{12}$ denote the Euclidean distance
between the centers of $\eta(C_1)$ and
$\eta(C_2)$ and let $d'_{34}$ denote the
Euclidean distance between the point of
maximum height on $\eta(C_2)$ and
the point of minumum height of
$\eta(C_1)$.  

\begin{lemma}
Equation \ref{crunch2} 
is true provided that
\begin{equation}
\label{crunch3}
\frac{d'_{34}}{d'_{12}}<\frac{4}{5}.
\end{equation}
\end{lemma}

\proof
Our proof here is very much like what we did in the
previous section, for Estimates 2 and 3.  Estimate 3
works exactly the same way, and Estimate 3 is replaced
by the statements that
\begin{equation}
\label{hard5}
d_{34}<d'_{34}; \hskip 40 pt
d_{12}>\frac{11}{12} d'_{12}.
\end{equation}
Thus, if Equation \ref{crunch3} holds then:
$$\frac{d_{34}}{d_{12}} <
\frac{12}{11} \times \frac{4}{5}<\frac{7}{8}$$
The first part of Equation \ref{hard5}  is immediate,
as above.
The second part of Equation \ref{hard5} is similar
to what we have already done, though we have to
think about our new normalization.
In the old normalization we knew that
$\eta(C_1)$ was centered at the real
axis, and had a center which was at
most $10/9$ from the $y$--axis.
(In the old normalization $\eta(C_2)$ was some other
circle below $\eta(C_1)$.)
Our new normalization is obtained by the
old one by applying a hyperbolic isometry $I$
translating along the geodesic joining
the endpoints of $\eta(Q_0)$.   The hyperbolic
isometry is chosen so that $\eta(C_1)$ and
$\eta(C_2)$ are moved into symmetric position.

Let $h_1$ denote the Euclidean distance
between the centers of $\eta(C_1)$ and
$\eta(C_1)$.
Let $h_2$ denote the Euclidean length of the
segment on the $y$--axis between
$\gamma_1$ and $\gamma_2$.
We want to show that $h_2/h_1<11/12$.
This is the same as showing that:
$$1-\frac{h_2}{h_1}<1/12$$
The argument given in subsection \ref{subsubsec7.1.2}
goes through, once a few changes are
noted.

By construction,
the $x$--coordinate of the
hyperbolic center of
the new $\eta(C_1)$ is less than
the $x$--coordinate of the hyperbolic
center of the old $\eta(C_1)$,
which is in turn less than the
$x$ coordinate of the Euclidean
center of the old $\eta(C_1)$.
We conclude that the Euclidean
distance from the hyperbolic
center of the new $\eta(C_1)$
to the $y$--axis
is at most $10/9$.
(This hyperbolic center is the intersection
point of $\gamma_1=\eta(\alpha_1)$ and
$\Lambda$, the curves from subsection \ref{subsubsec7.1.2}.)
The radius of the circle containing
$\gamma_1$ is still at least $R(c)$.
Having made these definitions, we
see that the rest of the argument
in subsection \ref{subsubsec7.1.2} is the same.
\endproof

By the previous result,
Equation \ref{crunch3} implies item 1
of Technical Lemma III.  Since $4/5<1$,
Equation \ref{crunch3} 
implies item 2 of Technical Lemma III.
Thus, to finish the proof of Technical Lemma III
we just have to verify Equation \ref{crunch3}.

\subsection{Back to Heisenberg space}\label{subsec7.3}

\begin{figure}[ht!]\anchor{fig7.3}
\cl{\small
\psfraga <-2pt,0pt> {T}{$T$}
\psfrag {R}{$R$}
\psfraga <-2pt,0pt> {S}{$S$}
\psfrag {0}{$0$}
\psfraga <-2pt,0pt> {1}{$1$}
\psfrag {q}{$\theta$}
\psfrag {R0}{$\rho_0$}
\psfrag {C1}{$C_1$}
\psfrag {C2}{$C_2$}
\psfrag {c1}{$c_1$}
\psfrag {axis1}{axis 1}
\includegraphics[width=.9\hsize]{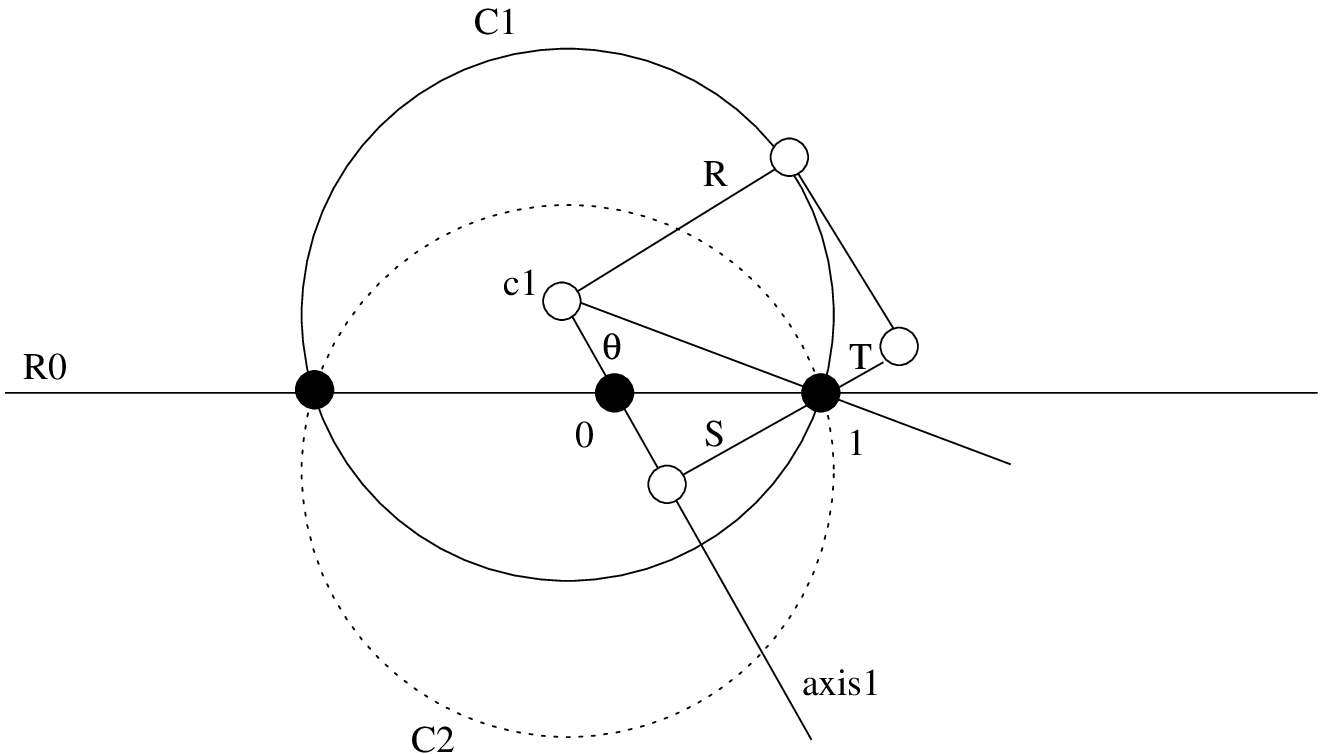}
}
\caption{}\label{fig7.3}
\end{figure}

\figref{fig7.3} shows the picture of 
$\pi(C_1)$ and $\pi(C_2)$ in $\C$.  
Referring to \figref{fig7.3}, we will show that
Technical Lemma III is true provided
that $\sin(\theta)>5/9$.

Let $c_1=z+it$ be the center of $C_1$.
The quantity $d_{12}'$ is just the vertical
distance between centers of the two circles.
The quantity $d_{34}'$ is the vertical distance
between the minimum height point on $C_1$ and
the maximum height point on $C_2$.

The outer edge of the rectangle in \figref{fig7.3}
is tangent to
$\pi(C_1)$ and parallel to axis 1, the
$\R$--circle through the center of $C_1$ which intersects
$C_1$ in two other points.  The number $R$ is
the radius of $C_1$.

The vertical distance between the two centers
is twice the distance from the center of
$C_1$ to $\C \times \{0\}$.  This latter distance
is $2|z|S$ because $S$ runs
perpendicular to axis 1
and joins a point on $C_1$ to
the point $(1,0)$. 
The vertical distance between the max of $C_2$
and the min of $C_1$ is twice the vertical
distance from $\C \times \{0\}$ to the min of
$C_1$.  This distance is 
$2|z|T$. Therefore
$$\frac{d'_{34}}{d'_{12}}=\frac{T}{S}=
\frac{R-R\sin(\theta)}{R\sin(\theta)}=\frac{1}{\sin(\theta)}-1.$$
If $\sin(\theta)>5/9$ then the quantity on the right
is less than $4/5$.

\subsection{The end of the proof}\label{subsec7.4}

We will suppose that $\sin(\theta)<5/9$ 
and derive a contradiction.  

\subsubsection{Siegel model computation}\label{subsubsec7.4.1}  
The center of
$C_1$ is the point $c=z+it$.
In Section~\ref{sec5} we
saw that $C_1$ has aspect less than $10$.
Hence 
\begin{equation}
\frac{R}{|z|}<\sqrt{10}.
\end{equation}
In the Siegel model,
the vectors representing $0,\infty,c$ are
respectively:
\begin{equation}
\label{Vs}
V_1=\left[\matrix{0 \cr 0 \cr 1}\right] \hskip 30 pt
V_1=\left[\matrix{1 \cr 0 \cr 0}\right] \hskip 30 pt
V_1=\left[\matrix{-|z|^2+it \cr \sqrt 2 z \cr 1}\right]
\end{equation}
Letting $\delta$ be the invariant from Equation \ref{ai}
we have
\begin{equation}
\label{schwartz3}
\label{badd}
\delta(V_1,V_2,V_3)=\frac{t}{|z|^2}=
\frac{2|z|R\sin(\theta)}{|z|^2}=
\frac{R}{|z|} \times 2 \sin(\theta)<
\sqrt{10} \times \frac{10}{9}<3.52.
\end{equation}

\subsubsection{Ball model computation}\label{subsubsec7.4.2}

In the ball model the corresponding vectors are:
\begin{equation}
U_1=\left[\matrix{b\cr\overline b\cr1}\right] \hskip 30 pt
U_2=\left[\matrix{a\cr\overline a\cr1}\right]\hskip 30 pt
U_3=\left[\matrix{-1 & 0 & 0 \cr
               0 & 3 & -4\overline \beta \cr
               0 & 4\beta & -3} \right] 
\left[\matrix{a\cr\overline a\cr1}\right]
\end{equation}
Note that $\langle U_1,U_2 \rangle \in \R$.
Also $\langle U_2,U_3 \rangle \in \R$ because
these two vectors are exchanged by $I_1$.
Therefore
$\delta(V_1,V_2,V_3)=|{\rm Im\/}(\tau)|/|{\rm Re\/}(\tau)|$,
where
\begin{equation}
\tau=\langle U_1,U_3 \rangle=
3+3(a\overline b+b\overline a)-
4 a \overline b - 4 b \overline \beta - 4 \beta \overline a.
\end{equation}
Using Equations \ref{parker199}, \ref{massive}, 
\ref{massive3}, and \ref{finale} we get:
\begin{equation}
\label{great3}
\delta(V_1,V_2,V_3)=\sqrt{\frac{(x-1)(2x-3)}{x-3}}
\end{equation}
A bit of calculus now shows that
\begin{equation}
\delta(V_1,V_2,V_3)>3.83
\end{equation}
in our range of parameters.  This contradicts
Equation \ref{badd} and we are done.


\begin{thebibliography}

\bibitem{AGG} \textbf{S Anan'in}, \textbf{H Grossi}, \textbf{N Gusevskii},
{\it Complex Hyperbolic Structures on Disc Bundles over
Surfaces I\/}, preprint (2003)

\bibitem{E}
\textbf{D\,B\,A Epstein}, \emph{Complex hyperbolic geometry}, from:
  ``Analytical and geometric aspects of hyperbolic space (Coventry/Durham,
  1984)'', London Math. Soc. Lecture Note Ser. 111, Cambridge Univ. Press,
  Cambridge (1987)  93--111 \MR{903851}

\bibitem{FP}
\textbf{E Falbel}, \textbf{J\,R Parker}, \emph{The moduli space of the modular
  group in complex hyperbolic geometry}, Invent. Math. 152 (2003) 57--88
  \MR{1965360}

\bibitem{G}
\textbf{W\,M Goldman}, \emph{Complex hyperbolic geometry}, Oxford Mathematical
  Monographs, The Clarendon Press Oxford University Press, New York (1999)
  \MR{1695450}

\bibitem{GP}
\textbf{W\,M Goldman}, \textbf{J\,R Parker}, \emph{Complex hyperbolic ideal
  triangle groups}, J. Reine Angew. Math. 425 (1992) 71--86 \MR{1151314}

\bibitem{S0}
\textbf{R\,E Schwartz}, \emph{Ideal triangle groups, dented tori, and numerical
  analysis}, Ann. of Math. (2) 153 (2001) 533--598 \MR{1836282}

\bibitem{S1}
\textbf{R\,E Schwartz}, \emph{Degenerating the complex hyperbolic ideal
  triangle groups}, Acta Math. 186 (2001) 105--154 \MR{1828374}

\bibitem{S2} \textbf{R\,E Schwartz}, {\it Spherical CR Geometry and
Dehn Surgery\/},  preprint of Research Monograph (2004)

\bibitem{S3} \textbf{R\,E Schwartz}, {\it An Interactive Proof
of the G-P Conjecture\/},  
Java applet (2004)  
\url{http://www.math.brown.edu/~res/Java/App45/test1.html} 

\bibitem{W}
\textbf{S Wolfram}, \emph{The {M}athematica book}, Wolfram
  Media, Inc. Champaign, IL (1999) \MR{1721106}


\end{thebibliography}
\end{document}